\documentclass[11pt,twoside,leqno]{aomamlt2e} 
  \pageno{423}
\received{June 23, 2003}

\begingroup \makeatletter
\@ifundefined{theorem@style}{\input{theorem.sty}}{}
             [\FMithmInfo]
\gdef\th@plain{\normalfont\itshape
  \def\@begintheorem##1##2{%
        \item[\hskip\labelsep \theorem@headerfont \hskip22pt ##2\ ##1.]}%
\def\@opargbegintheorem##1##2##3{%
   \item[\hskip\labelsep \theorem@headerfont\hskip22pt  ##2\ ##1 {\rm (##3)}.\ ]}}
\endgroup

   \newtheorem{thm}[equation]{Theorem}
\newtheorem{cor}[equation]{Corollary}
\newtheorem{lem}[equation]{Lemma}
\newtheorem{prop}[equation]{Proposition}
\theorembodyfont{\upshape}

\newcommand{\lemref}[1]{Lemma~\ref{#1}}
\newcommand{\propref}[1]{Proposition~\ref{#1}}
\newcommand{\corref}[1]{Corollary~\ref{#1}}

\newcommand{\secref}[1]{Section~\ref{#1}}

\DeclareMathOperator{\Hom}{Hom}

\DeclareMathOperator{\sg}{sgn}
\DeclareMathOperator{\md}{mod}
\DeclareMathOperator{\imag}{Im}

\numberwithin{equation}{section}

\newcommand{\gobble}[1]{}
  \newcommand{\rangeref}[2]{%
    \ref{#1}--\afterassignment\gobble\fam 0\ref{#2}%
  }

\renewcommand\a{\alpha}         

\newcommand\g{\gamma}
\renewcommand\d{\delta}
\newcommand\e{\varepsilon}
\renewcommand\l{\lambda}

\newcommand\G{\Gamma}
\newcommand\f{\frac}
\newcommand{\N}{{\mathbb{N}}}
\newcommand{\Z}{{\mathbb{Z}}}
\newcommand{\R}{{\mathbb{R}}}

\newcommand{\C}{{\mathbb{C}}}
\newcommand{\A}{{\mathbb{A}}}
\newcommand{\Q}{{\mathbb{Q}}}

\newcommand{\Sch}{{\mathcal{S}}}

\renewcommand\Re{\mbox{Re~}}

\newcommand{\ttwo}[4]{\left(\begin{array}{cc}
{#1} & {#2} \\ {#3} & {#4} \end{array} \right)}
\newcommand{\ttt}[9]{\left(\begin{array}{ccc} {#1} & {#2} & {#3} \\
    {#4} & {#5} & {#6} \\   {#7} & {#8} & {#9} \end{array}  \right)}

\renewcommand\i{^{-1}}

\renewcommand\({\left(}         
\renewcommand\){\right)}

\begin{document}
\currannalsline{164}{2006} 

\title{Automorphic distributions, $L$-functions,\\ and Voronoi summation
for ${\rm GL}(3)$}

 \acknowledgements{Supported by NSF grant
DMS-0122799 and an NSF post-doctoral fellowship.

\hglue7pt$^{\textstyle\ast\ast}$Supported in part by NSF grant DMS-0070714.
}
\twoauthors{Stephen D. Miller$^{\textstyle\ast}$}{Wilfried
Schmid\raise2.75pt\hbox{${\textstyle\ast}$}}

 \institution{Rutgers University, Hill Center-Busch Campus, Piscataway, NJ
\\
\email{miller@math.rutgers.edu}\\
\vglue-9pt
Harvard University, Cambridge, MA
\\
\email{schmid@math.harvard.edu}}

 \shorttitle{Automorphic distributions} 
 \shortname{Stephen D. Miller and Wilfried Schmid}

\section{Introduction}\label{introduction}

In 1903 Voronoi \cite{Voronoi:1903} postulated the existence of explicit formulas for sums of the form
\begin{equation}
\label{v1904}
{\sum}_{n\geq 1}\ a_n\,f(n) \,,
\end{equation}
for any ``arithmetically interesting" sequence of coefficients $(a_n)_{n \geq 1}$ and every $f$ in a large class of test functions, including characteristic functions of bounded intervals. He actually established such a formula when $a_n = d(n)$ is the number of positive divisors of $n$ \cite{Voronoi:1904}. He also asserted a formula for
\begin{equation}
\label{vba}
a_n \ = \ \#\{(a,b) \in \Z^2 \mid Q(a,b)=n\}\,,
\end{equation}
where $Q$ denotes a positive definite integral quadratic form \cite{Voronoi:1905}; Sierpi\'nski \cite{Sierpinski:1906} and Hardy \cite{Ha} later proved the formula rigorously. As Voronoi pointed out, this formula implies the bound
\begin{equation}
\label{vbd} \left|\,\#\{\,(a,b)\, \in \,\Z^2 \mid a^2 + b^2 \le x\,\}\, -\, \pi x \, \right| \ = \ O(x^{1/3})
\end{equation}
for the error term in Gauss' classical circle problem, improving greatly on Gauss' own bound $O(x^{1/2})$.  Though Voronoi originally deduced his formulas from Poisson summation in $\R^2$, applied to appropriately chosen test functions, one nowadays views his formulas as identities involving the Fourier coefficients of modular forms on ${\rm GL}(2)$, i.e., modular forms on the complex upper half plane. A discussion of the Voronoi summation formula and its history can be found in our expository paper \cite{MilSch:2003a}.

The main result of this paper is a generalization of the Voronoi summation formula to 
${\rm GL}(3,\Z)$-automorphic representations of ${\rm GL}(3,\R)$. Our technique is quite general; we plan to extend
the formula to the case of ${\rm GL}(n,\Q)\backslash {\rm GL}(n,\A)$ in the future. The arguments make heavy use of
representation theory. To illustrate the main idea, we begin by deriving the well-known generalization of
the Voronoi summation formula to coefficients of modular forms on ${\rm GL}(2)$, stated below in
(\ref{gl2voronoi})--(\ref{gl2voronoi4}). This formula is actually due to Wilton -- see \cite{Huxley:1993} --
and is not among the formulas predicted by Voronoi.  However, because it is quite similar in style one
commonly refers to it as a Voronoi summation formula. We shall follow this tradition and regard our ${\rm GL}(3)$
formula as an instance of Voronoi summation as well. The ${\rm GL}(2)$ formula is typically derived from
modular forms via Dirichlet series and Mellin inversion; see, for example, \cite{ConIwan:2002},
\cite{KowMichVan:2002}.  We shall describe the connection with Dirichlet series later on in this
introduction. Since we want to exhibit the analytic aspects of the argument, we concentrate on the case of
modular forms invariant under $\G={\rm SL}(2,\Z)$. The changes necessary to treat the case of a congruence
subgroup can easily be adapted from \cite{ConIwan:2002}, \cite{KowMichVan:2002}, for example.

We consider a cuspidal, ${\rm SL}(2,\Z)$-automorphic form $\Phi$ on the upper half plane $H = \{ z \in \C \mid \imag z > 0\}$. This covers two separate possibilities: $\Phi$ can either be a holomorphic cusp form, of -- necessarily even -- weight $k$,
\begin{equation}
\label{holomexpan}
\Phi(z)\ =\ {\sum}_{n=1}^\infty\, a_n\, n^{(k-1)/2}\,e(nz),\qquad (\,e(z)\ =_{\rm def}e^{2\pi i z}\,)\,,
\end{equation}
or a cuspidal Maass form -- i.e., $\Phi \in C^\infty(H)$, $\, y^2\(\f{\partial^2}{\partial x^2}+\f{\partial^2}{\partial y^2}\) \Phi = -\l \, \Phi\,$ with $\l =\f14-\nu^2$, $\,\nu \in i\R\,$, and
\begin{equation}
\label{maassexpan} \Phi(x+iy)\ = \ {\sum}_{n\neq 0}\, a_n \, \sqrt{y} \, K_{\nu}(2\pi|n|y)\, e(nx)
\end{equation}
\cite{Maass:1949}. In either situation, $\Phi$ is completely determined by the distribution
\begin{equation}
\label{gl2taudef} \tau(x)\ =\ {\sum}_{n\neq 0}\, a_n\, |n|^{-\nu}\, e(nx)\,,
\end{equation}
with the understanding that in the holomorphic case we set both $a_n = 0$ for $n<0$ and $\nu = -\f{k-1}{2}$. One can also describe $\,\tau$ as a limit in the distri\-bution topology: $\,\tau(x)=\lim_{y\to 0^+}\Phi(x+iy)\,$ when $\,\Phi$ is a holomorphic cusp form; the analogous formula for Maass forms is slightly more complicated \cite{Schmid:2000}. As a consequence of these limit formulas, $\tau$ inherits automorphy from $\Phi$,
\begin{equation}
\label{tauaut} \tau(x)\ =\ |cx+d|^{2\nu-1}\,\tau\({\textstyle\f{ax+b}{cx+d}}\)\,,\qquad \text{for any}\ \ \ttwo abcd \in {\rm SL}(2,\Z).
\end{equation}
This is the reason for calling $\tau$ the automorphic distribution attached to $\Phi$. The regularity properties of automorphic distributions for ${\rm SL}(2,\R)$ have been investigated in \cite{BernRez:1999}, \cite{Lewis:1978}, \cite{Schmid:2000}, but these properties are not important for the argument we are about to sketch.

If $c\neq 0$ in (\ref{tauaut}), we can substitute $x - d/c$ for $x$, which results in the equivalent equation
\begin{equation}
\label{tauaut2}
\tau\(x-{\textstyle\f dc}\)\ = \ |cx|^{2\nu-1}\,\tau\({\textstyle\f{a}{c}-\f{1}{c^2x}}\).
\end{equation}
We now integrate both sides of (\ref{tauaut2}) against a test function $g$ in the Schwartz space $\mathcal S(\R)$. On one side we get
\begin{equation}
\label{tauaut3}
\begin{aligned}
\int_\R \tau(x-{\textstyle\f dc})\,g(x)\,dx \ &= \ \int_\R \sum_{n\neq 0}\, a_n\,|n|^{-\nu}\, e(nx - {\textstyle\f{nd}c})\,g(x)\,dx
\\
&=\ \sum_{n\neq 0}\, a_n\,|n|^{-\nu}\, e(-{\textstyle\f{nd}c})\,\, \widehat{g}(-n)\,.
\end{aligned}
\end{equation}
On the other side, arguing formally at first, we find  
\begin{equation} 
\label{tauaut4}\begin{aligned}
{}& \int_{\R} |cx|^{2\nu-1}\,\tau\textstyle{\left(\frac{a}{c}-\frac{1}{c^2x} \right)}\,g(x)\,dx 
\\
{}&\qquad\qquad = 
\ |c|^{2\nu-1}\,\int_{\R} |x|^{2\nu-1} \, \sum_{n\neq 0}\, a_n\,|n|^{-\nu}\, e({\textstyle\frac{na}c}
-{\textstyle\frac{n}{c^2 x}})\, g(x)\, dx 
\\
{}&\qquad\qquad = \ |c|^{2\nu-1}\,\sum_{n\neq 0}\,a_n\,|n|^{-\nu}\, e({\textstyle\frac{na}c})
\int_{\R} |x|^{2\nu-1} \, e(-{\textstyle\frac{n}{c^2 x}})\, g(x)\, dx\,.
\end{aligned}
\end{equation}  
To justify this computation, we must show that (\ref{tauaut2}) 
can be interpreted as an identity of tempered distributions defined on all of $\,\R\,$. A tempered distribution,
we recall, is a continuous linear functional on the Schwartz space $\mathcal S(\R)$, or equivalently, a
derivative of some order of a continuous function having at most polynomial growth. Like any periodic
distribution, $\tau$ is certainly tempered. In fact, since the Fourier series (\ref{gl2taudef}) has no constant
term, $\tau$ can even be expressed as the $n$-th derivative of a {\it bounded} continuous function, for every
sufficiently large $n\in \N$. This fact, coupled with a simple computation, exhibits $|cx|^{2\nu - 1} \, \tau
\({\textstyle\f{a}{c}-\f{1}{c^2x}}\)$ as an $n$-th derivative of a function which is continuous, even at
$x=0$. Consequently this distribution extends naturally from $\,\R^*$ to $\,\R\,$. Using the cuspidality of
$\Phi$, one can show further that the identity (\ref{tauaut2}) 
  holds in the strong sense -- i.e., the extension of $|cx|^{2\nu-1}\,\tau\({\textstyle\f{a}{c}-\f{1}{c^2x}}\)$ which was just described coincides with $\tau(x-{\textstyle\f dc})$ even across the point $x=0$. The fact that $\tau$ is the $n$-th derivative of a bounded continuous function, for all large $n$, can also be used to justify interchanging the order of summation and integration in the second step
 of (\ref{tauaut4}). In any event, the equality (\ref{tauaut4}) is legitimate, and the resulting sum converges
absolutely. For details see the analogous argument in Section \ref{voronsec} for the case of ${\rm GL}(3)$,
as well as \cite{MilSch:2003b}, which discusses the relevant facts from the theory of distributions in some
detail.

Let $f \in \mathcal S(\R)$ be a Schwartz function which vanishes to infinite order at the origin, or more generally, a function such that $|x|^\nu f(x) \in \mathcal S(\R)$. Then $g(x) =\int_\R f(t) |t|^{\nu}\, e(-xt)\, dt$ is also a Schwartz function, and $f(x) = |x|^{-\nu} \, \widehat{g}(-x)$. With this choice of $g$,
(\ref{tauaut2}) to (\ref{tauaut4}) imply
\begin{eqnarray*}
&&{\sum}_{n\neq 0} \, a_n \, e(-nd/c)\,f(n)  
\\
&&\ \   = {\sum}_{n\neq 0}\,a_n \, e(na/c)\, \f{|c|^{2\nu - 1}}{|n|^{\nu}}\!\! \int_{x=-\infty}^\infty \!\!\!\!
|x|^{2\nu-1}e(-\textstyle \f{n}{c^2 x}) \displaystyle \int_{t=-\infty}^\infty  \!\!\!\! f(t) |t|^{\nu}\,e(-x t)\, dt
\,dx
\\
&&\ \  = {\sum}_{n\neq 0}\,a_n \, e(na/c)\, \f{|c|^{2\nu - 1}}{|n|^{\nu}} \!\! \int_{x=-\infty}^\infty
\int_{t=-\infty}^\infty  \!\!\!\! |x|^{-2\nu-1}\, |t|^{\nu}\,f(t)\, e(-\textstyle\f tx - \textstyle\f {nx}{c^2})\, dt \,dx
\\
&&\ \  = {\sum}_{n\neq 0}\,a_n \, e(na/c)\, \f{|c|^{2\nu - 1}}{|n|^{\nu}} \!\! \int_{x=-\infty}^\infty
\int_{t=-\infty}^\infty \!\! |x|^{-\nu}|t|^{\nu}\,f(xt)\,e(-\,t-\textstyle\f{nx}{c^2})\,dt\, dx
\\
&&\ \  = {\sum}_{n\neq 0}\,a_n \, e(na/c)\, \f{|c|}{|n|}  \int_{x=-\infty}^\infty \int_{t=-\infty}^\infty \!\! 
|x|^{-\nu}|t|^{\nu}\,f(\textstyle\f{xtc^2}{n})\,e(-t-x)\,dt\, dx\,.
\end{eqnarray*}
 In this derivation, the integrals with respect to the variable $t$ converge absolutely, since they represent the
Fourier transform of a Schwartz function. The integrals with respect to $x$, on the other hand, converge only
when $\operatorname{Re} \nu > 0$, but have meaning for all $\nu \in \C$  by holomorphic continuation.

So far, we have assumed only that $a, b, c, d$ are the entries of a matrix in ${\rm SL}(2,\Z)$, and $c \neq 0$. We now fix a pair of relatively prime integers $a$, $c$, with $c \neq 0$, and choose a multiplicative inverse $\bar a$ of $a$ modulo $c$:
\begin{equation}
\label{voronoihypothesis} a,\, c,\, \bar a \in \Z\,,\ \ \ (a,c) = 1\,,\ \ \ c\neq 0\,,\ \ \ \bar aa \equiv 1\ \ (\md c)\,.
\end{equation}
Then there exists $b \in \Z$ such that $a\bar a - bc = 1$. Letting $\bar a, b, c, a$ play the roles of $a, b, c, d$ in the preceding derivation, we obtain the {\em Voronoi Summation Formula} for ${\rm GL}(2)$:
\begin{equation}
\label{gl2voronoi}
{\sum}_{n\neq 0}\, a_n\, e(-na/c)\, f(n)\ =\ |c|{\sum}_{n\neq 0}\,\f{a_n}{|n|}\,e(n\bar{a}/c)\,F(n/c^2)\,.
\end{equation}
Here $a_n$ and $\nu$ have the same meaning as in (\ref{holomexpan})--(\ref{gl2taudef}), $f(x) \in
|x|^{-\nu}\mathcal S(\R)$, and
\begin{equation}
\label{gl2voronoi1}
F(t)\ =\
\int_{\R^2}f\!\(\f{x_1x_2}{t}\)|x_1|^{\nu}\,|x_2|^{-\nu}\,e(-x_1-x_2)\, dx_2\, dx_1\,.
\end{equation}
One can show further that this function $F$ vanishes rapidly at infinity, along with all of its derivatives, and has identifiable potential singularities at the origin:
\begin{equation}
\label{gl2voronoi2}
F(x)\ \ \in \ \ \begin{cases}
\ |x|^{1-\nu}\,\Sch(\R)\, + \, |x|^{1+\nu}\,\Sch(\R) &\ \text{if}\ \ \nu\notin \Z
\\
\ |x|^{1-\nu}\,\log |x|\,\Sch(\R)\, + \, |x|^{1+\nu}\,\Sch(\R) &\ \text{if}\ \ \nu\in \Z_{\leq 0}
\end{cases}
\end{equation}
\cite[(6.58)]{MilSch:2003b}; the case $\nu\in\Z_{>0}$ never comes up. The formula (\ref{gl2voronoi1})
for $F$ is meant symbolically, of course: it should be interpreted as a repeated integral, via holomorphic
continuation, as in the derivation. Alternatively and equivalently, $F$ can be described by Mellin inversion,
in terms of the Mellin transform of $f$, as follows. Without loss of generality, we may suppose that $f$ is
either even or odd, say $f(-x) = (-1)^\eta f(x)$ with $\eta \in \{0,1\}$. In this situation,
\begin{equation}
\label{gl2voronoi3}
F(x)  =  \f{\sg(-x)^\eta}{4\pi^2 i} \! \int_{\operatorname{Re} s = \sigma}\!\!\!\!\! \pi^{- 2s} \f {\G(\frac{1 + s + \eta  + \nu } {2})\,\G(\frac{1 + s + \eta  - \nu }{2})}{\G(\frac{-s + \eta + \nu }{2})\,\G(\frac{-s + \eta  - \nu }{2})}\,M_\eta f(-s)|x|^{-s} ds\,,
\end{equation}
where $\,\sigma > | \operatorname{Re} \nu \,| - 1\,$ is arbitrary, and
\begin{equation}
\label{gl2voronoi4}
M_\eta f(s)\ =\ \int_\R f(x)\,\sg(x)^\eta \, |x|^{s-1}\, dx
\end{equation}
denotes the signed Mellin transform. For details see Section \ref{voronsec}, where the ${\rm GL}(3)$
analogues of (\ref{gl2voronoi2}) and (\ref{gl2voronoi3}) are proved.

If one sets $\,c=1\,$ and formally substitutes the characteristic function $\,\chi_{[\varepsilon,x+\varepsilon]}\,$ for $\,f\,$ in (\ref{gl2voronoi}), one obtains an expression for the sum $\sum_{0<n\leq x} a_n\,$;\, formulas of this type were considered especially useful in Voronoi's time. There is an extensive literature on the range of allowable test functions $f$. However, beginning in the 1930s, it became clear that ``harsh" cutoff functions like $\,\chi_{[\varepsilon,x+\varepsilon]}\,$ are no more useful from a technical point of view than the type of test functions we allow in  (\ref{gl2voronoi}).

The Voronoi summation formula for ${\rm GL}(2)$ has become a fundamental analytic tool for a number of deep results in analytic number theory, most notably to the sub-convexity problem for automorphic $L$-functions; see \cite{IwanSar:2000} for a survey, as well 
as 
\cite{DukFrieIwan:2001}, \cite{KowMichVan:2002}, \cite{Sarnak:2001}. In these applications, the
presence of the additive twists in (\ref{gl2voronoi}) -- i.e., the factors $e(-na/c)$ on the left-hand side -- has
been absolutely crucial. These additive twists lead to estimates for sums of modular form coefficients over
arithmetic progressions. They also make it possible to handle sums of coefficients weighted by Kloosterman
sums, such as $\sum_{n\neq 0} a_n f(n) S(n,k;c)$, which appear in the Petersson and Kuznetsov trace
formulas \cite{Goldfeld:1979}, \cite{Sarnak:2001}. In view of the definition of the Kloosterman sum
$S(m,k;c)$, which we recall in the statement of our main theorem below,
\begin{equation}
\begin{aligned}
\label{kuz}
\sum_{n\neq 0}\, a_n\, f(n)\, S(n,k;c)\ &=\ \sum_{d\in(\Z/c\Z)^*} \!\! e(kd/c)\,\,\sum_{n\neq 0} a_n\, f(n)\, e(n \bar{d}/c)
\\
&=\ |c|\,\sum_{n\neq 0}\, \f{a_n}{|n|}\, F(n/c^2)\!\!\!\sum_{d\in(\Z/c\Z)^*}\!\! e((k-n)d/c)\,.
\end{aligned}
\end{equation}
The last sum over $d$ in this equation is a {\it Ramanujan sum}, which can be explicitly evaluated; see, for example,
\cite[p.~55]{Iwaniec:1997}. The resulting expression for $\sum_{n\neq 0} a_n f(n) S(n,k;c)$ can often be
manipulated further.

We should point out another feature of the Voronoi formula that plays an important role in applications. Scaling the argument $x$ of the test function $f$ by a factor $T^{-1}$, $T>0$, has the effect of scaling the argument $t$ of $F$ by the \pagebreak reciprocal factor $T$. Thus, if $f$ approximates the characteristic function of an 
interval, more terms enter the left-hand side of (\ref{gl2voronoi}) in a significant way as the scaling parameter $T$ tends to infinity. At the same time, fewer terms contribute significantly to the right-hand side. This mechanism of lengthening the sum on one side while simultaneously shortening the sum on the other side is known as ``dualizing". It helps detect cancellation in sums like $\sum_{n\leq x} a_n f(n) e(-na/c)$ and has become a fundamental technique in the subject.

We mentioned earlier that our main result is an analogue of the ${\rm GL}(2)$ Voronoi summation formula for cusp forms on ${\rm GL}(3)$:

\begin{thm}\label{mainthm} Suppose that $a_{n,m}$ are the Fourier coefficients of a cuspidal ${\rm GL}(3,\Z)$-automorphic representation of ${\rm GL}(3,\R)$, as
 in {\rm (\ref{cnorm}),} with representation parameters $\l,\d${\rm ,} as in {\rm (\ref{lambda})}. Let $f \in
\Sch(\R)$ be a Schwartz function which vanishes to infinite order at the origin{\rm ,}
 or more generally{\rm ,} a function
on $\,\R - \{0\}\,$ 
such that $\,(\sg x)^{\d_3}|x|^{-\l_3}f(x) \in \Sch(\R)$. Then for $(a,c)=1${\rm ,} $c \neq 0${\rm ,}
$\bar a a \equiv 1 \ (\operatorname{mod} c)$ and $q > 0${\rm ,}
$$
\sum_{n\neq 0}\, a_{q,n}\,e(-na/c)\,f(n)\ =\ \sum_{d|cq}\,\left|\f{c}{d}\right| \, \sum_{n\neq 0}\, \f{a_{n,d}}{|n|}\, S(q\bar{a},n;qc/d) \, F\!\(\f{nd^2}{c^3q}\),
$$
where $\,S(n,m;c)=_{\rm def}\sum_{x\in(\Z/c\Z)^*} e\(\f{nx+m\bar{x}}{c}\)\,$ denotes the Kloosterman sum and{\rm ,}
 in symbolic notation{\rm ,}
$$
F(t)=\int_{\R^3}f\!\(\f{x_1x_2x_3}{t}\)\ {\prod}_{j=1}^3 \((\sg x_j)^{\d_j}\,|x_j|^{-\l_j}\,e(-x_j)\) dx_3\,dx_2\,dx_1\,.
$$
This integral expression for $F$ converges when performed as repeated integral in the indicated order \/{\rm --}\/
 i.e.{\rm ,} with $x_3$ first{\rm ,} then $x_2$, then $x_1$ \/{\rm --}\/ and provided $\,\operatorname{Re}\l_1 >
\operatorname{Re}\l_2 > \operatorname{Re}\l_3$\,{\rm ;} it has meaning for arbitrary values of
$\,\l_1,\,\l_2,\,\l_3\,$ by analytic continuation. If $f(-x) = (-1)^\eta f(x)${\rm ,} with $\eta \in \{0,1\}${\rm ,} one can
alternatively describe $F$ by the identity
$$
F(x) \, = \, \f{\sg (-x)^\eta}{4\pi^{5/2}i}\int_{\operatorname{Re}s = \sigma} \!\!\!\! \pi^{-3s} \! \left( {\prod}_{j=1}^3 i^{\d_j'} \, \pi^{\l_j} \, \f{\G(\f{s+1-\l_j+\d_j'}{2})}{\G(\f{-s+\l_j+\d_j'}{2})} \right)\!  M_\eta f(-s)\,|x|^{-s}\,ds\,;
$$
here $M_\eta f (s)$ denotes the signed Mellin transform {\rm (\ref{gl2voronoi4}),} the $\d_j' \in \{0,1\}$
are characterized by the congruences $\,\d_j' \equiv \d_j + \eta\,$ {\rm ($\operatorname{mod} 2$)\,,} and
$\sigma$ is subject to the condition $\sigma > \operatorname{max}_j(\operatorname{Re}\lambda_j-1)$ but
is  otherwise arbitrary. The function $F$ is smooth except at the origin and decays rapidly at infinity{\rm ,}
along with all its derivatives. At the origin{\rm ,} $F$ has singularities of a very particular type{\rm ,}
 which are described in {\rm
(\ref{vort})} to {\rm (\ref{vorw})} below.
\end{thm}

Only very special types of cusp forms on ${\rm GL}(3,\Z)\backslash {\rm GL}(3,\R)$ have been constructed explicitly; these all come from the Gelbart-Jacquet symmetric square functorial lift of cusp forms on ${\rm SL}(2,\Z)\backslash H$ \cite{GelJac:1978}, though nonlifted forms are known to exist and are far more abundant \cite{Miller:2001}. When specialized to these symmetric square lifts, our main theorem provides a nonlinear summation formula involving the coefficients of modular forms for ${\rm GL}(2)$. The relation between the Fourier coefficients of ${\rm GL}(2)$-modular forms and the coefficients of their symmetric square lifts is worked out in
 \cite[\S5]{MilSch:2003a}.

Our main theorem, specifically the resulting formula for the symmetric squares of ${\rm GL}(2)$-modular forms, has already been applied to a problem originating from partial differential equations and the Berry/Hejhal random wave model in Quantum Chaos. Let $X$ be a compact Riemann surface and $\{\phi_j\}$ an orthonormal basis of eigenfunctions for the Laplace operator on $X$.  A result of Sogge \cite{Sog} bounds the $L^p$-norms of the $\phi_j$ in terms of the corresponding eigenvalues $\l_{\phi_j}$, and these bounds are known to be sharp. However, in the case of $X={\rm SL}(2,\Z)\backslash H$ -- which is noncompact, of course, and not even covered by Sogge's estimate -- analogies and experimental data suggest much
 stronger bounds \cite{HejRack:1992}, \cite{Sarnak:1995}: when the orthonormal basis $\{\phi_j\}$ consists
of Hecke eigenforms, one expects
\begin{equation}
\label{berry}
\|\phi_j\|_p\ \ = \ \ O(\l_{\phi_j}^\e)  \qquad  \(\,\e > 0\,,\ \ 0 < p < \infty\,\)\,.
\end{equation}
Sarnak and Watson \cite{SarWat:2003} have announced (\ref{berry}) for $p=4$, at present under the assumption of the Ramanujan conjecture for Maass forms,
 whereas \cite{Sog} gives the bound $O(\l_{\phi_j}^{1/16})$ in the compact case, for $p=4$. Their
argument uses our Voronoi summation formula, among other ingredients. To put this bound into context, we
should mention that a slight variant of (\ref{berry}) would imply the Lindel\"of Conjecture:
$|\zeta(1/2+it)|=O(1 +|t|^\e)$, for any $\e>0$ \cite{Sarnak:1995}.

There is a close connection between $L$-functions and summation formulas. In the prototypical case of the Riemann $\zeta$-function, the Poisson
summation formula -- which should be regarded as the simplest instance of Voronoi summation -- not only implies, but is equivalent to analytic properties of
the\break $\zeta$-function, in particular its analytic continuation and functional equation. The ideas involved
carry over quite directly to the ${\rm GL}(2)$ Voronoi summation formula (\ref{gl2voronoi}), but encounter
difficulties for ${\rm GL}(3)$.

To clarify the nature of these difficulties, let us briefly revisit the case of ${\rm GL}(2)$.  For simplicity, we suppose $\Phi$ is a holomorphic cusp form, as in (\ref{holomexpan}).  A formal computation shows that the choice of $f(x)=|x|^{-s}$ corresponds to $F(t)=R(s)|t|^s$ in (\ref{gl2voronoi1}), with
\begin{equation}
\label{rs}
R(s)\ = \ i^k\,(2\pi)^{2s-1}\,\f{\G(1-s+\f{k-1}2)}{\G(s+\f{k-1}2)}\ .
\end{equation}
Inserting these choices of $f$ and $F$ into (\ref{gl2voronoi}) results in the equation
\begin{equation}
\label{gl2vwithxs}
{\sum}_{n=1}^\infty \, a_n \, e(-na/c) \, n^{-s} \ = \ R(s)\, |c|^{1-2s}\, {\sum}_{n=1}^\infty \,a_n\,e(n\bar{a}/c)\,n^{s-1}\,,
\end{equation}
which has only symbolic meaning because the regions of convergence of the two series do not intersect. We should remark, however, that the methods of our companion paper \cite{MilSch:2003b} can be used to make this formal argument rigorous.  When $c = 1$, (\ref{gl2vwithxs}) reduces to the functional equation of the standard $L$-function $L(s,\Phi)=\sum_{n=1}^\infty a_n n^{-s}.$  Taking linear combinations over the various $a\in(\Z/c\Z)^*$ for a fixed $c>1$ gives the functional equation for the multiplicatively twisted $L$-function
\begin{equation}
\label{gl2twist}
L(s,\Phi\otimes \chi)\ = \ {\sum}_{n=1}^\infty \,a_n\,\chi(n)\, n^{-s}
\end{equation}
with twist $\chi$, which can be any primitive Dirichlet character mod $c$.

The traditional derivation of (\ref{gl2voronoi}), in \cite{ConIwan:2002}, \cite{KowMichVan:2002} for example, argues in reverse.  It starts with the functional equations for $L(s,\Phi)$ and expresses the left-hand side of the Voronoi summation formula through Mellin inversion,
\begin{equation}
\label{mellinversvor}
\sum_{n=1}^\infty a_n f(n) = \f{1}{2\pi i}\int_{\operatorname{Re}s=\sigma}\!\! L(s,\Phi) Mf(s) ds\,,\ \ \ Mf(s)=\int_0^\infty \!\!\! f(t)t^{s-1}dt\,,
\end{equation}
with $\sigma > 0$. The functional equation for $L(s,\Phi)$ is then used to conclude $\sum_{n=1}^\infty a_n f(n)=\sum_{n=1}^\infty a_n F(n)$, where $MF(s)=r(1-s)Mf(1-s)$. To deal with additive twists, one applies the same argument to the multiplicatively twisted $L$-functions $L(s,\Phi \otimes \chi)$. A combinatorial argument makes it possible to express the additive character $e(-na/c)$ in terms of the multiplicative Dirichlet characters modulo $c$; this not particularly difficult. An analogous step appears already in the classical work of Dirichlet and Hurwitz on the Dirichlet $L$-functions $\sum_{n=1}^\infty \chi(n)n^{-s}$. For ${\rm GL}(3)$, the same reasoning carries over quite easily, but only until this point: the combinatorics of converting multiplicative information to additive information on the right-hand side of the Voronoi formula becomes far more complicated. For one thing, the functional equation for the $L(s,\Phi\otimes \chi)$ only involves the coefficients 
$a_{1,n}$ and $a_{n,1}\,$, whereas the right-hand side of the Voronoi formula involves also the other
coefficients. It is possible to express all the $a_{n,m}$ in terms of the $a_{1,n}$ and $a_{n,1}\,$, but this
requires Hecke identities and is a nonlinear process. The Voronoi formula, on the other hand, is a purely
additive, seemingly nonarithmetic statement about the $a_{n,m}$. In the past, the problem of converting
multiplicative to additive information was the main obstacle to proving a Voronoi summation formula for
${\rm GL}(3)$. Our methods bypass this difficulty entirely by dealing with the automorphic representation
directly, without any input from the Hecke action.

The Voronoi summation formula for ${\rm GL}(3,\Z)$ encodes information about the additively
 twisted $L$-functions $\sum_{n \neq 0} e(na/c) a_{n,q} |n|^{-s}$. It is natural to ask
  if this information is equivalent to the functional equations for the multiplicatively
   twisted $L$-functions $L(s,\Phi\otimes \chi)$. The answer to this question is yes: in
   Section 6 we derive the functional equations for the $L(s,\Phi\otimes \chi)$, and in
   Section~7, we reverse the process by showing that it is possible after all to recover
   the additive information from these multiplicatively twisted functional equations.
   It turns out that our analysis of the boundary distribution -- concretely,
   the ${\rm GL}(3)$ analogues of (\ref{tauaut}) to (\ref{tauaut4})~-- presents the additive
   twists in a form which facilitates conversion to multiplicative twists.
    Section 7 concludes with a proof of the ${\rm GL}(3)$ converse theorem of
    \cite{JacPiatShal:1979}. Though this theorem has been long known,
    of course, our arguments provide the first proof for ${\rm GL}(3)$ that
     can be couched in classical language, i.e., without ad\`eles. To
     explain why this might be of interest, we recall that Jacquet-Langlands
      gave an adelic proof of the converse theorem for ${\rm GL}(2)$ under the
      hypothesis of functional equations for all the multiplicatively twisted
       $L$-functions \cite{JacLan:1970}. However, other arguments demonstrate
        that only a finite number of functional equations are needed
       \cite{Piatetski:1975},   \cite{Weil:1967}. In particular, for the full-level
         subgroup $\G={\rm SL}(2,\Z)$, Hecke proved a converse theorem requiring
         the functional equation merely for the standard $L$-function.\footnote{
         Booker has recently shown \cite{booker} that a single functional
         equation also
         suffices for\break 2-dimensional Galois representations, regardless of the level (see also
\cite{conreyfarm}).} Until
          now it was not clear what the situation for ${\rm GL}(3)$ would be. Our arguments demonstrate that automorphy under $\G={\rm GL}(3,\Z)$ is equivalent to the functional equations for all the twisted $L$-functions. Since the various twisted $L$-functions are generally believed to be analytically independent -- their zeroes are uncorrelated \cite{RudSar:1996}, for example~-- our analysis comes close to ruling out a purely analytic proof using fewer than all the twists.

Our paper proves the Voronoi summation formula only for cuspidal forms, automorphic with respect to the full-level subgroup $\G={\rm GL}(3,\Z)$. It is certainly possible to adapt our arguments to the case of general level $N$, but the notation would become prohibitively complicated. For this reason, we intend to present an adelic version of our arguments in the future, which will also treat the case of ${\rm GL}(n)$, and not just ${\rm GL}(3)$. Extending our formula to noncuspidal automorphic forms would involve some additional technicalities. We are avoiding these because summation formulas for Eisenstein series can be derived from formulas for the smaller group from which the Eisenstein series in question is induced. In fact, the Voronoi summation formula for a particular Eisenstein series on ${\rm GL}(3)$, relating to sums of the triple divisor function $d_3(n)=\#\{x,y,z\in \N \mid n=xyz \}$, has appeared in \cite{Bellman:1949} and in \cite{ConIwan:2000}, in a somewhat different form.

Some comments on the organization of this paper: in the next section we present the representation-theoretic results on which our approach is based, in particular the notion of automorphic distribution. Automorphic distributions for ${\rm GL}(3,\Z)$ restrict to $N_\Z$-invariant distributions on the upper triangular unipotent subgroup $N \subset {\rm GL}(3,\R)$, and they are completely determined by their restrictions to $N$. We analyze these restrictions in \secref{heis}, in terms of their Fourier expansions on $N_\Z\backslash N$. Proposition \ref{prop1} gives a very explicit description of the Fourier decomposition of distributions on $N_\Z\backslash N$; we prove the proposition in \secref{proofofprop1}. Section \ref{voronsec} contains the proof of our main theorem, i.e., of the Voronoi summation formula for ${\rm GL}(3)$. The proof relies heavily on a particular analytic technique -- the notion of a distribution vanishing to infinite order at a point, and the ramifications of this notion.
Since the technique applies in other contexts as well, we are developing it in a separate companion paper
\cite{MilSch:2003b}. We had mentioned already that we derive the functional equations for the
$L$-functions $L(s,\Phi\otimes \chi)$ in \secref{lfuncsec}, using the results of the earlier sections, and that
\secref{convsec} contains our proof of the Converse Theorem of \cite{JacPiatShal:1979}.

It is a pleasure to thank James Cogdell, Dick Gross, Roger Howe, David Kazhdan, Peter Sarnak, and Thomas Watson for their encouragement and helpful comments.

\section{Automorphic distributions}\label{autodist}

For now, we consider a unimodular, type I Lie group\footnote{The type I condition is a technical hypothesis, satisfied in particular by reductive and nilpotent Lie groups; these are the two cases of interest for our investigation.} $G$ and a discrete subgroup $\G\subset G$.  Then $G$ acts on $L^2(\G \backslash G)$ by the right regular representation,
\begin{equation}
\label{rightregular}
(r(g)f)(h)\ = \ f(hg),\ \ \ (\, g \in G, \ \ h \in \G \backslash G\,)\,.
\end{equation}
If $(\G \cap Z_G )\backslash Z_G$, the quotient of the center $Z_G$ by its intersection with $\G$, fails to be compact -- as is the case for $G = {\rm GL}(n,\R)$, $\G = {\rm GL}(n,\Z)$, for example\nobreak{ --} one needs to fix a unitary character $\omega : Z_G \to \C^*$ and work instead with the right regular representation on
\begin{equation}
\label{rightregular1}
\begin{aligned}
L^2_\omega(\G\backslash G)& =   \text{space of all}\ \, f\in L^2_{\text{loc}}(\G\backslash G)\ \ \text{such
that}  \\   f(gz) &= \omega(z)f(g) \, \ \text{for $g\in G$, $z\in Z_G$}\,,\ \text{and}\ \ \int_{\G\backslash
G/Z_G}|f|^2\,dg < \infty \,.
\end{aligned}
\end{equation}
The resulting representation can be decomposed into irreducible constituents \cite{Dixmier:1969}: a direct sum of irreducibles if the quotient $\G\backslash G/Z_G$ is compact, or a ``continuous direct sum" -- i.e., a direct integral -- in general. Even in the case of a direct integral decomposition, direct summands may occur. It is these direct summands we are concerned with. We recall some standard facts.

Let $(\pi,V)$ be an irreducible, unitary representation of $G$, embedded as a direct summand in $L_\omega^2(\G \backslash G)$,
\begin{equation}
\label{embed}
i : V \hookrightarrow L_\omega^2(\G\backslash G).
\end{equation}
One calls $v \in V$ a $C^\infty$ vector if $v\mapsto\pi(g)v$ defines a $C^\infty$ map from $G$ to the Hilbert space $V$. The totality of $C^{\infty}$ vectors constitutes a dense subspace $V^\infty\subset V$, which carries a natural Fr\'echet topology via the identification
\begin{equation}
\label{topology}
V^{\infty} \simeq \{\, f \in C^{\infty}(G,V) \mid f(g) = \pi(g)f(e) \ \text{for all}\ g \in G \,\}\,,\ \ v \leftrightarrow f(e)\,.
\end{equation}
Note that $\pi$ restricts to a continuous representation on this Fr\'echet space.  Dually, $V$ lies inside $V^{-\infty}$, the space of distribution vectors; by definition, the distribution vectors are continuous linear functionals on $(V')^\infty$, the space of $C^{\infty}$ vectors for the irreducible unitary representation $(\pi',V')$ dual to $(\pi,V)$. Thus
\begin{equation}
\label{diffdist}
V^{\infty} \ \subset \ V \ \subset \ V^{-\infty}\,,
\end{equation}
which is consistent with the following convention: we define distributions on a manifold as continuous linear functionals on the space of compactly supported smooth measures. This makes continuous functions, or $L^2$ functions, particular examples of distributions, in analogy to (\ref{diffdist}).

The inclusion (\ref{embed}) sends $C^\infty$-vectors to $C^\infty$ functions, resulting in a continuous, $G$-invariant linear map
\begin{equation}
\label{smoothembed}
i : V^{\infty} \hookrightarrow C^{\infty}(\G\backslash G).
\end{equation}
Since $i(v)$, for $v \in V^{\infty}$, is $\G$-invariant on the left, the composition of $i$ with evaluation at the identity determines a $\G$-invariant, continuous linear functional on $V^{\infty}$ -- in other words, a $\G$-invariant distribution vector for the dual representation $(\pi',V')$:
\begin{equation}
\label{tau}
\tau \in ((V')^{-\infty})^{\G}, \qquad \langle \, \tau\, , \, v \, \rangle \,= \ i(v)(e)\ \ \ \ \text{for}\ v \in V^{\infty}\,.
\end{equation}
This is the {\it automorphic distribution} corresponding to the embedding (\ref{embed}). We remark that $\tau$ completely determines the embedding. Indeed, for $v \in V^{\infty}$ and $g \in G$, $i(v)(g) = (r(g)i(v))(e) = i(\pi(g)v)(e) = \langle \, \tau , \pi(g)v \, \rangle$,
and so $\tau$ does determine the restriction of (\ref{embed}) to $V^{\infty}$, which is dense in $V$, and
hence determines the embedding itself.

Since we work with the automorphic distribution rather than the embedding, it will be more convenient to interchange the roles of $\pi$ and the dual representation $\pi'$. Thus, from now on,
\begin{equation}
\label{newtau}
i : V' \hookrightarrow L_{\bar\omega}^2(\G\backslash G)\,,\ \ \ \tau \in (V^{-\infty})^{\Gamma}\,, \ \ \ \langle \, v \, , \, \tau \, \rangle = \, i(v)(e) \ \ \text{for}\ v \in (V')^{\infty}.
\end{equation}
The natural duality between $L_{\bar\omega}^2(\G\backslash G)$ and $L_{\omega}^2(\G\backslash G)$ makes this reversal of roles legitimate. Even if $\pi$ has central character $\omega$, not every $\tau \in (V^{-\infty})^{\Gamma}$ arises from an embedding of $V'\hookrightarrow L_{\bar\omega}^2(\G\backslash G)$. However, for any such $\tau$ and $v \in (V')^\infty$, the map $g \mapsto \langle \,\pi'(g)v,  \tau \, \rangle$ defines a $\G$-invariant $C^\infty$ function on $G$, so continuous, $G$-invariant homomorphisms from $(V')^\infty$ to $C^\infty(\G\backslash G)$ do correspond bijectively to distribution vectors $\tau \in (V^{-\infty})^{\Gamma}$.

We now specialize our discussion to the case $G = {\rm GL}(3,\R)$. Loosely speaking, any irreducible unitary representation can be realized as a subrepresentation of a not-necessarily-unitary principal series representation \cite{Casselman:1980}. 
To make this precise, we consider the subgroups
\begin{equation}
\label{minpar}
\begin{gathered}
A = \left\{\left.\ttt{a_1}{ }{ }{ }{a_2}{ }{ }{ }{a_3}\ \right| \ a_j > 0 \ \right\} \,,\ \ \ N_- = \left\{ \ttt{1}{ }{ }{*}{1}{ }{*}{*}{1} \right\} \,,\\
M = \left\{\left.\ttt{\e_1}{ }{ }{ }{\e_2}{ }{ }{ }{\e_3}\, \right| \, \e_j \in \{\pm 1\} \, \right\} \ .
\end{gathered}
\end{equation}
Here, as elsewhere, we do not explicitly write out zero matrix entries. Then $MA$ is the full diagonal subgroup of $G$, which normalizes $N_-$. The semidirect product $P = MAN_-$ constitutes a minimal parabolic subgroup. We fix parameters
\begin{equation}
\label{lambda}
\begin{gathered}
\l = (\,\l_1\,,\,\l_2\,,\,\l_3\,)\ \in \C^3 \ \ \ \text{such that}\ \ \ \textstyle{\sum}_{j=1}^3 \,\l_j \, = \, 0\,,
\\
\text{and}\ \ \ \d = (\,\d_1\,,\,\d_2\,,\,\d_3\,)\ \in (\Z/2\,\Z)^3\,,
\end{gathered}
\end{equation}
which we use to define the character
\begin{equation}
\label{chi}
\omega_{\l,\d} : P \to \C^*\,, \qquad \omega_{\l,\d}\ttt{\e_1 a_1}{ }{ }{*}{\e_2 a_2}{ }{*}{*}{\e_3 a_3}\, =\ \frac{a_3}{a_1}\, \prod_{j=1}^3\,a_j^{\l_j} \e_j^{\d_j}\,.
\end{equation}
Via left translation, $G$ acts on
\begin{equation}
\label{vlambda}
\gathered
V_{\l,\d}^{\infty} = \{\, f \in C^{\infty}(G) \mid f(gp) = \omega_{\l,\d}(p^{-1})f(g)\,\ \text{for all}\,\ g \in G, p \in P\,\},\\
(\pi_{\l,\d}(g)f)(h) \ = \ f(g^{-1}h)\,.
\endgathered
\end{equation}
The significance of the factor $a_3/a_1$ in the definition of the inducing character $\omega_{\l,\d}$ will become apparent presently. The hypothesis $\sum_j \l_j =0\,$ means that the identity component $Z_G^0$ of the center of $G$ acts as the identity on $V_{\l,\d}^\infty$. In effect, we are restricting our attention to the case when the central character $\omega$ in (\ref{rightregular1}) is trivial on $Z_G^0$. We can do so without essential loss of generality: since ${\rm SL}(3,\Z)\cap Z_G^0 = \{e\}$, any ${\rm SL}(3,\Z)$-automorphic representation can be twisted by a character of $Z_G^0$ to make $Z_G^0$ act trivially.

In geometric terms, $V_{\l,\d}^{\infty}$ can be regarded as the space of $C^\infty$ sections of a $G$-equivariant $C^\infty$ line bundle $\mathcal L_{\l,\d} \to G/P$. The quotient $G/P$ is compact -- as follows, for example, from the Iwasawa decomposition $G = KAN_-$, with $K = O(3,\R)$. Since $AN_-$ fails to be unimodular, $G/P$ does not admit
a\break $G$-invariant measure. However, any product $f_1f_2$, with $f_1 \in V_{\l,\d}^{\infty}$ and $f_2
\in V_{-\l,\d}^{\infty}$, transforms under $G$ as a smooth measure on $G/P$. Since integration of smooth
measures over the compact manifold $G/P$ has invariant meaning, it follows that there exists a canonical,
$G$-invariant pairing
\begin{equation}
\label{pairing}
V_{\l,\d}^{\infty} \times V_{-\l,\d}^{\infty} \ \longrightarrow \ \C\,.
\end{equation}
This duality between representations with parameters $\l$ and $-\l$ depends on the presence of the factor $a_3/a_1$ in the parametrization of $\omega_{\l,\d}$ in (\ref{chi}).

To make the pairing explicit, we note that $K = O(3,\R)$ acts transitively on $G/P\,$; indeed, $G/P \cong K/M$ since $G = KP$ and $K \cap P = M$. The action of $K$, in particular, preserves the pairing, so
that  it can be described concretely as integration over $K$,
\begin{equation}
\label{pairingoverK}
\langle \, f_1 \, , \, f_2 \, \rangle \ = \ \int_K \, f_1(k) f_2(k)\, dk \qquad (\,f_1 \in V_{\l,\d}^{\infty}\,,\ f_2 \in V_{-\l,\d}^{\infty} \,)\,,
\end{equation}
up to a positive constant which reflects the normalization of measures. For $\l \in i\,\R^3$, the complex conjugate $\bar\l$ coincides with $-\l$. In this situation,
\begin{equation}
\label{innerproductoverK}
(f_1,f_2) \ = \ \int_K \, f_1(k) \overline{f_2}(k)\, dk \qquad (\,f_1, f_2 \in V_{\l,\d}^{\infty}\,,\,\ \l \in i\,\R^3 \,)
\end{equation}
defines a $G$-invariant inner product, and $V_{\l,\d}^{\infty}$ is the space of $C^{\infty}$ vectors for a unitary representation $(\pi_{\l,\d} , V_{\l,\d})$, on the Hilbert space completion of $V_{\l,\d}^{\infty}$. Even without the hypothesis $\l \in i\,\R^3$, there exists a representation $(\pi_{\l,\d} , V_{\l,\d})$ on a Hilbert space -- though not necessarily a unitary representation -- whose space of $C^{\infty}$ vectors coincides with $V_{\l,\d}^{\infty}$.

We now consider an arbitrary irreducible unitary representation $(\pi,V)$ of~$G$. The result of Casselman
\cite{Casselman:1980} that we alluded to before, combined
 with Theorem 5.8 of \cite{Wallach:1983} and specialized to the case at hand, guarantees the existence of
parameters $(\l,\d)$ such that
\begin{equation}
\label{embed2}
V^\infty \ \hookrightarrow \ V_{\l,\d}^{\infty}\,,
\end{equation}
continuously and $G$-invariantly. A deeper result of Casselman-Wallach \cite{Casselman:1989}, \cite{Wallach:1983} implies that this embedding extends continuously, and of course equivariantly, to the spaces of distribution vectors,
\begin{equation}
\label{embed3}
V^{-\infty} \ \hookrightarrow \ V_{\l,\d}^{-\infty}\,.
\end{equation}
Here $V_{\l,\d}^{-\infty}$ can be interpreted in three equivalent ways. On the one hand, it is the space of distribution vectors for the -- possibly nonunitary -- representation which has $V_{\l,\d}^{\infty}$ as the space of $C^\infty$ vectors. It can also be characterized as the space of distribution sections of the line bundle $\mathcal L_{\l,\d} \to G/P$ whose $C^\infty$ sections constitute the space $V_{\l,\d}^{\infty}$. Lastly, it is the space of distributions on $G$ which transform on the right under $P$ according to the character $\omega_{\l,\d}$, as in the distribution analogue of the definition (\ref{vlambda}).

This brings us closer to the idea of an automorphic distribution as a distribution in the usual sense. The subgroup
\begin{equation}
\label{uppern}
N\ = \ \left\{  \, \ttt{1}{*}{*}{}{1}{*}{}{}{1} \,  \right\} \ \subset \ G
\end{equation}
acts freely on its open orbit in $G/P$, i.e., on the image of $NP$ in $G/P$. Since this open Schubert cell is dense, restriction from $G$ to $N$ defines an $N$-invariant inclusion
\begin{equation}
\label{embed4}
j : V_{\l,\d}^{\infty} \ \hookrightarrow \ C^{\infty}(N)\,.
\end{equation}
Via $j$, the representation $\pi_{\l,\d}$ acts on $C^{\infty}(N)$: for $g \in G$ and a generic $n \in N$, we write $g^{-1}n = n_g\, m_g \, a_g \, n_{-,g}\,$, with $n_g \in N$, $m_g \in M$, $a_g \in A$, $n_{-,g} \in N_-\,$; then
\begin{equation}
\label{action}
j(\pi_{\l,\d}(g)v)(n) \ = \ \omega_{\l,\d}((m_g a_g)^{-1})\,jv(n_g)\,,
\end{equation}
as follows from (\ref{vlambda}). When $g^{-1}n$ fails to lie in the open Schubert cell, the right-hand side is undefined, so this equation must be interpreted as the equality of two $C^{\infty}$ functions on their common domain, which is dense. The embedding $j$ extends continuously to the space of distribution vectors,
\begin{equation}
\label{embed5}
j : V_{\l,\d}^{-\infty} \ \longrightarrow \ C^{-\infty}(N)\,,
\end{equation}
but no longer as an injection, since a distribution cannot be reconstructed from its restriction to a dense open subset. However, (\ref{action}) remains valid for distribution vectors $v \in V_{\l,\d}^{-\infty}$, as long as the right-hand side is well-defined.

The composition of the inclusion (\ref{embed3}) and the map (\ref{embed5}) defines a continuous, $N$-invariant linear map $V^{-\infty} \to C^{-\infty}(N)$. If $\tau \in V^{-\infty}$ arises from an embedding $i : V' \hookrightarrow L_{\bar\omega}^2(\G\backslash G)$ as in (\ref{newtau}), we tacitly identify $\tau$ with its image in $C^{-\infty}(N)$, which is $\G \cap N$-invariant:
\begin{equation}
\label{tauagain}
\tau  \in  C^{-\infty}(\G \cap N \backslash N)\,.
\end{equation}
The concrete interpretation of the automorphic distribution $\tau$ with a $\G \cap N$-invariant distribution on $N$ takes notational license in two ways. First of all, it depends on the choice of the embedding (\ref{embed3}), and secondly, the image of a distribution vector in $C^{-\infty}(N)$ does not determine the vector. We deal with the former ambiguity by fixing the embedding throughout the discussion; this is legitimate since the $L$-function we attach to $\tau$ will turn out to be an invariant of $i$. As for the latter, when the discrete subgroup $\G$ is sufficiently large -- e.g., a congruence subgroup of ${\rm GL}(3,\Z)$ -- the $\G$-translates of the open Schubert cell cover all
 of $G/P$, so that  any $\G$-invariant distribution vector is determined by its restriction to $N$, after all.

\section{Fourier series on the Heisenberg group}\label{heis}

We now apply Fourier analysis on $\G \cap N \backslash N$ to automorphic distributions $\tau$ as in (\ref{tauagain}). To simplify the discussion, we let ${\rm GL}(3,\Z)$ play the role of $\G$:
\begin{equation}
\pagebreak
\label{nz}
\G \cap N \backslash N \ = \ N_\Z \backslash N\,,\ \ \text{with}\ \ N_\Z \ = \ {\rm GL}(3,\Z)\cap N\,.
\end{equation}

\noindent
We should remark, however, that the results of this section can be easily extended to the case of a congruence subgroup $\G \subset {\rm GL}(3,\Z)$. It will be convenient to use coordinates on $N$,
\begin{equation}
\label{coordinates}
\R^3 \ \overset{\sim}{\longrightarrow} \ N\,, \qquad (x,y,z) \ \mapsto \ttt{1}{x}{z}{}{1}{y}{}{}{1}\,.
\end{equation}
Then $\,\Z^3$ corresponds to $N_\Z$ and $\{x=y=0\}$ to the center of $N$. In terms of the coordinates, the group law is given by the formula
\begin{equation}
\label{multiplication}
(x_1,y_1,z_1) \cdot (x_2,y_2,z_2) \ = \ (x_1 + x_2,y_1 + y_2,z_1 + z_2 + x_1 y_2)\,.
\end{equation}
Left and right translation on $N$ preserves the measure $dx \, dy \, dz$. Since the inequalities $\,0 \leq x \leq
1\,$, $0 \leq y \leq 1\,$, $0 \leq z \leq 1\,$ cut out a fundamental domain for the action of $N_\Z$ on $N$,
\begin{equation}
\label{measure}
\int_{N_\Z \backslash N}  dx \, dy \, dz \ = \  \int_0^1 \! \int_0^1 \! \int_0^1 dx \, dy \, dz \, = \, 1\,;
\end{equation}
in other words, $dx \, dy \, dz$ represents Haar measure normalized so as to assign total measure one to the quotient $N_\Z \backslash N$.

The irreducible unitary representations of the three dimensional Heisenberg group $N$ are well known \cite{Mackey:1958}. First of all there are the one-dimensional unitary representations
\begin{equation}
\label{abelian}
(x,y,z)\ \mapsto \ e(ax + by)\,, \ \ \ \text{with}\ \ (a,b) \in \R^2\,.
\end{equation}
Any such character, considered as a  function on $N$, is $N_\Z$-invariant if and only if $a$ and $b$ are
integers. It follows that the functions $e(rx + sy)$, with $(r,s) \in \Z^2$, constitute a Hilbert space basis of the
largest subspace of $L^2(N_\Z \backslash N)$ on which the center of $N$ acts trivially.

Next we fix a nontrivial character of the center. Since we are interested in $N_\Z$-invariant functions, we only consider nontrivial characters that restrict trivially to the intersection of the center with $N_\Z$. These are precisely the central characters
\begin{equation}
\label{centralcharacter}
(x,y,z)\ \mapsto \ e(n z)\,, \ \ \ \text{with}\ \ n \in \Z - \{0\}\,.
\end{equation}
Up to isomorphism, there exists exactly one irreducible unitary representation $(\pi_n,V_n)$ of $N$ with this central character. It has two different, but equally natural models: in both cases on the Hilbert space $V_n = L^2(\R)$, one with action
\begin{equation}
\label{nonabelian1}
(\pi_n(x,y,z)f)(t) = e(n(z + ty))f(x + t)\,,
\end{equation}
the other, with action
\begin{equation}
\label{nonabelian2}
(\widehat\pi_n(x,y,z)h)(t) = e(n(z - xy + tx))\, h(t - y)\,.
\end{equation}
The two actions are intertwined by the Fourier transform and a scaling of the argument,
\begin{equation}
\label{intertwine}
f(t) \ \longleftrightarrow \ h(t) \ = \ \widehat f(nt) \,.
\end{equation}
Here, as always, we normalize the Fourier transform according to Laurent Schwartz' convention,
\begin{equation}
\label{fourierint}
\widehat f(t) \ = \ \int_\R f(u) \, e(-ut)\,du\,.
\end{equation}
There is another relation between the two models: the outer automorphism
\begin{equation}
\label{outer}
(x,y,z) \ \longmapsto \  (\,-y\,,\, -x \,, -z + xy\,)
\end{equation}
of $N$ conjugates $\widehat\pi_n$ into $\pi_{-n}$. Note that the $N$-invariant pairing
\begin{equation}
\label{hbpairing}
V_{-n} \times V_n \ \longrightarrow \ \C\,, \qquad \langle \, f_1 \, , \, f_2 \, \rangle \ = \ \int_\R \, f_1(t) f_2(t)\, dt \,,
\end{equation}
exhibits $(\pi_{-n},V_{-n})$ as the dual of $(\pi_n,V_n)$, and simultaneously $(\widehat\pi_{-n},V_{-n})$ as the dual of $(\widehat\pi_n,V_n)$.

The partial derivatives $\frac{\partial\ }{\partial x}$, $\frac{\partial\ }{\partial y}$, $\frac{\partial\ }{\partial z}$ at the origin in $\R^3$ span the Lie algebra of $N$. If $f$ is a $C^\infty$ vector for $\pi_n$, the identities
\begin{equation}
\label{partials}
\textstyle \pi_n(\frac{\partial\ }{\partial x})f\, =\, \frac{\partial\ }{\partial t}f\,,\ \ \ \pi_n(\frac{\partial\ }{\partial y})f \, = \, 2\,\pi\, i\, n \, t f\,,\ \ \ \pi_n(\frac{\partial\ }{\partial z})f \, =\, 2\,\pi\, i\, n f
\end{equation}
imply the square-integrability of the function $t \mapsto t^k f^{(\ell)}(t)$ for all $k, \ell \in \mathbb N$, so $f$ must be a Schwartz function. Conversely, for any Schwartz function $f$, $(x,y,z) \mapsto \pi_n (x,y,z)f$ visibly defines a $\C^\infty$ map from $N$ to the Schwartz space $\mathcal S(\R)$, hence in particular a $\C^\infty$ map from $N$ to $L^2(\R)$. One can argue the same way in the case of $\widehat\pi_n$. Thus, for both actions,
\begin{equation}
\label{schwartz}
V_n^\infty \ \simeq \ \mathcal S(\R)\,,\ \ \ \text{and dually,} \ \ \ V_{-n}^{-\infty} \ \simeq \ \mathcal S'(\R)\,,
\end{equation}
i.e., distribution vectors are tempered distributions. In analogy to (\ref{hbpairing}), we denote the pairing between $V_{-n}^{-\infty}$ and $V_{n}^{\infty}$ by integration.

We fix a nonzero integer $n \in \Z - \{0\}$ and a residue class $k \in \Z/ n\Z$. For any choice of $\sigma, \rho \in \mathcal S'(\R)$, the expressions
\begin{equation}
\label{nonabelian}
\begin{aligned}
(x,y,z) \ &\longmapsto \ {\sum}_{\ell\equiv k(\md n)}\ e(nz + \ell y)\,\sigma(x+{\textstyle\frac{\ell}{n}}),\\
(x,y,z) \ &\longmapsto \ {\sum}_{\ell\equiv k(\md n)}\ e(n(z-xy) + \ell x)\,\rho({\textstyle\frac{\ell}{n}}-y)
\end{aligned}
\end{equation}
define distributions on $N$; what matters here is the temperedness of $\sigma$, $\rho$ and the fact that the summation simultaneously involves a translation in one variable and multiplication by powers of a nontrivial character 
 in the other. Using (\ref{coordinates}), (\ref{multiplication}), one finds that these distributions are
$N_\Z$-invariant on the left, i.e., they lie in $C^{-\infty}(N_\Z \backslash N)$. Moreover, the first of the
two depends\break $N$-equivariantly on $\sigma$ when $N$ acts on $V_n^{-\infty} \cong \mathcal S'(\R)$
via
$\pi_n$ and on $C^{-\infty}(N_\Z\backslash N)$ via the right regular representation (\ref{rightregular}),
whereas the second expression depends $N$-equivariantly on $\rho$ relative to the action $\widehat\pi_n$
on $V_n^{-\infty} \cong \mathcal S'(\R)$. Our next statement involves the Fourier transform of tempered
distributions on $\R$ and the finite Fourier transform on the set $\Z/n\Z$. We define the former by the
identity
\begin{equation}
\label{fourierint2}
\int_\R \widehat\sigma(t)\, f(t) \, dt \ = \ \int_\R \sigma(t) \, \widehat f(t) \, dt \ \ \ \text{for all}\ \ \sigma \in \mathcal S'(\R)\,,\ \, f \in \mathcal S(\R)\,,
\end{equation}
in accordance with our convention of regarding the notion of distribution as an extension of the notion of function. In the definition of the finite Fourier transform of $\,a = (a_k)_{k \in \Z/n\Z}\,$,
\begin{equation}
\label{fourierint3}
\widehat a_k \ =  \ {\sum}_{\ell \in \Z/n\Z}\ e({\textstyle\frac{k \ell}{n}})\, a_\ell \,,
\end{equation}
we follow a common convention that omits the normalizing factor and complex conjugation of the character customary in representation theory.

\begin{prop}\label{prop1} Any $\,\tau \in C^{-\infty}(N_\Z\backslash N)$ has Fourier expansion
\begin{eqnarray*}
\tau(x,y,z)   &=  &  {\sum}_{r,s\in\Z}\ c_{r,s}\ e(rx + sy) \\
&&+  {\sum}_{n\in \Z - \{0\}} \ {\sum}_{k\in\Z/n\Z} \left( {\sum}_{\ell\equiv k(\operatorname{mod} n)}
\  e(nz+ \ell y)\,\sigma_{n,k}(x + \ell/n)\right),
\end{eqnarray*}
with $\,\sigma_{n,k} \in \mathcal S'(\R)$. The series converges in the strong distribution topology on
$C^{-\infty}(N_\Z\backslash N)$. The contribution on the right indexed by any $\,n \in \Z - \{0\}$ can be
written alternatively as
\begin{multline*}
 {\sum}_{k\in\Z/n\Z}\ {\sum}_{\ell\equiv k(\md n)} \  e(n z+ \ell y)\,\sigma_{n,k}(x + \ell/n) 
\\
   = \ {\sum}_{k\in\Z/n\Z}\ {\sum}_{\ell\equiv k(\operatorname{mod} n)} \  e(n(z-xy) +  \ell
x)\,\rho_{n,k}(\ell/n-y)\,,
\end{multline*}
in terms of distributions $\rho_{n,k}\in \mathcal S'(\R)$ which are related to the Fourier transforms of the
$\sigma_{n,k}$ by the identities
$$
{\sum}_{k \in \Z/n\Z}\ a_k\,\rho_{n,k}(y) \ =  \ {\sum}_{k \in \Z/n\Z}\  \widehat a_k\,\widehat{\sigma}_{n,k}(ny)\,;
$$
here the coefficients $\,a_k$, $k \in \Z/n\Z${\rm ,} can be chosen arbitrarily{\rm ,} and $(\widehat a_k)$ denotes the finite Fourier transform of $(a_k)${\rm ,} normalized as in {\rm (\ref{fourierint3})}.
\end{prop}

We shall refer to the $c_{r,s}$ as the {\it abelian Fourier coefficients\/} of $\tau$ since they are the coefficients of the abelian characters of $N$ in the Fourier expansion. On the other hand, the expressions (\ref{nonabelian}), with $\sigma_{n,k}$ and $\rho_{n,k}$ in place of $\sigma$ and $\rho$, should be viewed as the {\it non-abelian Fourier components} of $\tau$.

The proof of the proposition occupies Section \ref{proofofprop1} below. We finish the current section with
some fairly immediate consequences of the statement of the proposition. Let us suppose now that $\tau$
arises from a discrete summand $i : V' \hookrightarrow L_{\bar\omega}^2(\G\backslash G)$ as in
(\ref{newtau}), with $G = {\rm GL}(3,\R)$ and $\G = {\rm GL}(3,\Z)$, via an embedding $V^{-\infty}
\hookrightarrow V_{\l,\d}^{-\infty}$ as in (\ref{embed3}). We can then regard $\tau$ as
an\break $N_\Z$-invariant distribution on $N$, as in (\ref{tauagain}), so the notation of
Proposition~\ref{prop1} applies.

Recall the parametrization (\ref{coordinates}) of $N$. It makes the linear subspaces $\{y=0\}$, $\{x=0\}$ correspond to subgroups $N_{x,z}$, $N_{y,z}$ of $N$. By definition, the inclusion $i$ is cuspidal if
\begin{equation}
\label{cuspidality}
\int_{N_\Z \cap  N_{x,z}\backslash N_{x,z}} \!\! i(v)(n)\, dn \ = \ 0\ = \ \int_{N_\Z \cap  N_{y,z}\backslash N_{y,z}} \!\! i(v)(n)\, dn
\end{equation}
for all $v \in V'$, or equivalently, for all $v$ in the dense subspace $(V')^{\infty}$.

\begin{lem}\label{lem1} If the inclusion $i : V' \hookrightarrow L_{\bar\omega}^2(\G\backslash G)$ 
corresponding to $\tau$ is cuspidal{\rm ,} the coefficients $c_{r,0}${\rm ,} $c_{0,s}$ vanish{\rm ,} for all $r, s\in\Z$.
\end{lem}

The lemma has a partial converse, which is far more subtle -- see the proof of \lemref{lemexist}. 

\Proof  Since $i(v)(n) = \langle\, \pi'(n)v\,, \, \tau \,\rangle = \langle\, v \,, \, \pi(n^{-1})\tau \,\rangle$, the vanishing of the two integrals, for every $v \in (V')^\infty$, is equivalent to the vanishing of the distribution vectors
\begin{equation}
\label{cuspidality1}
\tau_{x,z} \ = \ \int_{N_{x,z} / N_{x,z}\cap N_\Z}  \pi(n)\tau\, dn \ , \qquad \tau_{y,z} \ = \ \int_{N_{y,z}/N_{y,z}\cap N_\Z}  \pi(n)\tau\, dn \,.
\end{equation}
If the variable $n \in N_{x,z}$ in the first integral corresponds to $(t,0,u)$ under the parametrization (\ref{coordinates}), applying $\pi(n)$ to the distribution $\tau$ results in the distribution $\tau(x-t,y,z-u-ty)$. We now appeal to Proposition \ref{prop1} and find
\begin{equation}
\label{cuspidality2}
\tau_{x,z}(x,y,z)\, = \, \int_0^1\!\! \int_0^1 \!\! \tau(x-t,y,z-u-ty)\, dt\, du \, = \, {\sum}_{s\in \Z} \ c_{0,s}\,e(sy)\,.
\end{equation}
Thus, if $\tau_{x,z} = 0$, the coefficients $c_{0,s}$, $s \in \Z$, all vanish. Similarly, $\tau_{y,z} = 0$ implies $c_{r,0}=0$ for all $r \in \Z$.
\Endproof\vskip4pt  

We now look at the action of the finite group $M$, defined in (\ref{minpar}). Since ${\rm GL}(3,\Z)$ contains $M$, the distribution $\tau$ must be $M$-invariant. Recall that $\tau$ has meaning as
a  distribution on $N$ via the embedding (\ref{embed3}) and restriction of distributions from $G$ to $N$. In
view of the transformation law (\ref{action}), for $m \in M$ and $n \in N$,
\begin{equation}
\label{delta1}
\tau(n) \ = \ (\pi_{\lambda,\d}(m)\tau) (n)\ = \ \tau(m^{-1}n) \ = \ \omega_{\lambda,\delta}(m)\,\tau( m^{-1} n \, m)\,.
\end{equation}
Written in terms of the coordinates on $N$,
 conjugation by a diagonal matrix $m$ with 
diagonal entries $\e_j$ sends $(x,y,z)$ to $(\e_1 \e_2 x,\e_2 \e_3
y,\e_1 \e_3z)$; hence
\begin{equation}
\label{delta2}
\tau(\e_1 \e_2 x \,,\,\e_2 \e_3 y\,,\,\e_1 \e_3z) \ = \ {\prod}_{j=1}^3 \e_j^{\delta_j} \ \tau(x,y,z) \ \ \ \ (\, \e_1,\, \e_2,\, \e_3 \in \{\pm1\}\,).
\end{equation}
In particular, a nonzero $\tau$ can exist only if $\d_1 + \d_2 + \d_3 = 0$ -- this is analogous to the nonexistence of modular forms of odd weight for ${\rm SL}(2,\Z)$. Thus we explicitly require
\begin{equation}
\label{delta3}
\d_1 + \d_2 + \d_3 \ = \ 0\,.
\end{equation}
We now combine (\ref{delta2}) with Proposition \ref{prop1} and conclude:

\begin{lem}\label{lem2} For all choices of indices $k, n, r, s${\rm ,}
\begin{eqnarray*}
c_{-r,s} &=& (-1)^{\delta_1}\,c_{r,s}\,,\ \ \sigma_{-n,k}(x) = (-1)^{\delta_1}\,\sigma_{n,k}(-x)\,,\ \
\rho_{-n,-k}(y) = (-1)^{\delta_1}\,\rho_{n,k}(y)\,,
\\
c_{r,-s} &=& (-1)^{\delta_3}\,c_{r,s}\,,\ \ \sigma_{-n,-k}(x) = (-1)^{\delta_3}\,\sigma_{n,k}(x)\,,\ \
\rho_{-n,k}(y) = (-1)^{\delta_3}\,\rho_{n,k}(-y)\,.
\end{eqnarray*}
\end{lem}

Our next statement relates the non-abelian Fourier components to the abelian coefficients. We fix two relatively prime integers $a$, $c$, with $c \neq 0$, and choose an integer $\bar a \in \Z$ which represents the reciprocal of $a$ modulo $c$:
\begin{equation}
\label{acabar}
a,\, \bar a \in \Z\,,\qquad c \in \Z - \{0\}\,,\qquad a\,\bar a \equiv 1\ \,(\md c)\,.
\end{equation}

\begin{prop}\label{prop2} Under the hypotheses just stated{\rm ,} for any $\,q \in \Z - \{0\}${\rm ,}
\begin{alignat*}{5}
 {\rm a)}\quad    &\sigma_{cq,aq}(x) \ = \ (\sg cx)^{\d_3}\, |cx|^{\l_1 - \l_2 - 1}\,{\sum}_{r \in \Z} \ c_{r,q}
\, e(\, r\, \bar a \, c^{-1}  - r\, c^{-2}\, x^{-1})\,,\qquad\qquad\qquad\qquad
\\
 {\rm b)}   \quad &\rho_{cq,aq}(y) \ = \ (\sg cy)^{\d_1}\, |cy|^{\l_2 - \l_3 - 1}\,{\sum}_{s \in \Z} \ c_{q,s} \,
e( \, s \, c^{-2} \, y^{-1}  -  s\, \bar a \, c^{-1})\,.\\
\intertext{The coefficients $c_{r,0}${\rm ,} $c_{0,s}$ satisfy the relations}
 {\rm c)}   \quad &{\sum}_{r\in \Z}\ c_{r,0} \, e(rx) \ = \  \begin{cases}\ \ |x|^{\l_1 - \l_2 -
1}{\displaystyle{\sum}_{r \in \Z}} \ c_{r,0} \, e( \,-\, r x^{-1}) \ \ &\text{if}\ \ \d_3 = 0\\
\ \ 0 \ \ &\text{if}\ \ \d_3 =1 \,, \end{cases}
\\
 {\rm d)} \quad& {\sum}_{s\in \Z}\ c_{0,s} \, e(sy) \ = \  \begin{cases}\ \ |y|^{\l_2 - \l_3 -
1}{\displaystyle{\sum}_{s \in \Z}} \ c_{0,s} \, e( \,-\, s y^{-1}) \ \ &\text{if}\ \ \d_1 = 0\\
\ \ 0 \ \ &\text{if}\ \ \d_1 = 1 \,.\end{cases}
\end{alignat*}
\end{prop}

The equations a)--d) can be interpreted as identities between distributions on $\,\R^*$. However, the proof
establishes more: an equality between distribution vectors in appropriately defined representation spaces for
${\rm SL}(2,\R)$. Details will be given in \corref{cor1}, following the proof of the proposition. 

\Proof  We begin with a). Because of (\ref{acabar}), there exists $b \in \Z$ such that $a\bar a - bc = 1$. The first matrix factor on the left in the identity
\begin{multline*}
 \ttt{a}{-b}{}{-c}{\bar a}{}{}{}1 \ttt {1}{x + \bar a c^{-1}}{z}{}1{y}{}{}1 
\\
 = \ttt 1{-c^{-2} x^{-1} - a c^{-1}}{a z - b y}{}1{-cz + \bar a y}{}{}1  \ttt{- c^{-1}
x^{-1}}{}{}{-c}{-cx}{}{}{}{1}
\end{multline*}
lies in ${\rm SL}(3,\Z)$. In view of (\ref{vlambda}), (\ref{delta3}), and the ${\rm GL}(3,\Z)$-invariance of
$\tau$, the matrix identity implies
\begin{equation}
\label{taucondition1}
\begin{aligned}
\tau(\,x + \bar a \, c^{-1} ,\, y \, ,\, z \,)\  = &\ (\sg( -cx))^{\d_3}\, |cx|^{\l_1 - \l_2 - 1}  
\\[5pt]
&\times \, \tau( -c^{-2}\, x^{-1} - a\,c^{-1}\, ,\, \bar a\, y - c\, z \, ,\, -b y + a \, z \, )\,.
\end{aligned}
\end{equation}
We now equate the Fourier components on both sides that transform according to any nontrivial character in the variable $y$ and the trivial character in the variable $z$. Since $c \neq 0$, each of the terms $\,c_{r,s}\, \textstyle e(r( -\frac{1}{c^2 x} - \frac {a}{c}) + s(\bar a y - cz))$ either involves $z$ nontrivially or does not involve $y$ at all. Consequently these terms do not contribute. We apply Proposition \ref{prop1} and conclude
\begin{multline*}
 (\sg(-cx))^{\d_3}\, |cx|^{\l_2 - \l_1 + 1}\,{\sum}_{r \in \Z, \, s \neq 0}\ c_{r,s} \, e( r( {\textstyle x +
\frac{\bar a}{c}}) + s y) 
\\[4pt]
  =  \sum_{\stackrel{\scriptstyle{n \neq 0}}{k\in\Z/n\Z}} \ \sum_{\stackrel{\scriptstyle{\ell\equiv k(\md
n)}}{c \ell = n a,\ n b \neq \ell \bar a}} \!\!\!\!\! e( n(-b y + a z) + \ell(\bar a y - c z))\,\sigma_{n,k}( {\textstyle
-\frac{1}{c^2 x} - \frac {a}{c}} + {\textstyle\frac{\ell}{n}} )\,.
\end{multline*}
If $c\ell=na$, the identity $a\bar a - bc = 1$ implies $\ell \bar a  - n b = n/c$, which cannot vanish; in
particular, $c\ell=na$ implies $nb \neq \ell \bar a$. This allows us to replace the three sums on the right by a
single sum, over  nonzero integers $n$ such that $\ell = nac^{-1}$ is integral, with $k$ denoting the residue
class of $\ell$ modulo $n$. Since $a$ and $c$ are relatively prime, we must sum over $n = cq$, with $q \neq
0$, and set $k = \ell = aq$. For these values of $n$, $\ell \bar a - n b = q$ and $\ell / n - a / c = 0$, hence
\begin{equation}
\label{sigmacondition1}
\begin{aligned}
{}&(\sg(-cx))^{\d_3}\, |cx|^{\l_2 - \l_1 + 1}\,{\sum}_{r \in \Z, \, s \neq 0}\ c_{r,s} \, e( r(  x + \bar a c^{-1})
 + s y)   
 \\[4pt]
{}& \qquad\qquad\qquad\ \ \  = \ {\sum}_{q \neq 0}\ e(\,q\,y\,)\,\sigma_{cq,aq}(\, -c^{-2}\, x^{-1}\,
)\,.
\end{aligned}
\end{equation}
We equate the coefficients of $e(qy)$ on both sides and replace $x$ by $-c^{-2} x^{-1}$, to obtain part a) of the proposition.

The verification of b) proceeds quite analogously. However, instead of using the first identity in Proposition \ref{prop1} directly, we use the one
 obtained
 from it by expressing the $\sigma_{n,k}$ in terms of the $\rho_{n,k}$:
\begin{equation}\label{rhoexpansion}
\begin{aligned}   
\tau(x,y,z)   =   &\sum_{r,s\in\Z}\ c_{r,s}\ e( rx + sy) 
 \\{}
& + \ \sum_{n \neq 0} \ \sum_{k\in\Z/n\Z} \left( \sum_{\ell\equiv k(\md n)}  e(n(z-xy) + \ell
x)\,\rho_{n,k}(\ell / n - y)\right).
\end{aligned}
\end{equation}
Because of (\ref{vlambda}), (\ref{delta3}), and the ${\rm GL}(3,\Z)$-invariance of $\tau$, the identity
\begin{multline*}
 \ttt{1}{}{}{}{a}{b}{}{c}{\bar a}\ttt {1}{x}{z + xy - \bar a c^{-1}x}{}1{y - \bar a c^{-1}}{}{}1 
\\
 = \ttt 1{\ \bar a x - cz\ }{c^{-1}(zy^{-1} + x - \bar a c^{-1} xy^{-1})}{}1{a c^{-1} - c^{-2}y^{-1}}{}{}1 
\!\! \ttt 1{}{}{}{c^{-1} y^{-1}}{}{}{c}{cy}
\end{multline*}
implies
\begin{eqnarray}
 \qquad &\tau(\,x\,,\,  y - \bar a c^{-1}\, ,\,z + xy - \bar a c^{-1} x\,)\ = \ (\sg cy)^{\d_1}\, |cy|^{\l_2 -
\l_3 - 1}\
\times
\nonumber\\[-4pt]
&&\label{taucondition2}
\\[-6pt]
&\qquad\ \ \ \times \ \tau(\,\bar a x - cz \, , \, a c^{-1} - c^{-2} y^{-1} \,, \, c^{-1}(zy^{-1} + x - \bar a c^{-1} xy^{-1})\,)\,.
\nonumber
\end{eqnarray}
Next we express $\tau$ in terms of the $c_{r,s}$ and $\rho_{n,k}$, as above. The formulas simplify considerably because what enters as
an  argument of $\rho_{n,k}$ in the expression (\ref{rhoexpansion}) is not $z$ itself but $z - xy$. Applied to
the arguments of $\tau$ on the left and the right, respectively, this substitution gives
\begin{gather*}
 (z + xy - \bar a c^{-1} x)\, - \, x( y - \bar a c^{-1} ) \ = \ z\,,
\\
 c^{-1}(zy^{-1} + x - \bar a c^{-1}xy^{-1}) - ( \bar a x - c z) ( a c^{-1} - c^{-2} y^{-1} ) \ = \ -\,b x\, + \, a
z \,.
\end{gather*}
We now equate the terms on both sides which are constant in $z$ and transform in the variable $x$ according
to any nontrivial character. Arguing exactly as in the proof of a), we find\begin{equation}
\label{rhocondition1}
\begin{aligned}
(\sg cy)^{\d_1}\, |cy|^{\l_3 - \l_2 + 1}\!\! \sum_{r \neq 0, \, s \in \Z} 
c_{r,s} \, &e(\,r \, x + s(y - \bar a \, c^{-1})\,) 
\\
&\ \ =\ \sum_{q \neq 0}\,e(\,q\,x\,)\, \rho_{cq,aq}(\,c^{-2} \, y^{-1}\,)\,.
\end{aligned}
\end{equation}
Isolating the coefficients of $\,e(qy)\,$ on both sides and substituting $\,c^{-2}y^{-1}$ for $y$ gives the formula b).

According to Lemma \ref{lem2}, $\d_3 = 1$ implies $c_{r,0} = 0$ and $\d_1 = 1$ implies $c_{0,s} = 0$. This covers two of the four cases in c) and d). For the proof of the remaining two, we set $a = \bar a = 0$, $b = 1$, $c = -1$ in the identities 
(\ref{taucondition1}), (\ref{taucondition2}). In the former, we express $\tau$ as in Proposition \ref{prop1}
and equate the terms on both sides which transform according to the trivial character in both $y$ and $z$;
when
$\d_3 = 0$, this immediately gives the first case in c). Similarly, for the first case in d), we express $\tau$ in
the two sides of the equation (\ref{taucondition2}) as in (\ref{rhoexpansion}) and equate the terms which
transform trivially under both $x$ and $z$.
\Endproof\vskip4pt  

In order to extend the validity of the identities a)-d) in \propref{prop2}, we need to interpret the distributions $\sum_r c_{r,q} e(rx)$, $\sum_s c_{q,s} e(sy)$, $\sigma_{n,k}$ and $\rho_{n,k}$ as distribution vectors for certain representations of ${\rm SL}(2,\R)$. Corresponding to the
data of $\mu \in \C$ and $\eta \in \Z/2\Z$, we define
\begin{equation}
\label{wmuinf1}
W_{\mu,\eta}^\infty\ =\ \{\,f \in C^\infty({\rm SL}(2,\R))\, \mid\, f\left(g\left( \begin{smallmatrix} 1/a & 0 \\ c & a \end{smallmatrix}\right)\right) \equiv (\sg a)^\eta |a|^{\mu - 1}f(g)\,\}\,,
\end{equation}
on which $\,{\rm SL}(2,\R)$ acts by left translation. Alternatively and equivalently,
\begin{equation}
\label{wmuinf2}
W_{\mu,\eta}^\infty\ =\ \{\,f \in C^\infty(\R)\, \mid\, (\sg x)^\eta\, |x|^{\mu -1}f(-1/x) \in C^\infty(\R)\,\}\,,
\end{equation}
with action
\begin{eqnarray}
(\psi_{\mu,\eta}(g^{-1})f)(x)\ &=& (\sg (cx+d))^\eta\, \textstyle |cx+d|^{\mu -1}f(\frac{ax+b}{cx+d}) \,,
\nonumber\\[-4pt]
&&\label{wmuaction}\\[-9pt]
\text{for}\ \ \ g  &=& \ttwo{a}{b}{c}{d} \ \in \ {\rm SL}(2,\R)\,.
\nonumber
\end{eqnarray}
The space of distribution vectors for this representation is
\begin{equation}
\label{wmudist}
W_{\mu,\eta}^{-\infty}\ =\ \text{strong dual space of}\ \, W_{-\mu,\eta}^{\infty}\,,
\end{equation}
which can be defined as in (\ref{wmuinf1}), with the $C^\infty$ condition replaced by $C^{-\infty}$. In other words, each $\sigma \in W_{\mu,\eta}^{-\infty}$ can be regarded as a distribution $\sigma \in C^{-\infty}(\R)$, together with a specific extension of $(\sg x)^\eta\, |x|^{\mu -1}\sigma(-1/x)$ across $x=0$. The action of ${\rm SL}(2,\R)$ on this space is also given by the formula (\ref{wmuaction}). However, one needs to be
careful to interpret this equality at $x= -d/c$. For details see the discussion in \secref{autodist} of the
analogous construction of the representations $V_{\lambda,\d}^{-\infty}$ of ${\rm GL}(3,\R)$.

\begin{cor}\label{cor1} 
The distributions $\,\sigma_{n,k} \,,\, \sum_{r\in \Z}\, c_{r,q}\, e(rx) \in C^{-\infty}(\R)$ extend naturally to
vectors in the representation space $W_{\lambda_1 - \lambda_2,\d_3}^{-\infty}$. Similarly  $\rho_{n,k}$
and $\,\sum_{s\in \Z}\, c_{q,s}\, e(sy)$ extend to vectors in $W_{\lambda_2 - \lambda_3,\d_1}^{-\infty}$.
With this interpretation{\rm ,} equations {\rm a)--d)} in Proposition {\rm \ref{prop2}}
 can be stated as follows:
\begin{eqnarray*}
{\sum}_{r\in \Z}\, c_{r,q}\, e(rx)&= & \left(\psi_{\lambda_1 - \lambda_2,\d_3}{\ttwo {\bar
a}{-c^{-1}}{c}{0}}\,\sigma_{cq,aq} \right)(x)\ ,
\\
 {\sum}_{s\in \Z}\, c_{q,s}\, e(-sy)&=& \left(\psi_{\lambda_2 - \lambda_3,\d_1}{\ttwo {\bar
a}{-c^{-1}}{c}{0}}\,\rho_{cq,aq} \right)(y)\ ,
\\
 {\sum}_{r\in \Z}\, c_{r,0}\, e(rx)& =& \left(\psi_{\lambda_1 - \lambda_2,\d_3}{\ttwo
{0}{1}{-1}{0}}{\sum}_{r\in \Z}\, c_{r,0}\, e(r\,\cdot) \right)\!(x)\ \ \ (\,\d_3 = 0\,),
\\
 {\sum}_{s\in \Z}\, c_{0,s}\, e(sy)& = & \left(\psi_{\lambda_2 - \lambda_3,\d_1}{\ttwo
{0}{1}{-1}{0}}{\sum}_{s\in \Z}\, c_{0,s}\, e(s\,\cdot) \right)\!(y)\ \ \ (\,\d_1 = 0\,).
\end{eqnarray*} 
The first two of these identities   depend only on $\,\bar a \in \Z/c\Z\,${\rm ,} not on the par\-ti\-cu\-lar choice
of $\,\bar a \in \Z$. 
\end{cor} 

Recall the definition of the abelian subgroups $N_{x,z},\,N_{y,z} \subset N\,$ just before the identity (\ref{cuspidality}). To see how the proof of \propref{prop2} implies the corollary, we introduce the projection operators
\begin{equation}\label{proj1}
\begin{aligned}
 p_{x,k}\,,\ p_{z,\ell} : (V_{\lambda,\d}^{-\infty})^{N_{x,z}\cap N_\Z}& \longrightarrow  
(V_{\lambda,\d}^{-\infty})^{N_{x,z}\cap N_\Z}\,,  
\\
  p_{x,k}\tau& =   \int_0^1 e(kx)\,\pi_{\lambda,\d}\left(\begin{smallmatrix} 1 & x & \ \\ \ & 1 & \ \\ \ & \ &
1
\end{smallmatrix}\right)\tau\ dx\,,
\\
  p_{z,\ell}\tau& =  \int_0^1 e(\ell z)\,\pi_{\lambda,\d}\left(\begin{smallmatrix} 1 & \ & z \\ \ & 1 & \ \\ \ & \
& 1 \end{smallmatrix}\right)\tau\ dz\,,
\end{aligned}
\end{equation}
indexed by $k,\, \ell \in \Z$, and the analogously defined projections
\begin{equation}
\label{proj2}
p_{y,k}\,,\ p_{z,\ell} : (V_{\lambda,\d}^{-\infty})^{N_{y,z}\cap N_\Z}\  \longrightarrow \ (V_{\lambda,\d}^{-\infty})^{N_{y,z}\cap N_\Z}\,.
\end{equation}
Using the common notation $p_{z,\ell}$ in both instances is justified: both extend naturally to the space of invariants for $N_{x,z}\cap N_{y,z}\cap N_\Z$, on which they coincide. Since $N_{x,z}\cap N_{y,z}$ is the center of $N$, the projection $p_{z,\ell}$ maps $N_\Z$-invariants to $N_\Z$-invariants,
\begin{equation}
\label{proj3}
p_{\ell,z} : (V_{\lambda,\d}^{-\infty})^{N_\Z}\  \longrightarrow \ (V_{\lambda,\d}^{-\infty})^{N_\Z}\,.
\end{equation}
The centrality of $N_{x,z}\cap N_{y,z}$ in $N$ also implies the commutation relations
\begin{equation}
\label{proj4}
p_{k,x}\circ p_{\ell,z}\ = \ p_{\ell,z}\circ p_{k,x}\,,\ \ \ p_{k,y}\circ p_{\ell,z}\ =\ p_{\ell,z}\circ p_{k,y}\,.
\end{equation}
Finally, for future reference, we observe that
\begin{equation}
\label{proj5}
p_{k,x}\circ p_{0,z} : (V_{\lambda,\d}^{-\infty})^{N_\Z} \to  (V_{\lambda,\d}^{-\infty})^{N_\Z}\,,\ \ \ p_{k,y}\circ p_{0,z} : (V_{\lambda,\d}^{-\infty})^{N_\Z} \to  (V_{\lambda,\d}^{-\infty})^{N_\Z}\,,
\end{equation}
for all $k \in \Z$. What matters here is the fact that $N_{x,z}$ and $N_{y,z}$ commute modulo the center of $N$, which acts as the identity on the image of $p_{0,z}$.

When the restriction to $N$ of an automorphic distribution $\tau \in (V_{\lambda,\d}^{-\infty})^{G_\Z}$ is expressed as in \propref{prop1}, one finds
\begin{eqnarray}
 \left(p_{y,k}\circ p_{z,n}\, \tau\right)|_N \, (x,y,z) &=& e(nz + ky)\,\sigma_{n,k}(x + k/n)\ \ \ \ (n\neq 0)\,,
\nonumber\\[-4pt]
&&\label{proj6}
\\[-8pt]
 \left(p_{y,q}\circ p_{z,0}\, \tau\right)|_N \, (x,y,z) & = & {\sum}_{r \in \Z}\, c_{r,q}\, e(rx + qy)\,.
\nonumber
  \end{eqnarray}
We now restrict $\,p_{y,k}\circ p_{z,n}\, \tau\,$ and $\,p_{y,q}\circ p_{z,0}\, \tau\,$ to $\,S\cdot N_{y,z} = N_{y,z}\cdot S\,$, where
\begin{equation}
\label{proj7}
S \ = \ \left\{ \ttt{a}{b}{}{c}{d}{}{}{}{1} \in \ {\rm SL}(3,\R) \right\}\ \cong \ {\rm SL}(2,\R)\,.
\end{equation}
Note that $\,S\cdot N_{y,z} = N_{y,z}\cdot S\,$ has an open orbit in $\,G/P$, namely the union of the open Schubert cell traced out by $\,N\,$ and the codimension one Schubert cell which can be described symbolically by the equation $\,x = \infty\,$. In general, one  may not restrict a distribution vector $\,\tau \in V_{\lambda,\d}^{-\infty}\,$ from $\,G\,$ to the subgroup $\,S\,$. However, $\,p_{y,k}\circ p_{z,n}\, \tau\,$ and $\,p_{y,q}\circ p_{z,0}\, \tau\,$ can be restricted to any subgroup whose 
\pagebreak
orbit through the identity coset in $\,G/P$ is open -- in particular to $\,N_{y,z}\cdot S\,$. Since both $\,p_{y,k}\circ p_{z,n}\, \tau\,$ and $\,p_{y,q}\circ p_{z,0}\, \tau\,$ transform according to a character of $\,N_{y,z}$, they can be restricted to $\,S\,$ after all. This restriction transforms according to the character $\,\omega_{\lambda,\d}$ of $\,S\cap P$ on the right; cf.\ (\ref{chi}). At this point, (\ref{proj6}) and the distribution analogue of (\ref{wmuinf1}) imply
the first assertion of the corollary.

The identity (\ref{taucondition1}), with $x$ replaced by $x - \bar a/c$, equates two $G_\Z$-invariant vectors in $V_{\lambda,\d}^{-\infty}$, to which we can apply the projection $p_{y,k}\circ p_{z,0}$. After doing so, we substitute back $x + \bar a/c$ for $x$. The resulting equation extends the meaning of the equation following (\ref{taucondition1}), and the other equations derived from it, across $x = \infty$ and $x = 0$. In particular, a) and c) have meaning even at $x = \infty$ and $x = 0$, as equalities in $W_{\lambda_1 - \lambda_2,\d_3}$. Because of the $N_\Z$-invariance of $\tau$, the identity a) depends only on $\,\bar a$ modulo $c$, not on the integer $\bar a$. This, in effect, establishes the first and third of the identities in \corref{cor1}. The remaining two follow similarly from the proof of \propref{prop2}.

We should remark that the proof of parts a) and c) of \propref{prop2}, and of the first and third identities in
Corollary \ref{cor1}, depend only on the invariance of $\tau$ under the subgroup of $\G$ generated by
$N_\Z$ and the copy of ${\rm SL}(2,\Z)$ embedded as the top left $2\times 2$ block in ${\rm SL}(3,\Z)$,
whereas the other parts of the proposition and the corollary use invariance under $N_\Z$ and the copy of
${\rm SL}(2,\Z)$ embedded as the bottom right $2\times 2$ block.

\section{Proof of Proposition \ref{prop1}}\label{proofofprop1}

We shall deduce decomposition of $C^{-\infty}(N_\Z\backslash N)$ from 
the analogous decomposition of $L^2(N_\Z\backslash N)$. The $L^2$ statement we need can be deduced
from the results of Brezin \cite{Brezin:1970}. However, it is just as simple to establish it directly, using the
notion of automorphic distribution.

The discussion of   Section \ref{autodist} provides a canonical $N$-invariant inclusion $\Hom_N(V_n,
L^2(N_\Z \backslash N)) \hookrightarrow (V_{-n}^{-\infty})^{N_\Z}$ -- an isomorphism, in fact, since
$N_\Z$ is cocompact in $N$. We now construct an explicit basis of the space of $N_\Z$-invariant
distribution vectors. Any $\phi \in V_{-n}^{-\infty}$ is automatically $(0,0,1)$-invariant. Since $(1,0,0)$,
$(0,1,0)$, and $(0,0,1)$ generate $N_\Z$, the  $N_{\Z}$-invariance of $\phi \in V_{-n}^{-\infty} \simeq
\mathcal S'(\R)$ under the action $\pi_{-n}$ comes down to two conditions:
\begin{equation}
\label{invariance}
\phi( t + 1) \equiv \phi(t)\ \ \ \text{and} \ \ \ \textstyle e(-nt) \phi(t) \equiv \phi(t)\,.
\end{equation}
Since $\,e(-nt) - 1\,$ vanishes to first order on the set $\frac 1n \,\Z$ and nowhere else,
\begin{equation}
\label{basis}
\{\,\phi_{n,k}\,\mid \, k \in \Z/n\Z\,\}\,,\ \ \ \text{with} \ \ \ \phi_{n,k} \ =  \ \sum_{\ell \equiv k(\md n)} \, \d_{\frac {\ell}{n}}\,,
\end{equation}
constitutes a basis of $(V_{-n}^{-\infty})^{N_\Z} \subset V_{-n}^{-\infty} \simeq \mathcal S'(\R)$. In this formula $\delta_q$ denotes the delta function at $q \in \R$.

The embedding $i_{n,k} : V_{n}^\infty \hookrightarrow C^{\infty}(N_\Z \backslash N)$, which corresponds to the automorphic distribution $\tau = \phi_{n,k}$ via (\ref{newtau}), sends $f \in V_{n}^{\infty} \cong \mathcal S(\R)$ to
\begin{eqnarray}
\nonumber
(i_{n,k}f)(x,y,z) &=& \int_\R \phi_{n,k}(t)\, (\pi_{n}(x,y,z)f)(t) \, dt \\
&= & \int_\R \ \sum_{\ell \equiv k ( \md n)} \, \d_{\frac {\ell}{n}}(t)\, e(n(z + ty))f(t + x) \, dt\label{fs1} \\
&= &\sum_{\ell \equiv k ( \md n)} \, e(nz + \ell y)\, f(x + \textstyle\frac{\ell}{n}) \,.\nonumber
\end{eqnarray}
By construction, $i_{n,k}$ is $N$-equivariant with respect to the action $\,\pi_{n}$ on $V_{n}^\infty$ and
the right translation action on $C^{\infty}(N_\Z\backslash N)$.

We now describe alternate embeddings which are equivariant with respect to $\,\widehat\pi_{n}$. For the purposes of this argument, we let $\,\mathcal F$ denote the Fourier transform, normalized as in (\ref{fourierint}). The Poisson summation formula asserts that $\,\mathcal F (\sum_{\ell \in \Z} \,\d_\ell) = 
\sum_{\ell \in \Z} \,\d_\ell\,$; hence
\begin{equation}
\label{fs4}
(\mathcal F^{-1} \phi_{n,k})(t) \ = \mathcal F^{-1}\({\sum}_{\ell \in \Z}\, \d_{\frac{k}{n} + \ell}\)(t) \ = \ e({\textstyle{\frac{k}{n}}}t)\,{\sum}_{\ell \in \Z}\, \d_{\ell}(t)\,.
\end{equation}
According to (\ref{nonabelian1})--(\ref{intertwine}), when $f,\, h$ are related by the identity $h(t) =
(\mathcal F f)(nt)$,
\begin{equation}
\label{fs3}
\(\mathcal F\circ\pi_{n}(x,y,z)\,f\,\)(t) \ = \ \(\,\widehat\pi_{n}(x,y,z)\,h\,\)(t/n)\,.
\end{equation}
Combining (\ref{fs1})--(\ref{fs3}), we now find
\begin{eqnarray*}
(i_{n,k}f)(x,y,z)  &=& \int_\R \phi_{n,k}(t)\, (\pi_{n}(x,y,z)f)(t) \, dt
\\
&=& \int_\R (\mathcal F^{-1}\phi_{n,k})(t)\, (\mathcal F \circ \pi_{n}(x,y,z)f)(t) \, dt
\\
&=& \int_\R (\mathcal F^{-1}\phi_{n,k})(t)\, (\widehat\pi_{n}(x,y,z) h)(t/n) \, dt
\\
&=& \int_\R \ {\sum}_{\ell \in \Z} \ \d_\ell(t)\, \textstyle e(n(z-xy + {\textstyle{\frac{t}{n}}} x ) +
{\textstyle{\frac{tk}{n}}})\,h(\f{t}{n} - y) \, dt
\\
&=& {\sum}_{\ell \in \Z} \ \textstyle e(n(z-xy) +
\ell(x+{\textstyle{\frac{k}{n}}}))\,h({\textstyle\frac{\ell}{n}} - y) \,.
\end{eqnarray*}
This can be summarized as follows: 
\begin{equation}
\label{fs5}
(i_{n,k}f)(x,y,z)\ =  \ {\sum}_{\ell \in \Z/n\Z} \ \textstyle e({\textstyle\frac {\ell k}{n}})\, (j_{n,\ell}\,h)(x,y,z) \,,
\end{equation}
with $h \in \mathcal S(\R)$ related to $f \in \mathcal S(\R)$ by the equation
\begin{equation}
\label{fs6}
h(t) \ =  \ \widehat f(nt)\,,
\end{equation}
and with $\,j_{n,k} : V_{n}^\infty \hookrightarrow C^{\infty}(N_\Z \backslash N)$ defined by
\begin{equation}
\label{fs7}
(j_{n,k}h)(x,y,z)\ = \ {\sum}_{\ell \equiv k(\md n)}  \  \textstyle e(n(z-xy) + \ell x)\,h({\textstyle\frac{\ell}{n}} - y) \,.
\end{equation}
The embedding $j_{n,k}$ is $N$-equivariant with respect to the action $\widehat\pi_n$ on $V_n^\infty \cong \mathcal S(\R)$, as can be seen from the derivation of (\ref{fs5}) or by direct verification.

Recall the normalization (\ref{measure}) of Haar measure. For $k_1, k_2 \in \Z/n\Z$ and $f_1, f_2 \in V_{n}^\infty \ \simeq \ \mathcal S(\R)$,
\begin{equation}
\begin{aligned}
{}&(i_{n,k_1}f_1,i_{n,k_2}f_2)_{L^2(N_\Z\backslash N)}    \\
{}&\ \ \ = \ \displaystyle\int_0^1 \!\! \int_0^1 \!\! \int_0^1   {\sum}_{\stackrel{\scriptstyle{\ell_1 \equiv
k_1}}{\ell_2 \equiv k_2}}\ e((\ell_1 - \ell_2) y) f_1(x + {\textstyle\frac{\ell_1}{n}})\, \overline{f_2}(x +
{\textstyle\frac {\ell_2}{n}}) \, dx\, dy\, dz \label{fs8}\\  {}&\ \ \ = \ \begin{cases}
\ \ (f_1,f_2)_{L^2(\R)} \ \ \ &\text{if}\ k_1 = k_2\\
\ \ 0 \ \ \ &\text{if}\ k_1 \neq k_2\ .
\end{cases}
\end{aligned}
\end{equation}
It follows that the $i_{n,k}$, for $k \in \Z/n\Z$, extend continuously to $L^2(\R)$, and that the extensions constitute an orthonormal basis of $\Hom_N(V_{n},L^2(N_\Z \backslash N))$. Different values of $n \in \Z - \{0\}$ correspond to different central characters, so the images of $i_{n,k}$ for different $n$ are perpendicular. Thus, for $F \in L^2(N_\Z \backslash N)$, there exist uniquely determined $b_{r,s} \in \C$ and $f_{n,k} \in L^2(\R)$ such that
\begin{equation}
\label{fs9}
\begin{aligned}
F(x,y,z)\ &= \ {\sum}_{r,s \in \Z}\  b_{r,s} e(rx + sy)\, + \, {\sum}_{\stackrel{\scriptstyle{n \neq 0}}{k \in \Z /n \Z}}\  (i_{n,k}f_{n,k})(x,y,z)\,,
\\
\text{and}\ \ \ \|F\|^2_{L^2(N_\Z \backslash N)} & =  \ {\sum}_{r,s}\  |b_{r,s}|^2 \ + \  {\sum}_{k,n}\ 
\|f_{n,k}\|^2_{L^2(\R)} \,.
\end{aligned}
\end{equation}
Quite analogously there exist $h_{n,k} \in L^2(\R)$ such that
\begin{equation}
\label{fs10}
\begin{aligned}
F(x,y,z) \ &= \ {\sum}_{r,s \in \Z}\  b_{r,s} \,e(rx + sy)\, + \, {\sum}_{\stackrel{\scriptstyle{n \neq 0}}{k \in \Z /n\Z}} \  (j_{n,k}h_{n,k})(x,y,z)\,,
\\
\text{and}\ \ \ \|F\|^2_{L^2(N_\Z \backslash N)} &  =  \ {\sum}_{r,s}\  |b_{r,s}|^2 \ + \  {\sum}_{k,n}\ 
\|h_{n,k}\|^2_{L^2(\R)} \,.
\end{aligned}
\end{equation}
The identities (\ref{fs5}) to (\ref{fs7}) relate the $h_{n,k}$ to the $f_{n,k}$:
$$
h_{n,k}(y) \ =  \ {\sum}_{\ell \in \Z/n\Z} \ e({\textstyle\frac {k \ell}{n}})\,\widehat f_{n,\ell}(ny) \,,
$$
or in equivalent, but more symmetric form,
\begin{equation}
\label{fs11}
{\sum}_{k \in \Z/n\Z}\ a_k\,h_{n,k}(y) \ =  \ {\sum}_{k \in \Z/n\Z}\  \widehat a_k\, \widehat f_{n,k}(ny) \,,
\end{equation}
with any choice of coefficients $a_k$, $k \in \Z/n\Z$;  cf.\ (\ref{fourierint3}).

When the explicit formulas (\ref{fs1},\ref{fs7}) for $i_{n,k}$ and $j_{n,k}$ are substituted,
(\ref{fs9})--(\ref{fs11}) amount to the $L^2$ analogue of Proposition \ref{prop1}. Let us suppose now that
the function $F$ in (\ref{fs9}) has derivatives of order up to $\ell \geq 1$. 
\pagebreak
In view of the equivariance of
$i_{n,k}$ and the identities (\ref{partials}), the non-abelian Fourier components of
$r({\textstyle\frac{\partial\ }{\partial x}})F$, $r({\textstyle\frac{\partial\ }{\partial y}})F$,
$r({\textstyle\frac{\partial\ }{\partial z}})F$ are
\begin{equation}\label{fs12}
\begin{aligned} 
r({\textstyle\frac{\partial\ }{\partial x}})\,i_{n,k}(f_{n,k}(t))   &=   i_{n,k}(f_{n,k}'(t))\,,
\\
r({\textstyle\frac{\partial\ }{\partial y}})i_{n,k}(f_{n,k}(t)) &=  2\,\pi \, i\, n\ i_{n,k}(tf_{n,k}(t))\,,
\\
r({\textstyle\frac{\partial\ }{\partial z}})i_{n,k}(f_{n,k}(t))   &=  2\,\pi \, i\, n\  i_{n,k}(f_{n,k}(t))\,.
\end{aligned}
\end{equation}
In these equations, $r$ denotes the right regular representation (\ref{rightregular}), and ${\textstyle\frac{\partial\ }{\partial x}}$, ${\textstyle\frac{\partial\ }{\partial y}}$, ${\textstyle\frac{\partial\ }{\partial z}}$ are regarded as generators of the Lie algebra of $N$, as in (\ref{partials}). In particular all the $f_{n,k}'(t)$ and $tf_{n,k}(t)$ must lie in $L^2(\R)$. Iterating this argument, one finds that the existence of derivatives up to order $\ell$ implies the finiteness\begin{footnote}{Note the presence of the factor $n$ in the second formula in (\ref{fs12}) and recall that $n \neq 0$.}\end{footnote} of
\begin{equation}
\label{fs13}
\|F\|_\ell^2 \, = \, 
\sum_{\stackrel{\scriptstyle{i,j \geq 0}}{i + j \leq \ell}} \(\sum_{r,s} \, r^{2i}\,s^{2j}\, |b_{r,s}|^2 \, + \,
\sum_{n,k}\,n^{2\ell-2i}\,\|t^j\,f_{n,k}^{(i)}(t)\|_{L^2(\R)}^2 \! \).
\end{equation}
Conversely the finiteness of $\|F\|_\ell^2$ ensures the existence of $L^2$ derivatives of order up to $\ell$. We conclude that the family of Sobolev norms $\,\|\ \|_\ell\,$, {\nobreak${\ell \in \N}$}, defines the topology of $C^{\infty}(N_\Z \backslash N)$. In particular, each non-abelian Fourier component $f_{n,k}$ of a $C^\infty$ function $F$ satisfies the finiteness conditions
\begin{equation}
\label{fs14}
{\sum}_{0 \leq i,j\,,\,i + j \leq \ell}\ \|t^j f_{n,k}^{(i)}(t)\|_{L^2(\R)}^2  \ < \ \infty
\end{equation}
for all $\ell > 0$; hence $f_{n,k} \in \mathcal S(\R)$.

Since the norms $\,\|\ \|_\ell\,$ increase with $\ell$, for any $\tau \in C^{-\infty}(N_\Z \backslash N)$, the linear map
\begin{equation}
\label{fs15}
C^{\infty}(N_\Z \backslash N) \ \ni \ F \longmapsto \ \int_{N_\Z \backslash N}\, \tau \, F\, dn
\end{equation}
must be bounded with respect to at least one of the norms, and therefore with respect to all but finitely many of them. Thus $\tau$ lies in the completion of $C^{\infty}(N_\Z \backslash N)$ with respect to a Sobolev norm $\,\|\ \|_{-\ell}\,$ with negative index, i.e. the norm dual to the norm $\,\|\ \|_\ell\,$, with $\ell \geq 0$. The existence of a series expansion as described in Proposition \ref{prop1} now follows from an essentially formal Hilbert space argument. Each of the non-abelian Fourier components $\sigma_{n,k}$ of $\tau$ pairs continuously against the component $f_{-n,-k} \in \mathcal S(\R)$ of an arbitrary test function {\nobreak${F \in C^{\infty}(N_\Z \backslash N)}$}; in other words, $\sigma_{n,k}$ is a tempered distribution. By construction, the series of $\tau$ converges with respect to the norm $\,\|\ \|_{-\ell}\,$, which means that
it converges in the strong distribution topology. The relation (\ref{fs11}), finally, is inherited by the
non-abelian Fourier
  components of $\tau$ because $\mathcal S(\R)$ is dense in $\mathcal S'(\R)$. This completes the proof.

\section{Voronoi summation for
${\rm GL}(3,\Z)$}\label{voronsec}

In this section we prove our main theorem using the machinery developed in Sections \ref{autodist} and
\ref{heis}, and the analytic tools developed in \cite{MilSch:2003b}. We continue with the hypotheses of
Section \ref{heis}; in particular the automorphic distribution $\tau$ is invariant under $\G={\rm GL}(3,\Z)$,
and $\sum \d_j = 0$; cf.\ (\ref{delta3}). We also suppose $\tau$ is cuspidal, so that 
\begin{equation}
\label{cuspc}
c_{r,0}\ = \ c_{0,s}\ = \ 0\ \ \ \text{for all}\ \ r,\, s \in \Z\,,
\end{equation}
as follows from Lemma \ref{lem1}.

In the paper \cite{MilSch:2003b} we introduce the notion of vanishing to infinite order for distributions: $\sigma \in C^{-\infty}(\R)$ vanishes to infinite order at $x=0$ if, for each $n \in \N$, there exists an open interval $I$ containing the origin, a collection of $L^\infty$ functions $f_j \in L^\infty(I)$, $1 \leq j \leq N_n$, and nonnegative integers $k_j\,$, $1 \leq j \leq N_n$, such that on the interval $I$,
\begin{equation}
\label{inf_order1}
\sigma(x) \ = \ x^n\,\textstyle{\sum}_{1\leq j \leq N_n}\, x^{k_j}\, \frac{d^{k_j}\ }{dx^{k_j}}\,f_j(x)\,.
\end{equation}
We show that a distribution $\sigma$ which satisfies this condition is uniquely determined, among all distributions with the same property, by its restriction to 
$\R^*$ \cite[Lemma 2.8]{MilSch:2003b}. This justifies the following terminology: a distribution $\sigma
\in C^{-\infty}(\R^*)$ has a canonical extension across $0$ if there exists a -- necessarily unique -- extension
to a distribution on $\,\R\,$ which vanishes to infinite order at the origin. By definition, $\sigma \in
C^{-\infty}(\R)$ vanishes to infinite order at $x_0 \in \R$ if $x \mapsto \sigma(x+x_0)$ vanishes to infinite
order at $0$; $\sigma \in C^{-\infty}(\R - \{x_0\})$ has a canonical extension across $x_0$ if $x \mapsto
\sigma(x+x_0)$ has a canonical extension across $0$; and $\sigma \in C^{-\infty}(\R)$ extends canonically
across $\infty$ if $x \mapsto \sigma(1/x)$ extends canonically across $0$ \cite{MilSch:2003b}. Very
importantly, if $\sigma$ has a canonical extension across either $0$ or $\infty$, 
 then so do the distributions $x \mapsto |x|^\a \sigma(x)$ and $x \mapsto |x|^\a \sg(x) \sigma(x)$, for every $\a \in \C$
\cite[Prop.~2.26]{MilSch:2003b}.

The assertion of the next lemma depends crucially on the cuspidality of the automorphic distribution $\tau$, which is a standing assumption in this section. We shall prove it at the end of the section, after using it to complete the proof of the Voronoi summation formula.

\begin{lem}\label{inforder} 
The distributions $\,\sum_{r\in \Z}\, c_{r,q}\, e(rx)  \,,\, \sigma_{n,k} \,,\, \sum_{s\in \Z}\, c_{q,s}\, e(sy)${\rm
,}\break
$\rho_{n,k}$ extend canonically across $\infty${\rm ,} and all of them vanish to infinite order at every
rational point.
\end{lem}

A comment on the connection between the lemma and \corref{cor1}: the former considerably strengthens the latter, from which it will be deduced. The corollary asserts the equality of unspecified extensions of distributions to distribution vectors, whereas the 
lemma implies the equality of the uniquely determined canonical extensions.

We begin the proof of the main theorem with some remarks about the operators $T_{\a,\d}$ which are the subject of
\cite[\S 6]{MilSch:2003b}: for $\a \in \C$, and $\d \in \Z/2\Z$,
\begin{equation}
\label{phixa1}
T_{\a,\d}f\ = \ \mathcal F\(\, x \mapsto f(1/x)\,(\sg x)^\d\,|x|^{-\a -1}\,\) \qquad(\,f \in \mathcal S(\R))\,.
\end{equation}
The integral that computes this Fourier transform converges absolutely when $\,\operatorname{Re} \alpha > 0\,$. Even without the restriction on $\a$, $T_{\a,\d}f$ is
well-defined as a function on $\R - \{0\}$ which depends holomorphically on $\alpha$, because the Fourier
kernel $x \mapsto e(-xy)$ vanishes to infinite order at infinity when $y \neq 0$. The function $T_{\a,\d}f$
tends to zero rapidly as $|x|\to \infty$, along with all its derivatives, and at the origin, $T_{\a,\d}f$ has
potential singularities, which can be described explicitly. Let $\Sch_{\text{sis}}(\R)$ denote the linear span
of all products $\,(\sg x)^\eta (\log |x|)^k |x|^\beta f(x)$, with $\beta \in \C$, $\eta\in\Z/2\Z$, $k\geq 0$, and
$f\in\Sch(\R)$. Then $T_{\a,\d}$ extends from $\Sch(\R)$ to $\Sch_{\text{sis}}(\R)$ and maps this space to
itself:
\begin{equation}
\label{phixa2}
T_{\a,\d} \, :\, \Sch_{\text{sis}}(\R) \ \longrightarrow \ \Sch_{\text{sis}}(\R)\,,
\end{equation}
\cite[Th.~6.6]{MilSch:2003b}. In particular, one can compose any two or more of these operators. Also
note that the integration pairing
\begin{equation}
\label{phixa3}
\begin{array}{l}
\displaystyle{ \int_\R f(x)\,\sigma(x)\,dx\ \ \ \text{is well-defined for any $f\in \Sch_{\text{sis}}(\R)$},}
\\
\qquad \text{
 provided $\sigma\in\Sch'(\R)$ vanishes to infinite order at the origin,}
\end{array}
\end{equation}
as follows from the definition of the space $\Sch_{\text{sis}}(\R)$.

If $\,\sigma \in \Sch'(\R)$ vanishes to infinite order at the origin, the Fourier transform $\,\widehat\sigma$ extends canonically across $\infty$ 
\cite[Th.~3.19]{MilSch:2003b}. Hence $\widehat\sigma(1/x)$ may be regarded as a tempered
distribution on $\,\R$ which vanishes to infinite order at $x=0$. The distribution $(\sg x)^\delta
|x|^{\alpha-1}\widehat\sigma(1/x)$ inherits this property, so that 
\begin{equation}
\label{phixa4}
T^*_{\a,\d}\,\sigma(x)\ =\ (\sg x)^\delta \, |x|^{\alpha-1}\, \widehat\sigma(1/x)
\end{equation}
defines a linear operator on the space of $\,\sigma\in\Sch'(\R)$ which vanishes to infinite order at the origin.
This is the adjoint of the operator (\ref{phixa2}) with respect to the pairing (\ref{phixa3}):
\begin{equation}
\label{phixa5}
\begin{aligned}
\int_\R T_{\a,\d}f(x)\,\sigma(x)\,dx \ = \ \int_\R f(x)\, T_{\a,\d}^*\, \sigma(x)\, dx\ \ \ \text{if $f\in \Sch_{\text{sis}}(\R)$}
\\
\text{and if $\sigma\in\Sch'(\R)$ vanishes to infinite order at $x=0$}
\end{aligned}
\end{equation}
\cite[Th.~6.9]{MilSch:2003b}.

Turning to the substance of the proof, we express the automorphic distri\-bution $\tau$ as in Proposition \ref{prop1}. It will be convenient to work with renormalized Fourier coefficients
\begin{equation}
\label{cnorm}
a_{r,s}\ =\ c_{r,s}\,|r|^{\l_1}\,|s|^{-\l_3}\,(\sg r)^{\d_1}\,(\sg s)^{\d_3}\,. 
\end{equation} \pagebreak
These can be described in terms of the Hecke action -- see (\ref{a1pp1})--(\ref{mudef}) below -- and
therefore have canonical meaning. The passage from the $\,c_{r,s}$ to the $\,a_{r,s}$ is analogous to the
normalization $a_n|n|^{-\nu}$ for the Fourier coefficients in (\ref{gl2taudef}). We should also point out that
the parameters $(\l,\d)$ in (\ref{embed3}) depend on the choice of a Casselman embedding, which is not
unique. Different embeddings give different Fourier coefficients $\,c_{r,s}\,$. On the other hand, the
description of the $\,a_{r,s}$ in terms of the Hecke action means that they are determined by $\tau$ itself,
except for a constant normalizing factor. Note that
\begin{equation}
\label{awitheps}
a_{r,s}\ = \ a_{-r,s}\ = \ a_{r,-s}\ = \ a_{-r,-s}\,,\ \ \ \ a_{r,0}\ = \ 0\,,\ \ \ \ a_{0,s}\ = \ 0\,,
\end{equation}
because of Lemma \ref{lem2}, the renormalization (\ref{cnorm}), and the cuspidality of $\tau$.

Fix integers $a,\, c,\, q\in\Z$ such that $a$, $c$ are relatively prime, $c \neq 0$, and $q > 0$. As usual, we let $\bar a\in\Z$ denote an inverse of $a$ modulo $c$; for emphasis,
\begin{equation}
\label{vora}
a,\, \bar a,\, c,\, q \in \Z\,, \ \ c \neq 0\,,  \ \ q>0\,, \ \ (a,c) \equiv 1 \ (\operatorname{mod} c)\,, \ \ a \bar a \equiv 1 \ (\operatorname{mod} c)\,.
\end{equation}
We now substitute $y = c^{-2}\,x$ into the second equation in \propref{prop2}, reverse the roles of $a$, $\bar a$, and use (\ref{cnorm}):
\begin{eqnarray}
&&\label{vorb}\\[-8pt] \nonumber
 \rho_{cq,\bar aq}(x/c^2)& =& (\sg cx)^{\d_1}\, \left|\f xc \right|^{\l_2 - \l_3 - 1} \sum_{n \neq 0} \ c_{q,n} \,
e( n/x -  n\, a / c)
\\ 
& =&  (\sg cx)^{\d_1} q^{-\l_1} \left|\f xc \right|^{\l_2 - \l_3 - 1} \sum_{n \neq 0} \, a_{q,n} \,(\sg n)^{\d_3}
|n|^{\l_3} e( n/x  -  n\, a / c)\,.
\nonumber
\end{eqnarray}
According to \propref{prop1},
\begin{equation}
\label{vorc}
\rho_{cq,\bar aq}(x/c^2) \ =  \ {\sum}_{\ell \in \Z/cq\Z}\ e(\ell \bar a/c)\,\widehat{\sigma}_{cq,\ell}(qx/c)\,.
\end{equation}
In view of \lemref{inforder} and the definition of the operator $\,T^*_{\alpha,\delta}\,$,
(\ref{vorb}), (\ref{vorc}) imply
\begin{equation}
\label{vord}
\begin{aligned}
{\sum}_{n \neq 0} \, a_{q,n} \,(\sg n)^{\d_3} |n|^{\l_3} e( n \, x  -  n\, a / c)\ = \ (T_{\l_2-\l_3,\d_1}^*\sigma_0)(x)\,, \ \ \text{with}
\\
\sigma_0(x)\ = \ q^{\l_1 - 1}\,(\sg c)^{\d_1}\, |c|^{\l_2-\l_3}\,{\sum}_{\ell \in \Z/cq\Z}\ e(\ell \bar a/c)\, \sigma_{cq,\ell}(cx/q)\,.
\end{aligned}
\end{equation}
For $\ell \in \Z/cq\Z$, we let $\bar\ell\in\Z$ denote a solution of the congruence $\bar \ell \ell \equiv (cq,\ell)$ modulo $cq\,$; the particular choice of $\bar\ell$ will not matter. With this convention,
\begin{equation} \label{vore}
\begin{aligned}
{}&  \sigma_{cq,\ell}(cx/q) 
 \\
{}&\quad   = \ (\sg\textstyle \frac{c^2 x}{(cq,\ell)})^{\d_3}\, |\textstyle \frac{(cq,\ell)}{c^2
x}|^{1-\l_1 +
\l_2}\,\displaystyle{\sum_{n\neq 0}}\,c_{n,\,(cq,\ell)}\, \textstyle e(\frac{n(cq,\ell)\bar \ell}{cq}  -
\frac{n(cq,\ell)^2}{c^3qx}) 
\\
{}& \quad   = \ {\sum}_{n\neq 0}\,\frac{(\sg x)^{\d_3}|(cq,\ell)|^{1 - \l_1 + \l_2 + \l_3}}{(\sg
n)^{\delta_1}|n|^{\lambda_1}|c^2x|^{1-\l_1 + \l_2}}\,a_{n,\,(cq,\ell)}\, \textstyle e(\frac{n(cq,\ell)\bar
\ell}{cq}  - \frac{n(cq,\ell)^2}{c^3qx})\,, 
\end{aligned}
\end{equation}
as follows from the first identity in \propref{prop2} and the definition of the $a_{r,s}$ in terms of the $c_{r,s}$. This implies
\begin{eqnarray}&&\label{vorf}
\\ 
\nonumber
 \quad \sigma_{cq,\ell}(cx/q)& =  &(T_{\l_1-\l_2,\d_3}^*\sigma_1)(x)\,, \ \ \text{with}
\\ 
 \sigma_1(x) &=  &{\sum}_{n\neq 0}\,\frac{|(cq,\ell)|^{1 - \l_1 + \l_2 + \l_3}}{(\sg
n)^{\delta_1}|n|^{\lambda_1}c^{2 - 2\l_1 + 2\l_2}}\,a_{n,\,(cq,\ell)}\, \textstyle e(\frac{n(cq,\ell)\bar
\ell}{cq})\  \d_{\frac{n(cq,\ell)^2}{c^3q}}(x), \nonumber
\end{eqnarray}
which does vanish to infinite order at the origin, as required.

We now fix a test function $\,f\in (\sg x)^{\delta_3}|x|^{\lambda_3}\Sch(\R)$, as in our statement of the Voronoi summation formula. Then
\begin{equation}
\label{vorg}
g(x)\ = \ \int_\R (\sg y)^{\delta_3}\, |y|^{-\lambda_3}\, f(y)\, e(-xy)\, dy
\end{equation}
is the Fourier transform of a Schwartz function, and hence is a Schwartz function itself. In view of Fourier inversion, (\ref{phixa5}), and (\ref{vord}),
\begin{equation}
\label{vorh}
\begin{array}{l}
\displaystyle{ {\sum}_{n\neq 0}\ a_{q,n}\,e(-na/c)\,f(n)\ = \ {\sum}_{n\neq 0}\ \frac{|n|^{\lambda_3}\,a_{q,n}}{(\sg
n)^{\delta_3}}\,e(-na/c)\, \widehat g(-n)}
\\
\qquad  \displaystyle{= \ \int_\R {\sum}_{n\neq 0}\  a_{q,n}\,|n|^{\lambda_3}(\sg n)^{\delta_3}\,e(nx-na/c)\,g(x)\,dx
}\\
\qquad  \displaystyle{= \  {\sum}_{\ell \in \Z/cq\Z}\ \frac{(\sg c)^{\d_1}\,e(\ell \bar a/c)}{q^{1 - \l_1}\,|c|^{-\l_2 +
\l_3}}\int_\R
\sigma_{cq,\ell}(cx/q)(T_{\l_2-\l_3,\d_1}g)(x)\,dx\,.
}\end{array}
\end{equation}
Similarly (\ref{phixa5}) and (\ref{vorf}) imply
\begin{equation}
\label{vori}
\begin{aligned}
\displaystyle{\int_\R \sigma_{cq,\ell}(cx/q)\,h(x)\,dx\ = \
 {\sum}_{n\neq 0}\ \frac{|(cq,\ell)|^{1 - \l_1 + \l_2 + \l_3}}{(\sg n)^{\delta_1}|n|^{\lambda_1}c^{2 - 2\l_1 +
2\l_2}}  }\\
\displaystyle{\times \ \ a_{n,\,(cq,\ell)}\, e(\textstyle\frac{n(cq,\ell)\bar \ell}{cq})\, (T_{\l_1-\l_2,\d_3}h)(\textstyle{\frac{n(cq,\ell)^2}{c^3q}})\,,
}\end{aligned}
\end{equation}
for any $\,h\in\Sch_{\text{sis}}(\R)$. We now substitute $h=T_{\l_2-\l_3,\d_1}g$ and combine the resulting equation with (\ref{vorh}):
\begin{equation}
\label{vorj}
\begin{array}{l}
\displaystyle{ {\sum}_{n\neq 0}\ a_{q,n}\,e(-na/c)\,f(n) }
\\
\displaystyle{ \quad = \ {\sum}_{n\neq 0}\displaystyle{\sum}_{\ell \in
\Z/cq\Z}\,\f{|c|}{|(cq,\ell)|}\,(\sg\textstyle\f{n(cq,\ell)^2}{c^3q})^{\d_1}\,
|\textstyle\f{n(cq,\ell)^2}{c^3q}|^{1-\l_1} }
\\
\qquad\displaystyle{ \times\  \f{a_{n,\,(cq,\ell)}}{|n|}\, \textstyle e(\frac{n(cq,\ell)\bar \ell}{cq} + \f{\ell \bar
a}{c})\,(T_{\l_1-\l_2,\d_3}\circ T_{\l_2-\l_3,\d_1}\,g)(\textstyle\frac{n(cq,\ell)^2}{c^3q})\,,}
\end{array}
\end{equation}
since $q>0$ by assumption.

As the index of summation $\ell$ in (\ref{vorj}) ranges over $\,\Z/cq\Z\,$, the greatest common divisor $d =_{\rm def} (cq,\ell)$ ranges over all divisors of $cq$ and the quotient $\ell' =_{\rm def} \ell/d$ over the units modulo $cq/d$. In this situation, any inverse $\bar\ell'$ to $\ell'$ modulo $cq/d$ can take the place of $\bar \ell$. When the sum over $\ell$ is broken up into a double sum over $d$ and $\ell'$, only the exponential term involves $\ell'$ and $\bar\ell'$.
  This observation leads to the following simplification of (\ref{vorj}):
\begin{equation}
\label{vork}
\begin{array}{l}
 \displaystyle{{\sum}_{n\neq 0}\ a_{q,n}\,e(-na/c)\,f(n) }
\\
\displaystyle{\quad= \ {\sum}_{n\neq 0}\displaystyle{\sum}_{d | cq}\ \f{|c|}{|n\,
d|}\,(\sg\textstyle\f{nd^2}{c^3q})^{\d_1}\, |\textstyle\f{nd^2}{c^3q}|^{1-\l_1} }
\\
 \qquad  \times\  \displaystyle{a_{n,d}\, S(\bar a q, n ;cq/d)\,(T_{\l_1-\l_2,\d_3}\circ
T_{\l_2-\l_3,\d_1}\,g)(\textstyle\frac{nd^2}{c^3q})}
\\
\displaystyle{\quad =\ {\sum}_{n\neq 0}\ {\sum}_{d|cq}\ \frac{|c|}{|n\,d|}\ a_{n,d}\ S(\bar a q, n;
cq/d)\,F(\textstyle\frac{nd^2}{c^3q})\,,}
\end{array}
\end{equation}
where by definition
\begin{equation}
\label{vorl}
F(t)  =  (\sg t)^{\d_1} |t|^{1 - \l_1}\!\( T_{\l_1 - \l_2,\d_3}\!\circ T_{\l_2 - \l_3, \d_1}\!\circ \mathcal F\!\((\sg x)^{\d_3}|x|^{-\l_3}f\)\)(t) \,,
\end{equation}
and where
\begin{equation}
\label{vorm}
{\sum}_{\ell' \in (\Z/cqd^{-1}\Z)^*} \, e(\textstyle\f{nd\bar\ell'}{cq} + \f{\ell' d\bar aq}{cq})\ \ = \ \ S(\bar a q, n ;cq/d)
\end{equation}
is the Kloosterman sum with parameters $\,\bar a q\,$, $\,n\,$, $\,cq/d\,$.

In the definition (\ref{vorl}) of $F$, the Fourier transform is computed by the absolutely convergent integral (\ref{vorg}). Since $g$ is a Schwartz function, the integral expressing $T_{\l_2-\l_3,\d_1}\,g$ in terms of $g$ converges absolutely, provided $\operatorname{Re}(\l_2 - \l_3)>0$. In that case $h(x)=T_{\l_2-\l_3,\d_1}g(x)$ is globally continuous and decays rapidly as $|x|\to\infty$. But then $\,T_{\l_1 - \l_2,\d_3}h\,$ is also computed by an absolutely convergent integral, provided $\,\operatorname{Re}(\l_1 - \l_2)>0$. We conclude: if $\,\operatorname{Re}\l_1 > \operatorname{Re}\l_2 > \operatorname{Re}\l_3\,$,
\begin{eqnarray}
&&\label{vorn}\\[-4pt] \nonumber
 F(t) &=& (\sg t)^{\d_1} |t|^{1 - \l_1}\!\!\int_{z=-\infty}^\infty \int_{y=-\infty}^\infty \int_{x=-\infty}^\infty
\!\!\!\! f(x)e(- x/y - y/z -tz) 
\\ & &\times\ (\sg x)^{\d_3}(\sg y)^{\d_1}(\sg z)^{\d_3}|x|^{-\l_3}|y|^{-\l_2+\l_3-1}|z|^{-\l_1+\l_2-1}
\,dx\,dy\,dz\,,
\nonumber
\end{eqnarray}
with all three integrals converging absolutely when performed in the indicated order. Since the operator (\ref{phixa1}) depends holomorphically on $\alpha$, the identity (\ref{vorn}) retains meaning for other values of $\l_1$, $\l_2$, $\l_3$ by analytic continuation. The same is true for the repeated integral
\begin{equation}
\label{voro}
\begin{aligned} 
{} F(t)   =& \displaystyle{ \int_{x_1=-\infty}^\infty\int_{x_2=-\infty}^\infty\int_{x_3=-\infty}^\infty\!
f\(\f{x_1 x_2 x_3}{t}\) }
\\[6pt]
 &\displaystyle{\times \ \ {\prod}_{j=1}^3\,\((\sg x_j)^{\d_j}\, |x_j|^{-\l_j}\,e(-x_j)\)
dx_3\,dx_2\,dx_1\,,}
\end{aligned}
\end{equation}
which results from (\ref{vorn}) by the change of variables $x_1 = tz$, $x_2 = y/z$, $x_3 = x/y$.

We now impose the parity condition $f(-x)=(-1)^\eta f(x)$, with $\eta \in \{0,1\}$, 
\pagebreak
in which case $F$ also satisfies this condition. The intermediate function $g$, defined in (\ref{vorg}), has parity $\eta + \delta_3$, and
\begin{equation}
\label{vorp}
M_{\eta + \delta_3}g \,(s) \ = \ (-1)^{\eta + \delta_3}\, G_{\eta + \delta_3}(s)\, M_\eta f\,(1-s-\lambda_3)\,,
\end{equation}
as follows from the formula \cite[4.51]{MilSch:2003b} for the Fourier transform of the signed Mellin
kernel; here $G_{\delta_3+\eta}(s)$ denotes the Gamma factor defined in (\ref{whyg}), (\ref{heresG1}). At
this point, two applications of \cite[Lemma 6.19]{MilSch:2003b} imply
\begin{equation}
\label{vorq}
M_\eta F\,(s) \ = \ (-1)^\eta \,\( {\prod}_{j=1}^3\ G_{\eta + \delta_j}(s - \lambda_j + 1)\) M_\eta f\,(-s)\,.
\end{equation}
We use the Gamma identities (\ref{gammarc}) and (\ref{heresG}) to express the product of Gamma factors in the more familiar form
\begin{equation}
\label{vorr}
{\prod}_{j=1}^3\ G_{\eta + \delta_j}(s - \lambda_j + 1) \ = \ \pi^{-3/2- 3s} \  {\prod}_{j=1}^3   \bigg( i^{\d_j'} \pi^{\l_j} \, \f{\G(\f{s+1-\l_j+\d_j'}{2})}{\G(\f{-s+\l_j+\d_j'}{2})} \bigg),
\end{equation}
with $\d_j' \in \{0,1\}$ determined by the congruence $\delta_j' \equiv \eta + \delta_j\ (\operatorname{mod}2$). An analysis of the poles of the Gamma factors and of $M_\eta f$ shows that $M_\eta F$ is meromorphic and regular for $\operatorname{Re}s > \operatorname{max}_j(\operatorname{Re} \l_j - 1)$. Moreover $M_\eta F\,(s)$ decays rapidly on vertical strips
\cite[\S 6]{MilSch:2003b}. Thus we can apply the signed Mellin inversion formula,
\begin{equation}
\label{vors}
F(x) \ \ = \  \  \f{(\sg{x})^\eta}{4\pi i} \int_{\operatorname{Re}s\,=\,\sigma} (M_\eta
F)(s)\,|x|^{-s}\,ds\,,
\end{equation}
with $\sigma$ chosen large enough to place the line of integration to the right of all the poles of $\,M_\eta F$. The alternate expression for $F$ in our statement of the Voronoi summation formula now follows from
 (\ref{vorq})--(\ref{vors}).

The singularities of $F$ at the origin reflect the location and order of the poles of $\,M_\eta F$. For the statement of the relevant results in \cite{MilSch:2003b}, we introduce the partial order
\begin{equation}
\label{vort}
(\alpha_1,\eta_1)\ \preceq \ (\alpha_2,\eta_2)\ \ \Longleftrightarrow \ \ \alpha_2 - \alpha_1 \in (2\Z + \eta_1 + \eta_2) \cap \Z_{\geq 0}
\end{equation}
on the set $\C \times \Z/2\Z$. According to \cite[6.54]{MilSch:2003b},
\begin{equation}
\label{voru}
F \, \in \, {\sum}_{j=1}^3\, (\sg x)^{\d_j}\, |x|^{1 - \l_j} \,\Sch(\R)\ \ \ \text{if}\ \ (\alpha_i,\eta_i) \npreceq  (\alpha_j,\eta_j)\ \ \text{for all}\ i\neq j\,.
\end{equation}
In the description of the remaining two cases, $\sigma$ denotes a permutation of the set $\{1,2,3\}$ such that $\operatorname{Re}\lambda_{\sigma_1} \leq \operatorname{Re}\lambda_{\sigma_2} \leq\operatorname{Re}\lambda_{\sigma_3}$. With this convention,
\begin{equation}
\label{vorv}
\begin{aligned}
F  \in  (\sg x)^{\d_{\sigma 1}}|x|^{1-\l_{\sigma 1}}\log{|x|}\,\Sch(\R)\, + \,  {\sum}_{j = 2,\, 3}\, (\sg x)^{\d_{\sigma j}}|x|^{1-\l_{\sigma j}}\Sch(\R)
\\
\text{if}\ \ (\l_{\sigma 1},\d_{\sigma 1}) \preceq (\l_{\sigma 2},\d_{\sigma 2})\ \ \text{and}\ \ (\l_{\sigma j},\d_{\sigma j}) \npreceq (\l_{\sigma 3},\d_{\sigma 3})\ \ \text{for}\ \ j=1,\,2
\end{aligned}
\end{equation}
\cite[6.55]{MilSch:2003b}, and
\begin{equation}
\label{vorw}
\begin{aligned}
F \ \in \ {\sum}_{1\leq j\leq 3}\ (\sg x)^{\d_{\sigma j}}\,|x|^{1-\l_{\sigma j}}\,(\log{|x|})^{3-j}\,\Sch(\R)
\phantom{overawholelot}\\
\phantom{overawholelot}\text{if}\ \ (\l_{\sigma 1},\d_{\sigma 1}) \preceq (\l_{\sigma
2},\d_{\sigma 2})
\preceq (\l_{\sigma 3},\d_{\sigma 3})
\end{aligned}
\end{equation}
\cite[6.56]{MilSch:2003b}. That completes the proof of the Voronoi summation formula for ${\rm
GL}(3)$, except for Lemma \ref{inforder}.

\demo{Proof  of Lemma {\rm \ref{inforder}}}  Recall the arguments and notation in the proof of
\corref{cor1}. Restated in terms of this notation, the first of the two cuspidality conditions (\ref{cuspidality})
asserts that
$p_{x,0}\circ p_{z,0}\,\tau = 0$. As pointed out in the proof of \corref{cor1}, $p_{x,0}$ and $p_{y,q}$
commute on the image of $p_{z,0}$; hence
\begin{equation}
\label{inforder1}
0 \ = \ p_{y,q}\circ p_{x,0}\circ p_{z,0}\,\tau \ = \ p_{x,0}\circ p_{y,q}\circ p_{z,0}\,\tau\,.
\end{equation}
We   argued that it is legitimate to restrict the distribution $p_{y,q}\circ p_{z,0}\,\tau$ to the subgroup $S
\cong {\rm SL}(2,\R)$, and we   identified this restriction with the extension of $\sum_r c_{r,q}\,e(rx)$ to
a vector in $W_{\l_1-\l_2,\d_3}^{-\infty}$ which is periodic, even at $\infty$. To simplify the notation, we
tacitly identify the sum with its extension. Then (\ref{inforder1}) is equivalent to
\begin{equation}
\label{inforder2}
\int_0^1\,\psi_{\l_1-\l_2,\d_3}\begin{pmatrix}1 & t \\ 0 & 1\end{pmatrix}\({\sum}_{r\in \Z}\ c_{r,q}\,e(rx)\) dt \ \ = \ \ 0\,.
\end{equation}
This in turn implies $c_{0,q} = 0$, as was argued already in the proof of \lemref{lem1}, and also has a consequence at infinity, which we shall now explicate.

The periodic distribution $\sum_r c_{r,q}\,e(rx)$ has zero constant term, thus extends canonically across $\infty$ by
\cite[Prop.~2.19]{MilSch:2003b}. At first glance, there are two notions of extension across $\infty\,$:
as a  distribution on $\R \cup \{\infty\}$, and as a vector in the representation space
$W_{\l_1-\l_2,\d_3}^{-\infty}$. The two extensions are related by the multiplicative factor $(\sg
x)^{\d_3}|x|^{\l_1-\l_2}$ -- see the discussion following (\ref{wmudist}) -- which does not affect the notion
of vanishing to infinite order at $\infty$ \cite[Prop.~2.26]{MilSch:2003b}. It follows that $\sum_r
c_{r,q}\,e(rx)$ can be extended to a vector in $W_{\l_1-\l_2,\d_3}^{-\infty}$ not only in the manner
described in \corref{cor1} -- we shall call this the {\it natural extension} -- but also as a distribution vector
which va\-nishes to infinite order at~$\infty$. We let $(\sum_r c_{r,q}\,e(rx))_\infty \in
W_{\l_1-\l_2,\d_3}^{-\infty}$ denote the difference between the two, which is supported at infinity by
construction. In view of \cite[Lemma 2.8]{MilSch:2003b} and \cite[Prop.~7.20]{MilSch:2003b}, the
canonical extension inherits both periodicity and   property (\ref{inforder2}) from the natural extension;
hence so does their difference:
\begin{equation}
\label{inforder3}
\begin{aligned}
 \psi_{\l_1-\l_2,\d_3}\begin{pmatrix}1 & 1 \\ 0 & 1\end{pmatrix}\({\sum}_r c_{r,q}\,e(rx)\)_\infty \ = \
\({\sum}_r c_{r,q}\,e(rx)\)_\infty\,
\\
 \text{and}\ \ \ \int_0^1\,\psi_{\l_1-\l_2,\d_3}\begin{pmatrix}1 & t \\ 0 & 1\end{pmatrix}\({\sum}_{r\in \Z}\
c_{r,q}\,e(rx)\)_\infty\ dt \ \ = \ \ 0\,.
\end{aligned}
\end{equation}
The space of distribution vectors in $W_{\l_1-\l_2,\d_3}^{-\infty}$ supported at $\infty$ has a natural increasing filtration by the {\it normal order\/}: by definition, the $k$-th derivative of the delta function at $\infty$ has order $k$. The upper triangular unipotent subgroup $N_\infty \subset {\rm SL}(2,\R)$ fixes the point at infinity, and therefore preserves this filtration. In fact,
\begin{equation}
\begin{gathered}
\label{inforder4}
\text{the action of}\,\  N_\infty \ \text{on\,}\ \{ \, \sigma \in W_{\l_1-\l_2,\d_3}^{-\infty}\, \mid \, \operatorname{suppt}\sigma\subset\{\infty\}\,\}\,\ \text{induces}
\\
\text{the identity on the successive quotients of the filtration}.
\end{gathered}
\end{equation}
To see this, note that $N_\infty$ is the unipotent radical of the isotropy subgroup at $\infty$. As such, it acts
trivially on the fiber at $\infty$ of any ${\rm SL}(2,\R)$-equivariant line bundle, in particular the line bundle
whose sections constitute $W_{\l_1-\l_2,\d_3}^{-\infty}$. For the same reason the vector field which
generates $N_\infty$ has a second order zero at $\infty$. Taken together, these two statements imply
(\ref{inforder4}).

Because of (\ref{inforder4}), $N_\infty$ acts unipotently on any finite dimensional subspace of $\{ \, \sigma \in W_{\l_1-\l_2,\d_3}^{-\infty} \mid \operatorname{suppt} \sigma\subset\{\infty\}\,\}$, for example the subspace of all sections of normal order at most $m$, for a given $m\geq 0$, which are invariant under some $n_\infty\in N_\infty$, $n_\infty\neq e$. This latter space is acted upon by the compact quotient $N_\infty/\{n_\infty^k\mid k\in\Z\}$. A compact group which acts unipotently must act trivially.
 Thus (\ref{inforder3}), (\ref{inforder4}) imply $(\sum_r c_{r,q}\, e(rx))_\infty = 0$ -- in other words, 
${\sum}_r c_{r,q}\, e(rx)$ vanishes to infinite order at $\infty$. Since the roles of $x$ and $y$ are related by
the outer automorphism (\ref{outer}), the analogous statement about ${\sum}_r c_{q,s}\,e(sy)$ is also
correct. Note that the identities in \corref{cor1} relate $\,\sigma_{cq,aq}\,$ and $\,\rho_{cq,aq}\,$ near
$x=0$ to, respectively, ${\sum}_r c_{r,q}\, e( rx)$ and ${\sum}_s c_{q,s}\, e(sy)$ near $x = \infty$. Hence,
\begin{equation}
\label{inforder6}
\text{the}\ \,\sigma_{n,k}\, \ \text{and}\ \,\rho_{n,k}\,\ \text{all vanish to infinite order at}\ \,0\,.
\end{equation}
Therefore, by \cite[Th.~3.19]{MilSch:2003b}, their Fourier transforms have canonical extensions across
infinity:
\begin{equation}
\label{inforder7}
\text{the}\ \,\sigma_{n,k}\, \ \text{and}\ \,\rho_{n,k}\,\ \text{all have canonical extensions across}\ \,\infty\,.
\end{equation}
As in the earlier discussion, one can interpret this alternatively and equivalently as a statement about distributions on $\R \cup \{\infty\}$ or about vectors in the appropriate representation spaces.

Like $\sum_r c_{r,q}\,e(rx)$, the distribution $\,\sigma_{n,k}\in C^{-\infty}(\R)$ has two potentially different extensions to vectors in $\,W_{\l_1-\l_2,\d_3}^{-\infty}$: the natural extension, i.e., in the manner specified in \corref{cor1}, and the canonical extension. We let $(\sigma_{n,k})_\infty \in W_{\l_1-\l_2,\d_3}^{-\infty}$ denote the difference of the two. Note that \corref{cor1} relates the behavior of $\,\sigma_{cq,\bar aq}$ at $\infty$ to that of $\,\sum_{r\in \Z}\, c_{r,q}\, e(rx)$ at $a/c$. In particular $\,\sum_{r\in \Z}\, c_{r,q}\, e(rx)$ vanishes to infinite order at $a/c$ if and only if $(\sigma_{cq,\bar aq})_\infty = 0$. But $(a,c)$ is an arbitrary pair of relatively prime integers with $c\neq 0$, so all the $\,\sum_{r\in \Z}\, c_{r,q}\, e(rx)$ vanish to infinite order at every rational point if and only if $(\sigma_{n,k})_\infty = 0$ for all $n\neq 0$ and $k$. When that is the case, another application of \corref{cor1} implies
also that the $\,\sigma_{n
 ,k}$ vanish to infinite order at all rational points. We can argue similarly in the case of $\,\rho_{n,k}$ and $\,\sum_{s\in \Z}\, c_{q,s}\, e(sy)$. Conclusion: to complete the proof of \lemref{inforder}, it suffices to show $(\sigma_{n,k})_\infty = 0$ and $(\rho_{n,k})_\infty = 0$ for all $n\neq 0$ and $k$.

We fix integers $\,a,\, c,\, q,\, \bar a\,$ as in (\ref{vora}), choose $b \in \Z$ so that $\,a \bar a - bc = 1$, and define
\begin{equation}
\label{inforder8}
\g \ = \ \ttt{1}{}{}{}{a}{b}{}{c}{\bar a} \ ,\ \ \ n_y \ = \ \ttt{1}{}{}{}{1}{1}{}{}{1}\ , \ \ \ n_{y,\g}\ = \ \g \,n_y\, \g^{-1}\,.
\end{equation}
Arguing as in the proof of \corref{cor1}, one finds that $p_{x,q}\circ p_{z,0}\,\tau\/$ is invariant under $n_y$ and $p_{y,0}\circ p_{x,q}\circ p_{z,0}\,\tau = 0$. Since $\g$ fixes $\tau$,
\begin{equation}
\label{inforder9}
\begin{aligned}
 \pi (\g) \circ p_{x,q}\circ p_{z,0}\circ \pi(\g^{-1})\tau\, = \, \pi(\g)\circ p_{x,q}\circ p_{z,0}\, \tau\ \, \text{is
invariant under}
\\
 n_{y,\g}\ \, \text{and}\ \, \(\pi(\g) \circ p_{y,0}\circ \pi(\g^{-1})\) \(\pi (\g) \circ p_{x,q}\circ p_{z,0}\circ
\pi(\g^{-1})\)\tau\,  = \, 0\,.
\end{aligned}
\end{equation}
The operator $p_{x,q}\circ p_{z,0}$ sends any vector in $V_{\l,\d}^{-\infty}$ to the component that transforms under $N_{x,z}$ according to the character $e(-qx)$. But $\g^{-1}$ normalizes $N_{x,z}$ and acts on this group via $(x,z) \mapsto ( a x + cz, bx + \bar az)$, from which one can deduce
\begin{equation}
\label{inforder10}
\pi (\g) \circ p_{x,q}\circ p_{z,0}\circ \pi(\g^{-1}) \ = \ p_{x,aq}\circ p_{z,cq}\,.
\end{equation}
Recall the definition of the codimension one Schubert cell $\{x=\infty\}$, following (\ref{proj7}). We use the matrix identity
\begin{equation}
\label{coordinates1}
\ttt{1}{x}{z}{}{1}{y}{}{}{1} \ = \ \ttt{1}{}{z}{}{1}{y}{}{}{1} \ttt{1}{x}{}{}{1}{}{}{}{1}\,,
\end{equation}

to extend the coordinates $(x,y,z)$ from the open cell $N$ to $N\cup\{x=\infty\}$, with $x$ expressed as $x=1/u$ near $x=\infty$, of course. This is legitimate because the second factor on the right lies in the group $S \simeq {\rm SL}(2,\R)$ -- cf.\ (\ref{proj7}) -- whose orbit through the identity coset is a copy of $\C\mathbb P^1 \simeq \R\cup\{\infty\}$. Proposition \ref{prop1} and the proof of \corref{cor1} imply
\begin{equation}
\label{inforder11}
\begin{array}{l}
\(p_{z,cq}\tau\)(x,y,z) 
\\
\qquad = \ \sum_{k\in \Z/cq\Z}\, \sum_{\stackrel{\scriptstyle{\ell\equiv
k}}{\operatorname{mod}(cq)}}\,e(cqz +
\ell y) \(\psi_{\l_1-\l_2,\d_3}\ttwo{1}{-\f{\ell}{cq}}{}{1}\sigma_{cq,k}\)\!(x)\,, \end{array}
\end{equation}

as an identity valid not only on $N$, but even on $N\cup\{x=\infty\}$. In view of (\ref{inforder10}) and the multiplication law (\ref{multiplication}),
\begin{equation} 
\label{inforder12}\begin{array}{l}
 \(\pi (\g) \circ p_{x,q}\circ p_{z,0}\circ \pi(\g^{-1})\,\tau\)(x,y,z) 
\\ \qquad
 = \int_0^1  \sum_{\stackrel{\scriptstyle{k\in \Z/cq\Z}}{\ell\equiv k\, (cq)}}  e(aqt  + cq(z-ty) + \ell
y)\!\(\! \psi_{\l_1-\l_2,\d_3} \! \ttwo{\!\! 1}{\!\!\f{cqt - \ell}{cq}\!}{\!\!}{\!\!1\!} \!
\sigma_{cq,k}\!\)\!\!(x)\,dt, \end{array}
 \end{equation}
\pagebreak
again as an identity between distribution sections of 
$\mathcal L_{\lambda,\delta}$ over $N\cup\{x=\infty\}$.

When we write $\sigma_{cq,k}=(\sigma_{cq,k})_{\text{can}}+(\sigma_{cq,k})_\infty$
 as the sum of the canonical extension of $\,\sigma_{cq,k}|_\R$ and a summand supported at infinity, the
integrand in (\ref{inforder12}) splits into the sum of two terms, one of which vanishes to infinite order along
$\{x=\infty\}$ \cite[Lemma 7.2]{MilSch:2003b}\footnote{To apply \cite[Lemma 7.2]{MilSch:2003b},
one needs to know that the canonical extensions of the summands in (\ref{inforder12}) vanish at $\infty$ 
uniformly \cite[Def.~7.1]{MilSch:2003b}. One can show this directly; alternatively, one can argue as
in the proof of the stronger statement (\ref{leminforder1}) in \S\ref{convsec}.}, and the other is supported
on
$\{x=\infty\}$. According to \cite[Prop.~7.20]{MilSch:2003b}, this remains the case after the
integration is performed:
\begin{equation}
\label{inforder13}
\pi (\g) \circ p_{x,q}\circ p_{z,0}\circ \pi(\g^{-1})\,\tau \ \ = \ \ \tau_{\text{inforder}} \ +\ \tau_{\{x=\infty\}}\,,
\end{equation}
where $\,\tau_{\text{inforder}}\,$ stands for a distribution section over $N\cup\{x=\infty\}$ which vanishes to infinite order along $\{x=\infty\}$,
with
\begin{equation}
\label{inforder13.5}
\begin{array}{l}
 \tau_{\{x=\infty\}}(x,y,z)  
\\
\qquad \displaystyle{= \,\int_0^1 \!\!\sum_{\stackrel{\scriptstyle{k\in \Z/cq\Z}}{\ell\equiv k\, (cq)}} \!\! e(aqt\! +
cq(z-ty)\! +
\ell y)\,\, \psi_{\l_1-\l_2,\d_3} \! \ttwo{\!\! 1}{\!\f{cqt - \ell}{cq}\!}{\!\!}{\!1\!}  (\sigma_{cq,k})_\infty\,dt
}\end{array}
\end{equation}
supported on $\{x=\infty\}$. We shall extend the identity (\ref{inforder13}) to a larger set. For this purpose, we need to know that
\begin{equation}
\label{inforder14}
\begin{array}{l}
 U =_{\rm def} \ N \cup \g N \cup \{x=\infty\} \cup \g\{x=\infty\}\ \ \text{is open in}\ \ G/P\,,
\\
 T =_{\rm def} \ \{x=\infty\} \cup \g\{x=\infty\}\ \ \text{is a closed submanifold of}\ \ U\,,
\\
  \text{and both $\,U\,$ and $\,T\,$ are $\, \g N_{y,z} \g^{-1}$-invariant}\,.
\end{array}
\end{equation}
This can be verified directly or, more easily, deduced from the explicit description of the Schubert cells in
Section \ref{convsec}, in the proof of \lemref{lemexist}. Note that $\g N_{y,z} \g^{-1}$ contains the
 one-parameter unipotent group through $n_{y,\g}$; the operator $\pi(\g)\circ p_{y,0}\circ \pi(\g^{-1})$
involves integration over this one-parameter group modulo the cyclic group generated by $n_{y,\g}$. We
claim:
\begin{equation}
\label{inforder15}
\begin{gathered}
\text{The identity (\ref{inforder13}) is valid on all of $\,U\,$, with $\,\tau_{\text{inforder}}\,$ now
denoting}
\\
\text{an expression which vanishes to infinite order along $\,T$, and}
\\
\text{with}\,\ \tau_{\{x=\infty\}}\,\ \text{supported on the complement of $\,\g\{x=\infty\}\,$ in $T$.}
\end{gathered}
\end{equation}
To see this, it suffices to show that $\pi (\g) \circ p_{x,q}\circ p_{z,0}\circ \pi(\g^{-1})\,\tau$ vanishes to infinite order along $\g\{x=\infty\}$, or equivalently, in view of (\ref{inforder9}), that
\begin{equation}
\label{inforder16}
\bigl(p_{x,q}\circ p_{z,0}\,\tau\bigr)(x,y,z)\ =\ {\sum}_{s\neq 0}\ c_{q,s}\, e(qx+sy)
\end{equation}
vanishes to infinite order along $\{x=\infty\}$. In the first part of the proof we showed that ${\sum}_{q\neq 0}\ c_{q,s}\, e(qx+sy)$ vanishes to infinite order along $\{x=\infty\}$, for each $s\neq 0$\,. Hence so does  ${\sum}_{q,s\neq 0}\ c_{q,s}\, e(qx+sy)$
 \cite[Lemmas 7.2 and~7.13]{MilSch:2003b}. Applying the projection operator $p_{x,q}$ does not destroy
this property \cite[Prop.~7.20]{MilSch:2003b}, so that (\ref{inforder15}) follows.

Because of (\ref{inforder14}), (\ref{inforder15}) and \cite[Prop.~7.20]{MilSch:2003b}, both
summands in (\ref{inforder13}) inherit from  $\pi(\g)\circ p_{x,q}\circ p_{z,0}\circ \pi(\g^{-1})\,\tau$ the
two invariance properties mentioned in (\ref{inforder9}). In particular,
\begin{equation}
\label{inforder17}
\begin{gathered}
\tau_{\{x=\infty\}}\,\ \text{is invariant under $\,n_{y,\g}\,$ and is annihilated}
\\
\text{by $\,\pi(\g)\circ p_{y,0}\circ \pi(\g^{-1})$}\,.
\end{gathered}
\end{equation}
A straightforward calculation shows that the two sets $\,\{x=\infty\,,\ y=a/c\}\,$ and $\,\g\{x=\infty\}\,$ cover $T$ disjointly. Hence, in view of (\ref{inforder15}),
\begin{equation}
\label{inforder18}
\tau_{\{x=\infty\}}\,\ \text{is supported on}\,\ \{x=\infty\,,\ y=a/c\}\ \subset\ T\,.
\end{equation}
We regard $\frac{\partial\ }{\partial y}$ as an element of the Lie algebra of $N$, as in (\ref{partials}). Then $\g \frac{\partial\ }{\partial y} \g^{-1}$ generates the
one-parameter unipotent subgroup of ${\rm GL}(3,\R)$ passing through $n_{y,\g}$. This subgroup fixes
$\,\{x=\infty\,,\ y=a/c\}\,$ pointwise -- equivalently, the generating vector field $\g \frac{\partial\ }{\partial
y} \g^{-1}$ vanishes along $\,\{x=\infty\,,\ y=a/c\}\,$. More precisely, in terms of the local coordinates
$\,u=1/x$\,, $\,v=y-a/c$\,, $\,z$\,,
\begin{equation}
\label{inforder19}
\begin{gathered}
\g \frac{\partial\ }{\partial y} \g^{-1}\ = \ h_1(u,v,z)\,\frac{\partial\ }{\partial u}\ + \ h_2(u,v,z)\,\frac{\partial\ }{\partial v}\ + \ h_3(u,v,z)\,\frac{\partial\ }{\partial z}\,,
\\
\text{with $h_1$, $h_2$ vanishing to second order along $\{x=\infty\,,\ y=a/c\}$}
\\
\text{and $h_3$ vanishing to first order}\,.
\end{gathered}
\end{equation}
One way to see this is to note that $g\g^{-1}$, with $g = \(\begin{smallmatrix} &  & {1} \\ {1} &  &  \\ & 1 &  \\
\end{smallmatrix}\)$, carries the submanifold $\,\{x=\infty\,,\ y=a/c\}\,$ onto $\{x=z=0\} \subset N$, and that
$$
g\ttt{1}{}{}{}{1}{t}{}{}{1} g^{-1}\ttt{1}{x}{z}{}{1}{y}{}{}{1} = \ttt{1}{\frac{x}{1+t(z-xy)}}{\frac{z}{1+tz}}{}{1}{\frac{y}{1+tz}}{}{}{1} \ttt{*}{}{}{*}{*}{}{*}{*}{*}\,,
$$
which implies $g \frac{\partial\ }{\partial y}g^{-1}= -x(z-xy)\frac{\partial\ }{\partial x}-yz\frac{\partial\ }{\partial y}-z^2\frac{\partial\ }{\partial z}$. 
This latter identity in\break\vskip-12pt\noindent turn implies (\ref{inforder19}).

The space of distribution sections of $\mathcal L_{\l,\d}$ over $U$ with support in the codimension two submanifold $\{x=\infty\,,\ y=a/c\}$ is filtered by the normal order. This is analogous to the normal order filtration discussed in the first part of the proof; it counts the total number of normal derivatives, in both conormal directions. Because of (\ref{inforder19}), the vector field $\g \frac{\partial\ }{\partial y} \g^{-1}$ lowers the normal order by one. In general the normal order is well defined only locally, but $\tau_{\{x=\infty\}}$ is periodic in the variable $z$ and therefore has a globally well defined normal order. It follows that $\tau_{\{x=\infty\}}$ lies in
 a finite-dimensional space of distribution sections, invariant under the action of the one-parameter subgroup
generated by $\g\frac{\partial\ }{\partial y} \g^{-1}$, on which this group acts unipotently. We now appeal to
(\ref{inforder17}) and argue the same way as in the first part of the proof, to conclude that
$\,\tau_{\{x=\infty\}}\,$ vanishes. Next we multiply $\,\tau_{\{x=\infty\}}=0\,$ by $e(-cqz)$, express
$\,\tau_{\{x=\infty\}}\,$ as in (\ref{inforder13.5}), substitute $v + a/c$ for $y$, note that $e(\ell
a/c)=e(ka/c)$ when $\ell \equiv k\,\operatorname{mod}(cq)$\,, and then unfold the integration with respect
to $t\,$:
\begin{equation}
\label{inforder20}
\begin{aligned}
0 \ &=\,  \int_0^1 \!\!\sum_{\stackrel{\scriptstyle{k\in \Z/cq\Z}}{\ell\equiv k\, (cq)}}\!\! e((\ell - cqt)v + ka/c)\,\, \psi_{\l_1-\l_2,\d_3} \! \ttwo{\!\! 1}{\!\f{cqt - \ell}{cq}\!}{\!\!}{\!1\!}  (\sigma_{cq,k})_\infty\,dt
\\
&= \sum_{k\in \Z/cq\Z} e(ka/c)\int_\R  e(- cqvt)\,\, \psi_{\l_1-\l_2,\d_3} \! \ttwo{\!\! 1}{\!t\!}{\!\!}{\!1\!}  (\sigma_{cq,k})_\infty\,dt\,.
\end{aligned}
\end{equation}
When we reinterpret this identity on the level of the top quotient of the filtration referred to in (\ref{inforder4}) -- i.e., the top level quotient at which there is a nonzero component -- $(\sigma_{cq,k})_\infty$ no longer gets translated, since $N_\infty$ acts tri\-vially on the quotient. That leaves the integral $\int_\R  e(- cqvt)\,dt \ = \ |cq|^{-1}\,\d_0(v)$ as a factor which can be pulled out. Conclusion:
\begin{equation}
\label{inforder21}
{\sum}_{k \in \Z/cq\Z}\,e(ka/c) \,(\sigma_{cq,k})_\infty \ = \ 0\,,
\end{equation}
since a nonzero ``top component" cannot occur at any level.

To simplify the notation, we write $n$ for the product $cq$. Any fraction with denominator $n$ can be expressed as $a/c$ with $a$ and $c$ as above, so
that (\ref{inforder21}) asserts the vanishing of the finite Fourier transform of $\,\Z/n\Z \ni k \mapsto
(\sigma_{n,k})_\infty\,$, and therefore implies $(\sigma_{n,k})_\infty=0$ for all $n\neq 0$ and $k\in
\Z/cq\Z$. The $\rho_{n,k}$ are related to the $\sigma_{n,k}$ by an outer automorphism of $G$ which
stabilizes $N$. It follows that the $(\rho_{n,k})_\infty$ must also be zero. As was pointed out before, this
completes the proof of \lemref{inforder}.\hfill\qed

\section{$L$-functions}\label{lfuncsec}

In this section we derive the analytic continuation and functional equation of the standard $L$-function of a cuspidal ${\rm GL}(3,\Z)$-automorphic representation of ${\rm GL}(3,\R)$, and more generally, of the standard $L$-function twisted by a primitive Dirichlet character.

Let us recall the existing methods of proof. Godement-Jacquet
first established the analytic continuation and functional
equation via an integral over the entire group $G={\rm GL}(3,\R)$
\cite{GodJac:1972}. Jacquet, Piatetski-Shapiro, and Shalika later
found an argument using Mellin transforms over one-dimensional
tori of $G$ \cite{JacPiatShal:1979}; see also \cite{Bump:1984}.
Their approach enabled them to give a necessary and sufficient
condition for automorphy in their converse theorem, a topic we
shall return to in \secref{convsec}. A third method, of Langlands
and Shahidi, uses the constant terms of maximal parabolic
Eisenstein series on ${\rm GL}(4)$ \cite{Shahidi:1985}.

The method we present here has some similarities to that of \cite{JacPiatShal:1979}, but totally avoids the use of Whittaker functions.  Our arguments apply equally to cusp forms on quotients ${\rm GL}(3,\R)$ by congruence subgroups of $\G={\rm GL}(3,\Z)$. We have chosen to limit the discussion
to the case of full level, to explain most transparently how our method works.

Let $\tau$ be a cuspidal automorphic distribution for $\G={\rm GL}(3,\Z)$ as in \secref{voronsec}, with normalized abelian Fourier coefficients $a_{n,m}$ as in (\ref{cnorm}). The standard and contragredient $L$-functions of $\tau$ are defined as
\begin{equation}
\label{stanlfuncdef}
L(s,\tau) \ \  = \ \ \sum_{n=1}^\infty
\,a_{1,n}\,n^{-s}\ \ \ \ \text{and}\ \ \ \ L(s,\tilde{\tau}) \ \ = \ \
\sum_{n=1}^\infty \, a_{n,1}\,{n^{-s}},
\end{equation}
respectively.  The latter is the standard $L$-function for the automorphic distribution $\tilde\tau$ obtained from $\tau$ by conjugation with the outer automorphism
\begin{equation}
\label{outer1}
g \ \longmapsto \ w_{\text{long}}\, (g^t)\i\, w_{\text{long}}^{-1}\,,\ \ \ \text{with}\ \ \ w_{\text{long}}\ = \ w_{\text{long}}^{-1}\ = \ \ttt{}{}{1}{}{1}{}{1}{}{}\,,
\end{equation}
whose restriction to $N$ we encountered before, in (\ref{outer}). When $\tau$ corres\-ponds to the embedding (\ref{newtau}) of the irreducible unitary representation $\pi'$, $\tilde\tau$ corres\-ponds to the embedding of the dual representation $\pi$, via the composition of (\ref{newtau}) with the outer automorphism (\ref{outer1}). The series (\ref{stanlfuncdef}) converge uniformly and absolutely for $\,\operatorname{Re}s\,$ sufficiently large, because the Fourier coefficients of a periodic distribution grow at most polynomially in terms of the index.

It is important to note that (\ref{stanlfuncdef}) agrees with the usual definition -- e.g., in
\cite[\S8]{Bump:1984}~-- of the $L$-function of an automorphic form $\phi \in
L^2_\omega(\G\backslash G)$ inside the direct summand determined by $\tau$. One can see this by relating
the Fourier-Whittaker expansion of $\phi$, as in \cite{Piatetski:1979} or \cite{Shalika:1974}, to the abelian
part of the Fourier expansion of $\tau$ --
  in effect, by proving a ${\rm GL}(3)$-analogue of \cite[Lemma 5.65]{Schmid:2000}. However, it is much
easier to argue using the Hecke action: when $\tau$ is an eigendistribution of the Hecke operators $T_{1,p}$
and $T_{p,1}$, the $a_{r,s}$ can be renormalized so that $a_{1,1}=1$, in which case the $a_{1,n}$ and
$a_{n,1}$ coincide with the Hecke eigenvalues of
$\phi$;
 see (\ref{t1pp1})--(\ref{eulerprod3}) in the next section.

A Dirichlet character modulo $q$ is a function $\chi:\Z \rightarrow \C$ obtained by lifting a character of $\,(\Z/q\Z)^*\,$ and then extending it to all of $\,\Z\,$ by defining $\,\chi(n)=0\,$ when $(n,q)\neq 1$. A Dirichlet character $\chi'$ modulo $q'$ is said to be {\em induced\,} from $\chi$ if $q$ divides $q'$, and $\chi$ and $\chi'$ agree on all integers relatively prime to $q'$. When a character $\chi$ modulo $q$ is not properly induced one calls $\chi$ {\em primitive} and $q$ the {\em conductor}. The standard $L$-functions twisted by the primitive Dirichlet character $\chi$ are defined as the series
\begin{equation}
\label{twistlfuncdef}
L(s,\tau\otimes \chi) \ =  \ \sum_{n=1}^\infty a_{1,n}\,\chi(n)\, n^{-s}\text{~~and~~} L(s,\tilde{\tau}\otimes \chi) \  = \ \sum_{n=1}^\infty a_{n,1}\,\chi(n)\, n^{-s} .
\end{equation}
These, too, converge for $\Re{s}$ sufficiently large. When $\chi \equiv 1$ is the trivial Dirichlet character, $L(s,\tau\otimes \chi)$ coincides with the standard $L$-function $L(s,\tau)$.

The statement and proof of the functional equation involve the $\G$-function and various $\G$-identities. We begin with the notational conventions
\begin{equation}
\label{gammarc}
\G_\R(s)\ = \ \pi^{-s/2}\,\G(s/2)\,,\ \ \ \G_\C(s)\ =\ 2\,(2\pi)^{-s}\,\G(s)\,,
\end{equation}
which are commonly used in the context of $L$-functions. The classical duplication formula relates these two modified $\G$-functions:
\begin{equation}
\label{duplication}
\G_\C(s)\ = \ \G_\R(s)\,\G_\R(s+1)\,.
\end{equation}
For $\d \in \Z/2\Z$, and initially only for $\, 0 < \operatorname{Re} s < 1\,$, we define
\begin{equation}
\label{whyg} G_\d(s) \ \ = \
\int_\R \, e(x)\, |x|^{s-1}\, \sg(x)^\d \, dx \qquad \( \, 0<\operatorname{Re}s < 1\,\) \,.
\end{equation}
Both $G_0(s)$ and $G_1(s)$ extend meromorphically to the entire complex plane:
\begin{equation}
\label{heresG1}
G_0(s) \ = \ \G_\C(s)\, \cos(\pi s / 2)\,, \qquad  G_1(s) \ = \ i\,\G_\C(s)\, \sin(\pi s / 2)\,,
\end{equation}
as follows from the standard identity $\,\int_0^\infty e^{iax} x^{s-1} dx = a^{-s}\G(s)e^{i \pi s/2}$ which holds when $0 < a\,$, $\,0 < \operatorname{Re}s < 1$\,. Since $\G(s)$ is never zero and has simple poles only, at the points $s=-n\,$, $\,n = 0,\, 1,\, 2,\, \dots$\,,
\begin{equation}
\label{dgpoleszeros}
\begin{array}{l}
 \text{the poles of $G_\d(s)$ are all simple and lie at}\, \ s \in (2\Z + \d)\cap \Z_{\leq 0} \,,
\\[4pt]
 \text{the zeros of $G_\d(s)$ are all simple and lie at}\, \ s \in (2\Z + \d+ 1)\cap \Z_{> 0} \,.
\end{array}
\end{equation}
In view of (\ref{heresG1}), the identity
\begin{equation}
\label{heresG3}
G_\d(s)\,G_\d(1-s)  \ \ = \ \ (-1)^\d
\end{equation}
paraphrases the well known functional equation for the $\G$-function,
\begin{equation}
\label{gfuncequ}
\G(s) \,\G(1 - s) \ = \ \pi\,\csc(\pi s)\,.
\end{equation}
This functional equation, in combination with the duplication identity (\ref{duplication}), also implies
\begin{equation}
\label{heresG}
G_\d(s)  \ \ = \ \ i^\d\, \f{\G_\R(\,s\,+\,\d\,)}{\G_\R(\,1\,-\,s\,+\,\d\,)}\ \ \ \ \(\, \d\, \in\, \{0,1\}\,\)\,.
\end{equation}
Unlike (\ref{whyg}) and (\ref{heresG3}), which visibly depend only on the parity of $\d$, (\ref{heresG}) becomes incorrect without the condition $\,\d \in \{0,1\}\,$. We conclude this discussion of $\G$-identities with a less well known formula,
\begin{equation}
\label{heresG2}
\begin{array}{l}
 \l_1 \, - \, \l_2 \ \in \ 2\,\Z + \d_1 - \d_2 +1\ \ \ 
\\
 \qquad \Longrightarrow G_{\d_1}(s+\l_1)\,G_{\d_2}(s+\l_2) \  = \  i^{\l_1-\l_2+1} \,
\displaystyle\f{\G_\C(s+\l_1)}{\G_\C(1-s-\l_2)} \ ,
\end{array}
\end{equation}
which can be deduced from the functional equation (\ref{gfuncequ}) and the standard 
\pagebreak
identity $\G(s+1) = s\,\G(s)$.

The functional equation for the twisted $L$-functions also involves the Gauss sum
\begin{equation}
\label{gausssumdef}
g_\chi \ = \ {\sum}_{k\in \Z/q\Z}\ \chi(k)\,e(k/q)\,.
\end{equation}
If $\chi$ is primitive,
\begin{equation}
\label{gausssumidentity}
\begin{array}{rl}
 \text{a)}&\quad \chi(-1)\,g_{\bar\chi} \ \ = \ \ \overline{g_\chi} \ \ = \ \ q\,g_\chi^{-1}\,,\ \ \text{and}
\\[4pt]
 \text{b)}&\quad {\sum}_{k\in \Z/q\Z}\ \chi(k)\,e(nk/q)\ \ = \ \ \bar\chi(n)\,g_\chi\,,\ \ \text{for all}\ \ n \in
\Z\,;
\end{array}
\end{equation}
see \cite[Lemma 3.63]{Shimura:1971}, for example.

\begin{thm}\label{funeqthm}  Let $\tau\,\in\,(V_{\l,\d}^{-\infty})^\G$ 
be a cuspidal automorphic distribution for $\G={\rm GL}(3,\Z)$ as in {\rm (\ref{newtau}),}
 and $\chi$ a primitive
Dirichlet character modulo\break $q>0$. Denote the parity of $\,\chi\,$ by $\,\chi(-1)=(-1)^\e\,${\rm ,} $\e \in
\Z/2\Z$. Then ${L(s,\tau\,\otimes\,\chi)}${\rm ,} $L(s,\tilde\tau\,\otimes\,\bar\chi)${\rm ,}
 $G_{\e+\d_3}(s+\l_3) L(s,\tau \otimes
\chi)$ and $G_{\e+\d_1}(s-\l_1)L(s,\tilde{\tau}\otimes\bar\chi)$ all have analytic continuations to entire
functions on $\C$ which are polynomially bounded on vertical strips. These functions obey the functional
equation
$$
L(1-s,\tilde\tau\otimes \bar\chi) \ \ = \ \ q^{3s}\,g_\chi^{-3} \left(\prod_{j=1}^3 G_{\e+\d_j}(s + \l_j)\right) L(s,\tau\otimes\chi)\,.
$$
\end{thm}

Our statement of the functional equation applies uniformly to all automorphic representations. The usual statements, such as those in \cite{GodJac:1972}, \cite{JacPiatShal:1979}, \cite{Shahidi:1985}, distinguish among representations depending on where they show up in the Langlands classification. It
is a  simple matter to translate between the different versions of the functional equation using the
$\G$-identities (\ref{gammarc})--(\ref{heresG2}). If an irreducible unitary representation $(\pi,V)$ occurs as
a cuspidal, ${\rm GL}(3,\Z)$-automorphic representation, there are only two cases to consider. Either
\begin{equation}
\label{psirreducibility}
\l_i \, - \, \l_j \, + \, \d_i \, - \, \d_j\ \notin \ 2\Z  + 1 \ \ \ \text{for}\ \ 1 \leq i < j \leq 3
\end{equation}
or $\l_i - \l_j + \d_i - \d_j$ is an odd integer for some $i<j$. In the former situation, the Langlands classification associates $(\pi,V)$ 
to the minimal parabolic subgroup, and the usual statements of the functional equation assert that
\begin{equation}
\label{psstandardversion}
\begin{aligned}
{}&\G_\R(s+\l_1+\d_1')\, \G_\R(s+\l_2+\d_2')\, \G_\R(s+\l_3 + \d_3')\,
L(s,\tau\otimes\chi)\,,
\\
&\hskip1.6in \text{with} \ \ \ \d_j' \in\{0,1\}\,,\ \ \ \ \d_j' \equiv\e + \d_j\ (\operatorname{mod} 2) \,,
\end{aligned}
\end{equation}
is an entire function, which agrees with the analogous expression invol\-ving $L(1-s,\tilde\tau\otimes\bar\chi)$, up to a constant factor and powers of $q$. The equality of the two expressions follows immediately from our version of the functional equation, in conjunction with (\ref{heresG}). The nonparity condition (\ref{psirreducibility}), together with (\ref{dgpoleszeros}), ensures that no poles of any one of the three factors $G_{\e + \d_j}(\dots)$ on the right-hand side of our functional equation are obliterated by a zero of one of the other factors. Another application of (\ref{dgpoleszeros}) shows that the poles of these three factors occur precisely at the points where the three $\G$-factors in (\ref{psstandardversion}) have poles. Thus, in view of our functional equation, the expression (\ref{psstandardversion}) is indeed entire.

If (\ref{psirreducibility}) fails for a cuspidal, ${\rm GL}(3,\Z)$-automorphic representation, the Langlands classification attaches it to a maximal parabolic subgroup. One can then choose the Casselman embedding (\ref{embed2}) so that
\begin{equation}
\label{psreducibility}
\begin{array}{rl}
\text{a)\ \ }& \l_1  -  \l_2 \, \in \, (2\Z + \d_1 -\d_2 + 1)\cap \Z_{\leq 0}\,,\ \ \ \l_2  -  \l_3 \, \notin \, \Z_{<0}\ \ \
\text{or}
\\[4pt]
\text{b)\ \ }& \l_2  -  \l_3 \, \in \, (2\Z + \d_2 -\d_3 + 1)\cap \Z_{\leq 0}\,,\ \ \ \l_1  -  \l_2 \, \notin \, \Z_{<0}\,.
\end{array}
\end{equation}
The two cases are related by the outer automorphism (\ref{outer1}), so we may as well suppose that we are in the situation a). The dual representation then corresponds to the case b), with the dual choice of $(\l,\d)$. According to the usual statements of the functional equation,
\begin{equation}
\label{nonpsstandardversion}
\begin{array}{l}
 \G_\C(s+\l_2)\, \G_\R(s+\l_3+\d_3')\,L(s,\tau\otimes\chi)\,,
\\[4pt]
 \qquad\qquad\qquad \ \ \ \text{with}\ \ \d_3' \in \{0,1\}\,,\ \ \d_3'\equiv \e + \d_3 \ \ (\operatorname{mod} 2)\,,
\end{array}
\end{equation}
is entire and agrees -- up to a factor involving $g_\chi$, $q$, and powers of $i$ -- with the analogous expression corresponding to the dual data. This follows from Theorem \ref{funeqthm}, in conjunction with
(\ref{dgpoleszeros}),  (\ref{heresG}), and (\ref{heresG2}).

When $q=1$ in the statement of Theorem \ref{funeqthm}, $\chi$ is identically 1, $g_\chi = 1$ and $\e = 0$. In that case, $L(s,\tau\otimes \chi)$ reduces to the standard $L$-function, which
satisfies the slightly simpler func\-tional equation
\begin{equation}\label{stanfuneq}
L(1-s,\tilde\tau) \ \ = \ \ G_{\d_1}(s+\l_1)\, G_{\d_2}(s+\l_2)\, G_{\d_3}(s+\l_3)\, L(s,\tau)\,.
\end{equation}
We remarked in the introduction that (\ref{stanfuneq}) can be formally derived from our main theorem, with
$f(x)=|x|^{-s}$. While the statement of the main theorem does not allow this choice of $f(x)$, the proof can
be adapted, and even simplifies with $f(x)=|x|^{-s}$. This is the strategy we follow below, not just for
$L(s,\tau)$ but also for $L(s,\tau\otimes \chi)$. The latter case requires some combinatorial arguments; we
present these in a form which is useful also for the proof of the converse theorem in \secref{convsec}. Our
point of departure is the symmetric form of the identity between the $\sigma_{n,k}$ and $\rho_{n,k}$ in
\propref{prop1},
\begin{equation}
\label{star2}
{\sum}_{k=1}^q \ a_k \ \rho_{q, \,k}(y) \ \  =  \ \ {\sum}_{\ell=1}^q \   \widehat{a}_\ell \  \widehat{\sigma}_{q, \, \ell}(qy)\,,
\end{equation}
rather than the asymmetric form used in (\ref{vorc}) to derive the Voronoi formula.

\demo{Proof  of Theorem {\rm \ref{funeqthm}}} We set $\a_k=\chi(k)$ in (\ref{star2}). Then by
(\ref{gausssumidentity}b) and the definition (\ref{fourierint3}) of the finite Fourier transform,
\begin{equation}
\label{ahatk}
{\sum}_{k \in \Z/q\Z}\ \chi(k) \, \rho_{q, \,k}(y)\ \ = \ \ g_\chi\, {\sum}_{\ell \in \Z/q\Z}\ \bar\chi(\ell)\, \widehat{\sigma}_{q, \, \ell}(qy)\,.
\end{equation}
According to \lemref{inforder}, both sides of (\ref{ahatk}) vanish to infinite order at\break  $y=0$ and have
canonical extensions across $\infty$, as do their Fourier transforms. Their signed Mellin transforms are
therefore entire functions of $\,s\,$ and are polynomially bounded on vertical strips \cite[\S
4]{MilSch:2003b}. We integrate (\ref{ahatk}) against the shifted Mellin kernel $\,|t|^{-s-\l_2}(\sg t)^{\e+\d_2}$\,:
\begin{equation}
\label{starmell1}
\begin{aligned}
{}&{\sum}_{k \in (\Z/q\Z)^*} \ \chi(k) \, \bigl( M_{\varepsilon+\delta_2}\rho_{q,k}\bigr)\,(1-s-\l_2)   
\\
{}&\qquad  =\ g_\chi \, q^{s+\l_2-1}\, {\sum}_{\ell \in (\Z/q\Z)^*} \ \bar\chi(\ell)\, \bigl(
M_{\varepsilon+\delta_2}\widehat\sigma_{q,\ell}\bigr)\,(1-s-\l_2)
\\
{}&\qquad = \ (-1)^{\e+\d_2}\, g_\chi \, q^{s+\l_2-1}\, G_{\e+\d_2}(1-s-\l_2)\ \ \times
\\
{}&\qquad\quad \times\ \  {\sum}_{\ell \in (\Z/q\Z)^*} \ \bar\chi(\ell) \, \bigl(
M_{\varepsilon+\delta_2}\sigma_{q,\ell}\bigr)\,(s+\l_2)
 \,;
\end{aligned}
\end{equation}
the second step is justified by \cite[Th.~4.12]{MilSch:2003b}. The index $k$ on the left-hand side is
relatively prime to $q$. We recall the hypothesis $q > 0$, and we let $\bar k$ denote the multiplicative
inverse of $k$ modulo $q$. Then
\begin{equation}
\label{starmell1.1}
\rho_{q,k}(t) \ = \ (\sg t)^{\d_1}\,|qt|^{\l_2 - \l_3 -1} \, {\sum}_{n \neq 0}\ a_{1,n} \, |n|^{\l_3} \,(\sg n)^{\d_3}\, e( \textstyle \f{n}{q^{2}t} - \f{n \bar k}{q}) \,,
\end{equation}
by \propref{prop2} and the relation (\ref{cnorm}) between the $c_{n,k}$ and $a_{n,k}$. We substitute this sum for $\rho_{q,k}$ on the left in (\ref{starmell1}), change variables from $t$ to $t^{-1}$, use the restriction (\ref{delta3}) on the $\delta_j$ and appeal to
\cite[Lemma 3.38]{MilSch:2003b}, which allows us to compute the Mellin transform of the resulting series
formally:
\begin{equation}
\label{starmell1.2}
\begin{aligned}
{}&\bigl( M_{\varepsilon+\delta_2}\rho_{q,k}\bigr)\,(1-s-\l_2)  
\\
{}&\ \ \ \ = \ q^{\lambda_2-\lambda_3 - 1}\,M_{\varepsilon+\delta_3} \bigl( \textstyle {\sum}_{n \neq 0}\
a_{1,n} \,|n|^{\l_3} \,(\sg n)^{\d_3} \, e(\textstyle \f{nt}{q^{2}} - \f{n \bar k}{q})\bigr)(s + \lambda_3)
\\
{}&\ \ \ \ = \ q^{2s + \l_2 + \l_3 - 1}\, G_{\e + \d_3}(s + \l_3)\, {\sum}_{n \neq 0}\ a_{1,n} \,|n|^{-s} \,(\sg
n)^{\e} \, e(-\textstyle \f{n \bar k}{q}) \,,
\end{aligned}
\end{equation}
provided $\operatorname{Re}s \gg 0$. Since $(\sg n)^{\e} = \chi(\sg n)$ and $(-1)^\e = \chi(-1)$, the identity (\ref{gausssumidentity}) implies
\begin{equation}
\label{starmell1.3}
\begin{aligned}
{}&(\sg n)^{\e}\ {\sum}_{k \in (\Z/q\Z)^*}\ \chi(k)\, e(-n\bar k/q) 
\\
{}&\qquad = \ {\sum}_{k \in (\Z/q\Z)^*}\ \bar\chi(-k \sg n)\, e(n k/q) \ = \ (-1)^\e\, \chi(|n|)\,g_{\bar\chi}\,.
\end{aligned}
\end{equation}
In view of (\ref{starmell1.2})--(\ref{starmell1.3}),
\begin{equation}
\label{starmellout3}
\begin{aligned}
{}&{\sum}_{k \in (\Z/q\Z)^*} \ \chi(k) \, \bigl( M_{\varepsilon+\delta_2}\rho_{q,k}\bigr)\,(1-s-\l_2)  
\\
{}&\qquad\qquad = \ 2\, (-1)^\e\, g_{\bar\chi}\, q^{2s + \l_2 + \l_3 - 1}\, G_{\e + \d_3}(s + \l_3)  \,
L(s,\tau\otimes \chi)\,.
\end{aligned}
\end{equation}
This is the left-hand side of (\ref{starmell1}), which we already know is an entire function, polynomially bounded on vertical strips. The same assertion can be made about $G_{\e + \d_3}(s + \l_3)  \, L(s,\tau\otimes \chi)$, since $\,g_{\bar\chi}\neq 0$ by (\ref{gausssumidentity}a).

With small modifications, the preceding argument applies also to the expression in parentheses on the right-hand side of (\ref{starmell1}). We merely replace $\,a_{1,n}\,$ by $\,a_{n,1}\,$, $\,\chi\,$ by $\,\bar\chi$, $\,s\,$ by $\,1-s\,$, and $(\l,\d)$ by their images under the outer automorphism (\ref{outer1}); also, the opposite signs of the arguments of the exponential terms in lines a) and b) in \propref{prop2} introduce the additional sign factor $(-1)^{\d_1}$:
\begin{equation}
\label{starmellout4}
\begin{aligned}
{}&{\sum}_{\ell \in (\Z/q\Z)^*}\ \bar\chi(\ell)\, \bigl( M_{\e+\d_2} \sigma_{q,\ell}(y)\bigr) (s+\l_2) 
\\
&\ \ \ \ = \ 2\, (-1)^{\e + \d_1}\, g_{\chi}\, q^{1-2s - \l_1 - \l_2}\, G_{\e + \d_1}(1 -s - \l_1)  \, L(1 - s,\tilde\tau\otimes \bar\chi)\,.
\end{aligned}
\end{equation}
This, too, is the shifted Mellin transform of a tempered distribution which vanishes to infinite order at $y=0$ and extends canonically across zero. Hence $\,G_{\e + \d_1}(s - \l_1) L(s,\tilde\tau\otimes \bar\chi)\,$ is entire and has polynomial growth along vertical strips. According to \lemref{inforder},
\begin{equation}
\label{starmell1.4}
\begin{aligned}
{}&g_\chi^{-1}\,\ {\sum}_{k\in (\Z/q\Z)^*}\ \chi(k)\,\ {\sum}_{n\neq 0}\ c_{n,1}\,e(nx + nk/q) 
\\
{}&= \ {\sum}_{n\neq 0}\ \bar\chi(n)\, c_{n,1}\,e(nx)\ = \ {\sum}_{n\neq 0}\ \bar\chi(n)\, (\sg
n)^{\d_1}\,|n|^{-\l_1}\,a_{n,1}\,e(nx)
\end{aligned}
\end{equation}
vanishes to infinite order at $x=0$. We let $\d_a(x)$ denote the delta function at\break $a \in \R$. Then
$\,\sum_{n\neq 0} \bar\chi(n) (\sg n)^{\d_1}|n|^{-\l_1}a_{n,1}\d_n(x)\,$ is the inverse Fourier transform of
(\ref{starmell1.4}), which therefore extends canonically across $\infty$
\cite[Th.~3.19]{MilSch:2003b}. This tempered distribution also vanishes identically near $x=0$, so that 
\begin{equation}
\label{starmell1.5}
\begin{aligned}
{}&2\,L(s,\tilde\tau\otimes\bar\chi)  
\\
{}&\ \ =\ \int_\R (\sg x)^{\e + \d_1} \,|x|^{\l_1 - s}\,{\sum}_{n\neq 0}\ \bar\chi(n)\, (\sg
n)^{\d_1}\,|n|^{-\l_1}\,a_{n,1}\,\d_n(x)\,dx
\end{aligned}
\end{equation}
is entire, of polynomial growth along vertical strips. Replacing $\tilde\tau\otimes\bar\chi$ by $\tau\otimes\chi$, we get the same conclusion about $L(s,\tau\otimes\chi)$, of course.

The functional equation is an immediate consequence of (\ref{starmellout3}), (\ref{starmellout4}) and
(\ref{starmell1}), coupled with the Gamma identity (\ref{heresG3}) and the relation
(\ref{gausssumidentity}a) between $g_\chi$, $g_{\bar\chi}$, $q$ and $\varepsilon$. 
\Endproof\vskip4pt

The preceding proof deduces the functional equation from the validity
of (\ref{star2}) with $a_k=\chi(k)$.  Reversing the steps we obtain a statement, \lemref{melinvconv} below, which we shall use in the proof of the converse theorem. We fix parameters
\begin{equation}
\label{mlconv3}
\begin{aligned}
{}&\l \ = \ (\l_1,\,\l_2,\,\l_3)\, \in \,\C^3\,, \ \ \ \text{with}\ \ \textstyle{\sum}_{1\leq j \leq 3}\ \l_j\, =\,0\,,
\\
{}&\d \ = \ (\d_1,\,\d_2,\,\d_3)\, \in \,(\Z/2\Z)^3\,, \ \ \ \text{with}\ \ \textstyle{\sum}_{1\leq j \leq 3}\ \d_j\ = \
0\,.
\end{aligned}
\end{equation}
In addition we consider the data of complex coefficients $(a_{1,n})_{n \geq 1}$, $(a_{n,1})_{n\geq 1}$,
which grow at most polynomially in $n$. If $\chi$ is a Dirichlet character, we define twisted Dirichlet series
\begin{equation}
\label{mlconv4}
L_\chi(s)\ =\ {\sum}_{n=1}^\infty\, a_{1,n}\,\chi(n)\,n^{-s},
\ \ \ \ \tilde{L}_{\bar\chi}(s)\ =\ {\sum}_{n=1}^\infty\, a_{n,1}\,\bar\chi(n)\,n^{-s} 
\end{equation}
which converge for $\,\operatorname{Re} s \gg 0\,$. We suppose that these series, corres\-ponding to every
choice of $q>0$ and primitive Dirichlet character modulo $q$, extend to the complex plane and satisfy the
regularity hypotheses
\begin{equation}
\label{mlconv2}
\begin{aligned}
{}&L_\chi(s)\,,\ \ G_{\e+\d_3}(s+\l_3)\,L_\chi(s)\,,\ \ \tilde{L}_{\bar\chi}(s)\ \ \text{and} \ \
G_{\e+\d_1}(s-\l_1)\,\tilde{L}_{\bar\chi}(s)
\\
{}&\text{are entire functions, of polynomial growth on vertical strips},
\end{aligned}
\end{equation}
as well as the functional equations
\begin{equation}
\label{mlconv1}
\tilde L_{\bar\chi}(1-s) \ \ =
 \ \ q^{3s} \, g_\chi^{-3} \, \( {\prod}_{j=1}^3\ G_{\e+\d_j}(s+\l_j)\) L_\chi (s)\, .
\end{equation}
Here $\e \in \Z/2\Z$ is defined by the parity condition $\chi(-1) = (-1)^\e$. We use (\ref{cnorm}) to
determine coefficients $(c_{1,n})_{n\neq 0}$ and $(c_{n,1})_{n\neq 0}$ in terms of the $a_{1,n}$ and
$a_{n,1}$. Recall the definition (\ref{wmudist}) of the representation spaces $W_{\mu,\eta}^{-\infty}$. As
periodic distributions without constant term, the tempered distributions
\begin{equation}
\label{mlconv5}
\tau_{x,1}(x)\ =\ {\sum}_{n\neq 0}\, c_{n,1}\,e(nx)\,,\ \ \ \tau_{y,1}(y)\ =\ {\sum}_{n\neq 0}\, c_{1,n}\,e(ny)
\end{equation}
extend canonically across $\infty$ \cite[Prop.~2.19]{MilSch:2003b}, and thus can be regarded as
vectors in $W_{\l_1-\l_2,\d_3}^{-\infty}$ and $W_{\l_2-\l_3,\d_1}^{-\infty}$, respectively. For $(c,a)=1$
only, we use the first two formulas in \corref{cor1} to define vectors
\begin{equation}
\label{mlconv6}
\sigma_{c,a}\, \in \, W_{\l_1-\l_2,\d_3}^{-\infty}\ ,  \ \ \ \ \ \ \rho_{c,a}\, \in \, W_{\l_2-\l_3,\d_1}^{-\infty}
\end{equation}
in terms of $\tau_{x,1}$ and $\tau_{y,1}$. By construction, both vanish to infinite order at $0$.

\begin{lem}\label{melinvconv} Let $\chi$ be a primitive Dirichlet character of conductor $q>0$.
 Under the hypotheses just stated{\rm ,} the identity {\rm (\ref{star2})} holds with $a_k=\chi(k)${\rm ,}
 as an identity of
tempered distributions on $\,\R$. Both $\,\sum_{k\in(\Z/q\Z)} \chi(k)\,\rho_{q,k}$ and
$\,\sum_{\ell\in(\Z/q\Z)} \bar\chi(\ell)\,\sigma_{q,\ell}$ vanish to infinite order at the point at infinity.
\end{lem}

At the very beginning of the proof of Theorem \ref{funeqthm}, we   remarked that $\widehat a_\ell =
g_\chi\, \bar\chi(\ell)$, $\ell\in \Z/q\Z\,$, corresponds to $a_k = \chi(k)$, $k\in \Z/q\Z\,$, via the finite
Fourier transform (\ref{fourierint3}). Thus (\ref{star2}) has meaning in the present context: $\rho_{q,k}$
and $\sigma_{q,\ell}$ occur with nonzero coefficients only when they are well defined, i.e., only when
$(q,k)=(q,\ell)=1$.

\Proof  In view of (\ref{starmell1.3}),
\begin{equation}
\label{melinvconv1}
\begin{aligned}
{}&{\sum}_{k \in (\Z/q\Z)^*} \ \chi(k) \,\tau_{y,1}(y-\bar k/q)  
\\
{}&\qquad = \ {\sum}_{k \in (\Z/q\Z)^*} \ \chi(k) \,\sum_{n\neq 0}\ a_{1,n}\, |n|^{\lambda_3}\, (\sg
n)^{\delta_3}\, e(ny-n\bar k/q)
\\
{}&\qquad  = \  (-1)^\varepsilon g_{\bar\chi}\, {\sum}_{n\neq 0} \ \chi(n) \, a_{1,n}\, |n|^{\lambda_3}\, (\sg
n)^{\delta_3}\, e(ny) \,,
\end{aligned}
\end{equation}
so $L_\chi(s-\lambda_3)$ is the Dirichlet series formed from the Fourier coefficients of this periodic distribution. It has parity $\varepsilon + \delta_3$, since $\chi(-n)=(-1)^\varepsilon \chi(n)$. At this point the hypothesis (\ref{mlconv2}) allows us to apply
 \cite[Prop.~5.37]{MilSch:2003b}, with $a_n = \chi(n)\, a_{1,n}$, $\,\nu=-\lambda_3$, and
$\,\delta=\varepsilon+\delta_3$, and conclude that
\begin{equation}
\label{melinvconv2}
{\sum}_{k \in (\Z/q\Z)^*} \ \chi(k) \,\tau_{1,y}(y-\bar k/q)\ \ \ \text{vanishes to infinite order at}\ \ y=0\,.
\end{equation}
The definition of $\rho_{q,k}$ relates the behavior of $\tau_{1,y}(y-\bar k/q)$ to that of $\rho_{q,k}(y)$ at $y=\infty$, and we already know that $\rho_{q,k}$ vanishes to infinite order at $y=0$. Arguing analogously in the case of the $\sigma_{n,k}$, we find
\begin{equation}
\label{melinvconv3}
\begin{aligned}
{}&{\sum}_{k \in (\Z/q\Z)^*} \ \chi(k) \,\rho_{q,k}(y)\ \ \ \text{and}\ \ \ {\sum}_{\ell \in (\Z/q\Z)^*} \
\bar\chi(\ell) \,\sigma_{q,\ell}(y)
\\
{}&\text{vanish to infinite order both at}\ \ y=0\ \ \text{and at}\ \ y=\infty\,.
\end{aligned}
\end{equation}
The sign factors in the defining relations for the $\rho_{q,k}$ and $\sigma_{q,k}$ show that both have parity $\varepsilon + \delta_2$. An even or odd tempered distribution which vanishes to infinite order at $0$ and extends canonically across $\infty$ is completely determined by its signed Mellin transform
\cite[Th.~4.8]{MilSch:2003b}. It follows that (\ref{ahatk})  -- i.e., the identity asserted by the lemma -- is
equivalent to the first half of (\ref{starmell1}). In effect, the proof of theorem \ref{funeqthm} shows that the
equality of the outer terms in (\ref{starmell1}) is equivalent to the functional equation (\ref{mlconv1}). To
complete the proof of the lemma, all that remains to be established is the second half of (\ref{starmell1}),
i.e.,
\begin{equation}
\label{melinvconv4}
\begin{aligned}
{}&(-1)^{\e+\d_2}\, {\sum}_{\ell \in (\Z/q\Z)^*} \ \bar\chi(\ell)\, \bigl(
M_{\varepsilon+\delta_2}\widehat\sigma_{q,\ell}\bigr)\,(1-s-\l_2) 
\\
{}&\ \ \ \ \ \ = \ G_{\e+\d_2}(1-s-\l_2)\  {\sum}_{\ell \in (\Z/q\Z)^*} \ \bar\chi(\ell) \, \bigl(
M_{\varepsilon+\delta_2}\sigma_{q,\ell}\bigr)\,(s+\l_2)
 \,.
\end{aligned}
\end{equation}
This follows from an application of \cite[Th.~4.12]{MilSch:2003b}.
\hfill\qed

\section{The converse theorem}\label{convsec}

The functional equation for a cuspidal automorphic representation involves only the coefficients $\,a_{n,1}\,$, $\,a_{1,n}\,$, but all the $\,a_{m,n}\,$ are needed to reconstruct the automorphic representation. As is well known, Hecke theory provides the bridge between the two types of data. Since Hecke operators for ${\rm GL}(n)$ have upper triangular representatives, they act transparently on automorphic distributions. We briefly review the relevant formulas, not only as a tool for the proof of the converse theorem, but also to identify our definition of the $L$-functions $L(s,\tau\otimes\chi)$ with the usual one. Details on the Hecke algebra can be found in
 \cite[\S3]{Shimura:1971} and \cite[\S9]{Bump:1984}, for example.

The Hecke algebra for $\,\G = {\rm GL}(3,\Z)\,$ acts on the $\,\G$-invariants in any representation space for ${\rm GL}(3,\R)$, in particular on $\,(V_{\l,\d}^{-\infty})^\G$. It is generated by operators $T_{1,p}$ and $T_{p,1}$ attached to each prime $p$. To shorten the formulas, we introduce the matrix
\begin{equation}
\label{habc}
h_{\,a,\,b,\,c;\,i,\,j,\,\ell}\ \ = \ \ {\ttt{a}{i}{\ell}{}{b}{j}{}{}{c}}^{\!-1} \ \in \ {\rm GL}(3,\Q)
\end{equation}
whose inverse has entries $\,a,\, b,\, c \in \Z_{>0}\,$\,. Then, for $\,\tau \in (V_{\l,\d}^{-\infty})^\G$,
\begin{equation}
\label{t1pp1def}
\begin{aligned}
{}&T_{1,p}\tau\ \ = \ \ p^{-1}\biggl( \pi_{\l,\d}\(\,h_{\,p,\,1,\,1;\,0,\,0,\,0} \)\tau \   \biggr.
\\
{}&\ \  + \biggl.\ {\sum}_{0\leq i <p}\, \pi_{\l,\d}\(\,h_{\,1,\,p,\,1;\,i,\,0,\,0}\)\tau\ + \  {\sum}_{0 \leq \,j,\,\ell
<p}\,
\pi_{\l,\d}\(\,h_{\,1,\,1,\,p;\,0,\,j,\,\ell}\)\tau \biggr),
\\
{}&T_{p,1}\tau\ \ = \ \ p^{-1}\biggl( \pi_{\l,\d}\(\,h_{\,p,\,p,\,1;\,0,\,0,\,0}\)\tau \    \biggr.
\\
{}&\ \  + \ \biggl.  {\sum}_{0 \leq j < p}\, \pi_{\l,\d}\(\,h_{\,p,\,1,\,p;\,0,\,j,\,0}\)\tau\ + \ {\sum}_{0\leq\,
i,\,\ell < p}\, \pi_{\l,\d}\(\,h_{\,1,\,p,\,p;\,i,\,0,\,\ell}\)\tau \biggr).
\end{aligned}
\end{equation}
Conjugation by the outer automorphism (\ref{outer1}) interchanges $T_{p,1}$ and $T_{1,p}$. They commute with each other and with the operators corresponding to other primes. We shall also need the formulas for the action of the Hecke operators in terms of the coordinates (\ref{coordinates}) on $N$,
\begin{equation}
\label{t1pp1}
\begin{aligned}
(T_{1,p}\tau)(x,y,z)\ = \ p^{-\l_1}\,\tau(px,y,pz) \, &+  \, p^{-\l_2-1}\!\!\! \sum_{i\in \Z/p\Z}\! \tau(\textstyle\f{x+i}{p},py,z+iy)
\\
&+ \, p^{-\l_3-2} \!\! \sum_{j,\ell\in \Z/p\Z}\! \tau(x,\textstyle\f{y+j}{p},\f{z+\ell}{p})\ ,
\\
(T_{p,1}\tau)(x,y,z) \ = \ p^{\l_3}\,\tau(x,py,pz) \, &+ \, p^{\l_2-1}\!\!\! \sum_{j\in \Z/p\Z}\! \tau(px,\textstyle\f{y+j}{p},z)
\\
&+ \, p^{\l_1-2} \!\! \sum_{i,\ell\in \Z/p\Z}\! \tau(\textstyle\f{x+i}{p},y,\f{z+\ell+i y}{p})\ ,
\end{aligned}
\end{equation}
which is a straightforward consequence of (\ref{action}).

Conjugation by any matrix in ${\rm GL}(3,\Q)$ maps $\G$ to a subgroup of ${\rm GL}(3,\Q)$ that is commensurate to $\G$. Using this observation one can show that the action of the Hecke algebra preserves the space of all cuspidal automorphic distributions $\,\tau \in (V_{\l,\d}^{-\infty})^\G$. If $v \in V_{-\l,\d}$ is an ${\rm SO}(3)$-finite vector, the map which assigns to any cuspidal $\tau \in (V_{\l,\d}^{-\infty})^\G$ the $L^2$ function $g \mapsto \langle \pi_{-\l,\d}(g)v , \tau \rangle$ is equivariant with respect to the Hecke actions on $L^2(\G\backslash {\rm GL}(3,\R))$ and on $(V_{\l,\d}^{-\infty})^\G$. This observation, in conjunction with the usual arguments, makes it possible to simultaneously diagonalize the action of the $T_{1,p}$ and $T_{p,1}$, corres\-ponding to all primes $p$, on the space of cuspidal automorphic distri\-butions $\tau \in (V_{\l,\d}^{-\infty})^\G$. Let us suppose then that $\tau\! \in \!(V_{\l,\d}^{-\infty})^\G$ is a cuspidal Hecke eigendistribution, say $T_{1,p}\tau!
  = \mu_{1,p}\tau$ and $T_{p,1}\tau = \mu_{p,1}\tau$ for all primes $p$. From (\ref{t1pp1}) one sees that $T_{1,p}$ and $T_{p,1}$ operate separately on the abelian and the non-abelian part of the Fourier decomposition of $\tau$, as in \propref{prop1}. Because of \corref{cor1} and \lemref{inforder}, the abelian Fourier coefficients $c_{r,s}$ completely determine $\tau$. Recall that the $a_{r,s}$, $1 \leq r,s < \infty$, determine the $c_{r,s}$ and are determined by them;
 cf.\ (\ref{cnorm}), (\ref{awitheps}). A straightforward analysis of (\ref{t1pp1}) lets us conclude that
$a_{1,p} = \mu_{1,p}a_{1,1}$ and $a_{p,1} = \mu_{p,1}a_{1,1}$, and that if these coefficients vanish,
then so do all the others. In particular, $a_{1,1}$ cannot vanish unless $\tau = 0$. We can therefore
renormalize $\tau$ so that $a_{1,1}=1$; in that case, $\mu_{1,p} = a_{1,p}$, $\mu_{p,1} = a_{p,1}$. To
summarize,
\begin{equation}
\label{a1pp1}
T_{1,p}\tau \ = \ a_{1,p}\,\tau  \,, \qquad T_{p,1}\tau \ = \ a_{p,1}\,\tau
\end{equation}
when $\tau$ is a cuspidal Hecke eigendistribution, normalized so that $a_{1,1}=1$. In this situation,
\begin{equation}
\label{anmmult}
\begin{aligned}
a_{1,p}\, a_{r,\,s} \ \ &=  \ \ a_{r\!/p,\,s} \ + \  a_{rp,\,s\!/p}\ + \ a_{r\!,\,sp}\,,
\\
a_{p,1}\, a_{r,\,s} \ \ &= \ \ a_{r,\,s\!/p} \ + \ a_{r\!/p,\,sp} \ + \ a_{rp,\,s}\,,
\end{aligned}
\end{equation}
subject to the convention that coefficients with a nonintegral index $r/p$ or $s/p$ are to be omitted. These identities can be deduced from (\ref{t1pp1}). The Euler factorizations
\begin{equation}
\label{eulerprod3}
\begin{aligned}
{\sum}_{n=1}^\infty\ a_{1,n}\,n^{-s} \  &= \, \  {\prod}_p \ (1\,-\,a_{1,p}\,p^{-s}\,+\,a_{p,1}\,p^{-2s}\,-\,p^{-3s})\i\,,
\\
{\sum}_{n=1}^\infty\ a_{n,1}\,n^{-s} \  &= \, \  {\prod}_p \ (1\,-\,a_{p,1}\,p^{-s}\,+\,a_{1,p}\,p^{-2s}\,-\,p^{-3s})\i\,,
\end{aligned}
\end{equation}
for $\,\operatorname{Re} s \gg 0\,$, then follow formally from (\ref{anmmult}). As   mentioned earlier,
(\ref{a1pp1}) and (\ref{eulerprod3}) imply in particular that our definition (\ref{stanlfuncdef}) of the
standard $L$-function agrees with the usual one. Recursive use of (\ref{anmmult}) leads to an explicit
formula for the $a_{r,s}$ in terms of the $a_{r,1}$ and $a_{1,s}$:
\begin{equation}
\label{anmrecurs}
a_{r,s} \ \ = \ \ {\sum}_{d|(r,s)}\ \mu(d)\ a_{r\!/d,\,1}\ a_{1,\,s\!/d} \,;
\end{equation}
here $\mu(\dots)$ denotes the M\"obius $\mu$-function,
\begin{equation}
\label{mudef}
\mu(p_1^{e_1}p_2^{e_2}\cdots p_r^{e_r})\ \ = \ \ \begin{cases}\ (-1)^r &\text{if all}\ \, e_j=1\ \text{and}\, \ p_i\neq p_j\ \text{for}\ i\neq j\,,
\\
\ 0 &\text{if at least one}\, \ e_j>1\,.\end{cases}
\end{equation}
A detailed discussion of the logical connections between (\ref{anmmult})--(\ref{anmrecurs}) can be found in
\cite[\S9]{Bump:1984}, for example.

For the statement of the converse theorem, we fix parameters $\,\l\,$, $\,\d\,$ as in (\ref{mlconv3}) and two Dirichlet series
\begin{equation}
\label{mlconv7}
L(s)\ =\ {\sum}_{n=1}^\infty\, a_{1,n}\,n^{-s},\ \ \ \ \tilde{L}(s)\ =\ {\sum}_{n=1}^\infty\, a_{n,1}\,n^{-s}\,,
\end{equation}
which we suppose converge for $\,\operatorname{Re} s \gg 0\,$. If $\chi$ is a primitive Dirichlet character with conductor $q>0$, we define the twisted $L$-functions $L_\chi$ and $\tilde{L}_{\bar\chi}$ as in (\ref{mlconv4}); these too converge for $\,\operatorname{Re} s \gg 0\,$. We suppose further that $L_\chi$ and $\tilde{L}_{\bar\chi}$ have meromorphic continuations to the entire complex plane, for every primitive Dirichlet character $\chi$.

\begin{thm}[${\rm GL}(3)\times {\rm GL}(1)$ converse theorem -- see
\cite{JacPiatShal:1979}] \label{conversethm} The function $L(s)$ is the standard $L$-function of a cuspidal
${\rm GL}(3,\Z)$-automorphic representation with scalar Hecke action provided\/{\rm :}\/
\begin{itemize}
\ritem{a)}  Both $L(s)$ and $\tilde L(s)$ have Euler product factorizations as in {\rm (\ref{eulerprod3}),}
and 
\ritem{b)} For every primitive Dirichlet character $\chi$ of conductor $q>0${\rm ,} the twisted
$L$-functions
$\,L_\chi\,${\rm ,} $\,\tilde{L}_{\bar\chi}\,$ satisfy the analytic hypotheses {\rm (\ref{mlconv2})} and the
functional equation {\rm (\ref{mlconv1})}.\end{itemize}
\end{thm}

The proof of the theorem occupies the remainder of this section. The idea is to construct a cuspidal automorphic distribution $\,\tau \in (V_{\l,\d}^{-\infty})^\G$ from the
data of the Dirichlet series (\ref{mlconv7}).

The identity (\ref{anmrecurs}) defines coefficients $\,a_{r,s}\,$, $\,r,\,s \geq 1\,$, in terms of the coefficients of the Dirichlet series (\ref{mlconv7}). The convergence of these series for $\,\operatorname{Re} s \gg 0\,$ implies that the $\,a_{1,n}\,$ and $\,a_{n,1}\,$ grow at most polynomially with $n$, which in turn bounds the growth of the $\,a_{r,s}\,$ in terms of powers of the indices. Equivalently the $\,a_{r,s}$ can be characterized by the identity
\begin{equation}
\label{double_dirichlet}
{\sum}_{m,n\geq 1}\ a_{m,n}\,m^{-s_1}\,n^{-s_2}\ = \ \frac{L(s_2)\, \tilde L(s_1)}{\zeta(s_1+s_2)}\qquad(\,\operatorname{Re} s_1,\,\operatorname{Re} s_2\gg 0\,)\,,
\end{equation}
which encodes the defining relation (\ref{anmrecurs}) for the $a_{r,s}$ \cite[\S9]{Bump:1984}. The Euler
pro\-ducts for $L(s_2)$, $\tilde L(s_1)$ and $\zeta(s_1+s_2)$ determine an Euler product for the double
Dirichlet series (\ref{double_dirichlet}). This latter Euler product implies that
\begin{equation}
\label{anmrecurs1}
\text{the $\,a_{r,s}\,$ satisfy the four term relations (\ref{anmmult})}\,,
\end{equation}
as can be seen by forming the double Dirichlet series from the four term relations and rearranging terms.

We now use (\ref{cnorm}) to introduce coefficients $\,c_{r,s}\,$, for all nonzero $\,r,\,s$. Like the $\,a_{r,s}\,$ they grow at most polynomially in the indices. Hence
\begin{equation}
\label{tauab}
\tau_{\text{abelian}}(x,y,z)\ \ = \ \ {\sum}_{r,\,s\, \neq\, 0}\ \,c_{r,s}\ e(rx+ sy) \ \in \ C^{-\infty}(N_\Z\backslash N)
\end{equation}
is a well defined $N_\Z$-invariant distribution on $N$. This is the abelian part of the automorphic distribution $\tau$ we want to construct. It is an eigendistribution of the Hecke operators (\ref{t1pp1}):
\begin{equation}
\label{a1pp1ab}
T_{1,p}\tau_{\text{abelian}} \ = \ a_{1,p}\,\tau_{\text{abelian}}  \,, \qquad T_{p,1}\tau_{\text{abelian}} \ = \ a_{p,1}\,\tau_{\text{abelian}}\,,
\end{equation}
since these identities are formally equivalent to the four term relations (\ref{anmmult}).

Our next task is to construct the non-abelian Fourier components. Because of 
\cite[Prop.~2.19]{MilSch:2003b}, the tempered distributions
\begin{equation}
\label{tauxqyq}
\tau_{x,q}(x) \ = \ {\sum}_{r\neq 0}\ c_{r,q}\,e(rx)\,,\ \ \ \  \tau_{y,q}(y) \ = \ {\sum}_{s\neq 0}\ c_{q,s}\,e(sy)
\end{equation}
have canonical extensions across infinity. Hence \cite[Prop.~2.26]{MilSch:2003b} allows us to regard them
as distribution vectors,
\begin{equation}
\label{tauxdyd}
\tau_{x,q}\, \in \, W_{\l_1-\l_2,\d_3}^{-\infty}\ ,  \ \ \ \ \ \ \tau_{y,q}\, \in \, W_{\l_2-\l_3,\d_1}^{-\infty}\,;
\end{equation}
cf.\ (\ref{wmudist}). By construction, $\tau_{x,q}(x)$ and $\tau_{y,q}(y)$ vanish to infinite order at $x=\infty$ and $y=\infty$, respectively. For $\,a,\, c,\, q \in \Z\,$, $\,c\neq 0\,$, $\,q> 0\,$ and $\,(a,c)=1$, we use the first two formulas in \corref{cor1} to define
\begin{equation}
\label{rhonksigmank}
\sigma_{cq,aq}\, \in \, W_{\l_1-\l_2,\d_3}^{-\infty}\ ,  \ \ \ \ \ \ \rho_{cq,aq}\, \in \, W_{\l_2-\l_3,\d_1}^{-\infty}\ .
\end{equation}
When $q=1$, the definitions (\ref{tauxqyq})--(\ref{rhonksigmank}) agree with
(\ref{mlconv5}), (\ref{mlconv6}). The hypotheses of the converse theorem include those of
\lemref{melinvconv}, so that  its conclusions apply in the present setting. The element of ${\rm SL}(2,\Q)$
that relates $\sigma_{cq,aq}$ to $\tau_{x,q}$ and
$\rho_{cq,aq}$ to $\tau_{y,q}$ depends only on $a$ modulo $c$, not on the particular choice of $a$. This
allows us to regard the distribution vectors (\ref{rhonksigmank}) as indexed by $n=cq$ and $k=aq \in
\Z/n\Z$, as in the statement of \propref{prop1}.

We shall now reduce the proof of the converse theorem to a sequence of five lemmas about the formal sums
\begin{equation}
\label{tau1tau2}
\begin{aligned}
\tau^{(1,2)}&(x,y,z)\ =\ \tau_{\text{abelian}}(x,y,z) 
\\
&+\ {\sum}_{n\neq 0} {\sum}_{k \in \Z/n\Z} {\sum}_{\ell\equiv k\, (\operatorname{mod}n)} \ e(nz+\ell y)\,\sigma_{n,k}(x+\ell/n)
\\
\tau^{(2,3)}&(x,y,z)\ =\ \tau_{\text{abelian}}(x,y,z) 
\\
&\!\!\!\!\!\!\! +\ {\sum}_{n\neq 0} {\sum}_{k \in \Z/n\Z} {\sum}_{\ell\equiv k\, (\operatorname{mod}n)} \ e(n(z-xy)+\ell x)\,\rho_{n,k}(\ell/n-y)\,.
\end{aligned}
\end{equation}
Recall the definition of the codimension one Schubert cell $\{x=\infty\}$ in the proof of \lemref{inforder}. Its union with the open Schubert cell $N$ is stable under the parabolic subgroup $P^{(1,2)} \subset {\rm GL}(3,\R)$ generated by $N$ and the copy of ${\rm GL}(2,\R)\times {\rm GL}(1,\R)$ which stabilizes the decomposition $\R^3 = \R^2 \oplus \R$. Similarly $N \cup \{y=\infty\}$ is stable under the parabolic subgroup $P^{(2,3)}$ generated by ${\rm GL}(1,\R)\times {\rm GL}(2,\R)$ and $N$. The space $V_{\lambda,\delta}^{-\infty}$ was defined as the space of distribution sections of the $G$-equivariant line bundle $\mathcal L_{\lambda,\delta} \to G/MAN_-$. Further notation: $P^{(1,2)}_\Z = P^{(1,2)}\cap \Gamma$ and $P^{(2,3)}_\Z = P^{(2,3)}\cap \Gamma$.

\begin{lem}\label{lemconv} The series $\tau^{(1,2)}$ converges in the strong distribution topology to a $P^{(1,2)}_\Z$-invariant distribution section of $\mathcal L_{\lambda,\delta}$ over
 $N\cup\{x=\infty\}$.\break Similarly $\tau^{(2,3)}$ converges to a $P^{(2,3)}_\Z$-invariant distribution
section over $N\cup\{y=\infty\}$.
\end{lem}

We do not make a notational distinction between the series $\tau^{(1,2)}$, $\tau^{(2,3)}$ and their limits. Note that the Hecke operators (\ref{t1pp1def}) have representatives in $P^{(1,2)}\cap P^{(2,3)}$, and thus can be applied to sections of $\mathcal L_{\lambda,\delta}$ over $N\cup\{x=\infty\}$ or $N\cup\{y=\infty\}$.

\begin{lem}\label{lemhecke} Both $\tau^{(1,2)}$ and $\tau^{(2,3)}$ are eigendistributions
 of the Hecke action\/{\rm :}\/ for all primes $p${\rm ,}
\begin{eqnarray*}
T_{1,p}\tau^{(1,2)}& =& a_{1,p}\tau^{(1,2)}\,,\ \ \ \ T_{p,1}\tau^{(1,2)}\ =\ a_{p,1}\tau^{(1,2)}\,,
\\
 T_{1,p}\tau^{(2,3)}& =& a_{1,p}\tau^{(2,3)}\,,\ \ \ \ T_{p,1}\tau^{(2,3)}\ =\ a_{p,1}\tau^{(2,3)}\,.
\end{eqnarray*}
\end{lem}

\begin{lem}\label{lemcoinc} $\tau^{(1,2)}$ and $\tau^{(2,3)}$ agree on their common domain $N$.
\end{lem}

\begin{lem}\label{leminforder} $\tau^{(1,2)}$ vanishes to infinite order along $\{x=\infty\}${\rm ,}  and
$\tau^{(2,3)}$ vanishes to infinite order along $\{y=\infty\}$.
\end{lem}

\begin{lem}\label{lemexist} There exists a unique $\,\tau\in (V_{\lambda,\delta}^{-\infty})^\Gamma$ which agrees with $\,\tau^{(1,2)}$ and $\,\tau^{(2,3)}$ on their respective domains. This distribution vector $\,\tau$ is cuspidal.
\end{lem}

Assuming the statements of the five lemmas, we can easily complete the proof of the converse theorem. The Hecke operators map
$\Gamma$-invariants to\break $\Gamma$-invariants. Thus, invoking Lemmas \ref{lemhecke} and
\ref{lemexist}, we find that $T_{1,p}\tau - a_{1,p}\tau$ is $\Gamma$-invariant and vanishes on the open
Schubert cell $N$, whose $\Gamma$-translates cover $G/MAN_-$. This difference must therefore vanish
globally. We argue similarly in the case of $T_{p,1}$ and conclude
\begin{equation}
\label{hecketau}
T_{1,p}\tau\ = \ a_{1,p}\tau\,,\qquad T_{p,1}\tau\ = \ a_{p,1}\tau\,,
\end{equation}
for all primes $p$. Thus $\tau$ is a cuspidal $\G$-invariant Hecke eigendistribution, as asserted by the theorem. We now turn to the proofs of the five lemmas.

\demo{Proof Lemma {\rm \ref{lemconv}}}  We shall only treat the case of $\tau^{(1,2)}$. One can argue
analogously in the case of $\tau^{(2,3)}$, or use the fact that the two series are related by the outer
automorphism (\ref{outer1}). We argue as in the proof of Lemma \ref{inforder} to extend the coordinates
$(x,y,z)$ from
$N$ to $N\cup\{x=\infty\}$, with $x$ replaced by its reciprocal along $\{x=\infty\}$, of course. This allows
us to identify sections of $\mathcal L_{\lambda,\delta}$ over $N\cup\{x=\infty\}$ with functions on
$N_{y,z}$ which take values in $W_{\lambda_1 - \lambda_2,\d_3}$; cf.\
(\ref{wmuinf1})--(\ref{wmudist}). Because of (\ref{tauxqyq})--(\ref{tau1tau2}), and
the definition of $\sigma_{cq,aq}$ in terms of the first identity in \corref{cor1}, the sum
\begin{equation}
\label{tau1}
\begin{aligned}
{}&\tau^{(1,2)}(x,y,z)\ =\ {\sum}_{q\neq 0}\ e(qy)\,\tau_{x,q}(x) 
\\
{}&\ \ +\ {\sum}_{q>0} {\sum}_{\stackrel{\scriptstyle{(a,c)=1}}{c\neq 0}} \ e(cqz+aq
y)\biggl(\!\psi_{\lambda_1 - \lambda_2,\d_3}\ttwo{1}{-a/c}{}{1}\sigma_{cq,aq}\!\biggr)\!(x)
\\
{}&=\ {\sum}_{q\neq 0}\ e(qy)\,\tau_{x,q}(x) 
\\
{}&\ \ +\ {\sum}_{q>0} {\sum}_{\stackrel{\scriptstyle{(a,c)=1}}{c\neq 0}} \ e(cqz+aq
y)\biggl(\!\psi_{\lambda_1 - \lambda_2,\d_3}\ttwo{a}{-b}{-c}{\bar a}\tau_{x,q}\!\biggr)\!(x)
\end{aligned}
\end{equation}
formally defines a distribution section of $\mathcal L_{\lambda,\delta}$ over $N\cup\{x=\infty\}$. On the right of this identity, $\bar a$ and $b$ are to be chosen subject to the condition $a\bar a -bc=1$; the particular choice does not matter since $\tau_{x,q}(x)$ is invariant under integral translations.

We establish convergence separately for ${\sum}_{q\neq 0}\, e(qy)\,\tau_{x,q}(x)$ and the remaining sum. When we regard the $\tau_{x,q}$ as scalar distributions, as we may,
\cite[Prop.~2.19]{MilSch:2003b} asserts that ${\sum}_{q\neq 0}\, e(qy)\,\tau_{x,q}(x)$ converges to a
distribution on $(\R\cup\{\infty\})\times \R/\Z$ which vanishes to infinite order along $\{\infty\}\times
\R/\Z$. The interpretation of the $e(qy)\,\tau_{x,q}(x)$ as distribution sections of $\mathcal
L_{\lambda,\delta}$ over $N\cup\{x=\infty\}$ involves multiplication by $(\sg x)^{\d_3}|x|^{\l_1-\l_2-1}$
and the introduction of the third coordinate $z$, which affect neither convergence of the series, nor vanishing
to infinite order along $x=\infty$ \cite[Prop.~2.26]{MilSch:2003b}. For future reference, we record
the observation that
\begin{equation}
\label{convergence1}
{\sum}_{q\neq 0}\ e(qy)\, \tau_{x,q}(x)\ \ \text{vanishes to infinite order along}\ \ \{x=\infty\}\,,
\end{equation}
as distribution section of $\mathcal L_{\lambda,\delta}$. Thus $\tau_{\text{abelian}}$ -- which coincides with the restriction of ${\sum}_{q\neq 0}\,e(qy)\, \tau_{x,q}(x)$ to $N$ -- extends canonically from $N$ to $N\cup\{x=\infty\}$.

The topology of $W_{\lambda_1 - \lambda_2,\d_3}^{-\infty}$ can be described as the inductive limit of the $C^{-k}$ topologies, $k\geq 0$, on the space of distribution sections of class $-k$, i.e., locally representable as $k$-th derivatives of continuous sections. Since $\R\cup\{\infty\}$ is compact, the topology on the space of $C^{-k}$ sections is definable by a Banach norm $\|\ \|_{-k}$. As was pointed out before, the $c_{m,n}$ grow polynomially in terms of the indices. This implies that all the $\tau_{x,q}$ are $C^{-k}$ sections for a common choice of $k$, and
\begin{equation}
\label{convergence2}
\|\tau_{x,q}\|_{-k}\ \ = \ \ O(\,|q|^{m_1})\qquad \text{for some $m_1\geq 0$.}
\end{equation}
The definition of the action $\psi_{\lambda_1 - \lambda_2,\d_3}$ in terms of the translation action of ${\rm SL}(2,\R)$ on $\R\cup\{\infty\}$ implies that the operators
\begin{equation}
\label{convergence3}
\psi_{\lambda_1 - \lambda_2,\d_3}(g) : W_{\lambda_1 - \lambda_2,\d_3}^{-\infty}\ \longrightarrow\ W_{\lambda_1 - \lambda_2,\d_3}^{-\infty}\qquad \bigl( g \in {\rm SL}(2,\R) \bigr)
\end{equation}
distort the norm $\|\ \|_{-k}$ by a factor which can be bounded in terms of some power of the matrix norm of $g$. For fixed $a,\,c \in \Z$ with $(a,c) =1$, the remaining integer entries $\bar a$, $-b$ of the matrix on the right in (\ref{tau1}) can be chosen so that the matrix norm becomes bounded by a multiple of $\,a^2 + c^2$. Hence
\begin{equation}
\label{convergence4}
\bigl\|\,\psi_{\lambda_1 - \lambda_2,\d_3}\ttwo{a}{-b}{-c}{\bar a} \tau_{x,q}\,\bigr\|_{-k}\ \ \leq \ \ C\,|q|^{m_1}\,(a^2+c^2)^{m_2}\,,
\end{equation}
for some $C>0$ and $m_1,\,m_2\in\N$. Since $\bigl((\frac{\partial\ }{\partial y})^2 + (\frac{\partial\ }{\partial z})^2\bigr)e(cqz+aqy)$ acts on $e(cqz+aqy)$ as multiplication by $-4\,\pi^2\,q^2\,(a^2 + c^2)$, the double sum indexed by $q$ and $(a,c)$ in (\ref{tau1}) converges in the strong distribution topology.

To establish the $P^{(1,2)}_\Z$-invariance of the sum (\ref{tau1}), 
we recall first of all that $\sum_{q\neq 0}\,e(qy)\tau_{x,q}(x)$ is the canonical extension of
$\tau_{\text{abelian}}(x,y,z)$ across\break $\{x=\infty\}$. Also,
\begin{equation}
\label{convergence6}
e(cqz+aq y)\bigl(\psi_{\lambda_1 - \lambda_2,\d_3}\!\ttwo{a}{-b}{-c}{\bar a}\tau_{x,q}\!\bigr)(x)\, = \, \pi_{\lambda,\delta}(\gamma)(e(qy)\,\tau_{x,q}(x))
\end{equation}
as an identity of the  distribution section of $\mathcal L_{\lambda,\delta}$ over $N\cup\{x=\infty\}$, where
\begin{equation}
\label{convergence7}
\gamma\ = \ttt{a}{-b}{}{-c}{\bar a}{}{}{}{1}\ .
\end{equation}
Thus $\,\tau^{(1,2)}$ can be identified with the sum of the $\,\gamma$-translates of the canonical extension of $\,\tau_{\text{abelian}}\,$, with $\,\gamma\,$ running over all the matrices (\ref{convergence7}) indexed by integers $a\,,c$ with $c\neq 0$ and $(c,a)=1$, in addition to $\gamma=e$. These constitute a complete set of representatives for 
$P^{(1,2)}_\Z/M  N_\Z$; recall (\ref{minpar}). The $P^{(1,2)}_\Z$-invariance of $\tau^{(1,2)}$ follows
because $\,\tau_{\text{abelian}}\,$ is $M$- and $N_\Z$-invariant, which also makes its canonical extension
$M$- and 
 $N_\Z$-invariant \cite[Lemma 2.8]{MilSch:2003b}.
\hfill\qed

\demo{Proof  of  Lemma {\rm \ref{lemhecke}}} Recall the definitions (\ref{habc}), (\ref{t1pp1def}) of the Hecke
operators $T_{1,p}$, $T_{p,1}$. The subsets
\begin{equation}
\label{habc2}
\begin{aligned}
{}&S_{1,p}\! = \{h_{p,1,1;0,0,0}\} \cup \{h_{1,p,1;i,0,0} \mid 0\!\leq \! i\! <\! p\} \cup \{h_{1,1,p;0,j,\ell}
\mid 0\!\leq\! j,\ell\! < p\},
\\
{}&S_{p,1}\! = \{h_{p,p,1;0,0,0}\} \cup \{h_{p,1,p;0,j,0} \mid 0\!\leq\! j\! < \!p\} \cup \{h_{1,p,p;i,0,\ell}
\mid 0\!\leq\! i,\ell\! < p\},
\end{aligned}
\end{equation}
of ${\rm GL}(3,\Q)$ satisfy the following two conditions:
\begin{equation}
\label{habc3}
\begin{aligned}
\text{a)}&\qquad g_1,\, g_2 \in S\ \ \text{and}\ \ g_1\G = g_2\G\ \ \Longrightarrow \ \ g_1 = g_2\,,
\\
\text{b)}&\qquad \G \bigl( \cup_{g\in S} \ g\G \bigr)\ = \  \cup_{g\in S} \ g\G \,,
\end{aligned}
\end{equation}
with either $S=S_{1,p}$ or $S=S_{p,1}$ \cite[\S 9]{Bump:1984}. These formally imply that the Hecke
operators preserve the space of $\G$-invariants. We claim:
\begin{equation}
\label{habc4}
\begin{aligned}
{}&\text{The conditions (\ref{habc3}), for both $S=S_{1,p}$ and $S=S_{p,1}$, remain}
\\
{}&\text{correct when $\G$ is replaced by either $\,N_\Z\,$ or $\,P^{(1,2)}_\Z=P^{(1,2)}\cap \G\,$}.
\end{aligned}
\end{equation}
This is obvious for the condition a), since $\G \supset P^{(1,2)}_\Z \supset N_\Z$. The group $N_\Z$ is generated by
\begin{equation}
\label{habc5}
\begin{aligned}
n_1 = \ttt{1}{1}{}{}{1}{}{}{}{1}\ ,\ \ n_2 = \ttt{1}{}{}{}{1}{1}{}{}{1}\ ,\ \ n_3 = \ttt{1}{}{1}{}{1}{}{}{}{1}\ .
\end{aligned}
\end{equation}
It is not difficult to show that conjugation by any one of these reshuffles the cosets $gN_\Z$, $g \in S$ -- in this connection, it helps to observe that in (\ref{habc2}), one can let the indices $i,\,j,\,\ell$ run over any complete set of representatives of $\Z/p\Z$ without changing the collection of right $N_\Z$-cosets generated by the members of the resulting set $S$. This establishes b) in the case of $N_\Z$. Modulo $N_\Z$, $\,P^{(1,2)}_\Z\,$ is generated by
\begin{equation}
\label{habc6}
\begin{aligned}
s_{1,2} = \ttt{}{\!\!-1}{}{1}{}{}{}{}{\!1},\ \ d_2 = \ttt{1}{}{}{}{\!\!-1}{}{}{}{1},\ \ d_3 = \ttt{1}{}{}{}{1}{}{}{}{\!\!-1} .
\end{aligned}
\end{equation}
Conjugation by either $d_2$ or $d_3$ maps the sets $S=S_{1,p}$ and $S=S_{p,1}$ to sets of the same type, with indices $i,\,j,\,\ell$ permuted modulo $p$. That, as was remarked already, does not alter the right $N_\Z$-cosets generated by $S$. Thus, to verify b) for the group $\,P^{(1,2)}_\Z$, it suffices to show that conjugation by $s_{1,2}$ reshuffles the right $\,P^{(1,2)}_\Z$-cosets generated by $S$. In the case of $S=S_{1,p}$, this comes down to the ${\rm GL}(2)$-analogue of (\ref{habc3}b), which is well known. To deal with $S=S_{p,1}$, note that conjugation by $s_{1,2}$ fixes $h_{p,p,1;0,0,0}$ and interchanges $\{h_{p,1,p;0,j,0}N_\Z \mid j \in \Z/p\Z\}$ with $\{h_{1,p,p;i,0,\ell}N_\Z\mid \ell \in \Z/p\Z,\,i=0\}$. On the other hand, when $i \notin p\Z$, there exists $a,\,b\in \Z$ such that $ap-ib=1$, which gives the matrix identity
\begin{equation}
\label{habc7}
\begin{aligned}
s_{1,2}\, h_{1,p,p;i,0,\ell}\, s_{1,2}^{-1} \, = \, h_{1,p,p;b,p\ell,b\ell}\ttt{a}{b}{}{i}{p}{}{}{}{1} \in  h_{1,p,p;b,p\ell,b\ell}\,P^{(1,2)}_\Z\, .
\end{aligned}
\end{equation}
As $i$ runs over $(\Z/p\Z)^*$, so does $b$; also, for fixed $i,\,b\in (\Z/p\Z)^*$, as $\ell$ runs over $\Z/p\Z$, so does $b\ell$, whereas $bp\equiv 0$ modulo $p$. We conclude that conjugation by $s_{1,2}$ permutes the cosets $\,h_{1,p,p;i,0,\ell}\,P^{(1,2)}_\Z$ indexed by $i\in (\Z/p\Z)^*$ and $\ell\in \Z/p\Z$. In view of what has already been shown, this implies that conjugation by $s_{1,2}$ permutes the $\,P^{(1,2)}_\Z$-cosets generated by $S_{p,1}$. The verification of (\ref{habc4}) is now complete.

For entirely formal reasons, (\ref{habc4}) implies that $T_{1,p}$ and $T_{p,1}$ preserve the space of $P^{(1,2)}_\Z$-invariants in any ${\rm GL}(3,\Q)$-module. In particular
\begin{equation}
\label{habc8}
(T_{1,p} - a_{1,p})\,\tau^{(1,2)}\ \ \text{and}\ \ (T_{p,1} - a_{p,1})\,\tau^{(1,2)}\ \ \text{are $P^{(1,2)}_\Z$-invariant}\,,
\end{equation}
since $\,\tau^{(1,2)}$ itself has this invariance property. We argued earlier that both operators preserve $N_\Z$-isotypic components. In the proof of \lemref{lemconv}, we had constructed $\tau^{(1,2)}$ as the sum of the translates $\pi_{\l,\d}(\g)\tilde\tau^{(1,2)}$, where $\,\tilde\tau^{(1,2)}$ denotes the canonical extension of $\,\tau_{\text{abelian}}$ from $N$ to $N\cup\{x=\infty\}$ and $\g$ runs over a complete set of representatives for $\,P^{(1,2)}_\Z/MN_\Z$. Among these summands, the one corresponding to $\g\in MN_\Z$ lies in the space on which the center of $N_\Z$ acts trivially, whereas all the others lie in spaces on which the center acts
 nontrivially; cf.\ (\ref{tau1}) and (\ref{convergence6}), (\ref{convergence7}). Appealing to
(\ref{a1pp1ab}), we now see that all the abelian Fourier coefficients of $(T_{1,p} - a_{1,p})\,\tau^{(1,2)}$
and $(T_{p,1} - a_{p,1})\,\tau^{(1,2)}$ vanish. The restriction to $N$ of a $P^{(1,2)}_\Z$-invariant section
of $\mathfrak L_{\l,\d}$ over $N\cup \{x=\infty\}$ is completely determined by the abelian Fourier
coefficients; this follows from the proofs of \propref{prop2} and Corollary \ref{cor1}, coupled with the
remark at the end of Section
\ref{heis}. Thus $(T_{1,p} - a_{1,p})\,\tau^{(1,2)}$ and $(T_{p,1} - a_{p,1})\,\tau^{(1,2)}$ vanish at least
on $N$. But the $P^{(1,2)}_\Z$-translates of $N$ cover $N\cup \{x=\infty\}$, forcing $(T_{1,p} -
a_{1,p})\,\tau^{(1,2)}$ and $(T_{p,1} - a_{p,1})\,\tau^{(1,2)}$ to vanish identically. The analogous
assertions about $\tau^{(2,3)}$ can be proved by the same reasoning or, alternatively, can be deduced from
those about $\tau^{(1,2)}$ by applying the outer automorphism (\ref{outer1}).\hfill\qed

\demo{Proof  of  Lemma {\rm \ref{lemcoinc}}}  In view of (\ref{tau1tau2}) and \propref{prop1}, we must
establish the identity (\ref{star2}), as an identity of tempered distributions on $\R$, for all $q\neq 0$ and all
$(a_k)_{k\in\Z/q\Z}$. In fact, we shall establish slightly more. The $\sigma_{q, \,k}$ and $\rho_{q, \,k}$
have been defined not only as tempered distributions, but as vectors in $W_{\l_1-\l_2,\d_3}^{-\infty}$ and
$W_{\l_2-\l_3,\d_1}^{-\infty}$ respectively. By construction, they vanish to infinite order at the origin.
Thus, once (\ref{star2}) has been established for all choices of coefficients $(a_k)$,
\cite[Th.~3.19]{MilSch:2003b} ensures that the $\sigma_{q, \,k}$ and $\rho_{q, \,k}$ can be canonically
extended from $\R$ to $\R\cup\{\infty\}$. It is not obvious, however, that the canonical extensions agree
with how the $\sigma_{q, \,k}$ and $\rho_{q, \,k}$ have been defined -- they might conceivably differ by
distributions supported at $\infty$. We shall show by induction that
\begin{equation}
\label{star2a}
\begin{aligned}
{\sum}_{k=1}^q \ a_k \ \rho_{q, \,k}(y) \ =  \ {\sum}_{\ell=1}^q \   \widehat{a}_\ell \  \widehat{\sigma}_{q, \, \ell}(qy)\ \ \text{as tempered distributions}\,,
\\
\text{and}\ \ {\sum}_{k=1}^q \ a_k \ \rho_{q, \,k}(y) \ \ \text{vanishes to infinite order at $y=\infty\,$}\,,
\end{aligned}
\end{equation}
for all $(a_k)_{k\in\Z/q\Z}$. This implies in particular that both the $\rho_{q, \,k}$ and the $\sigma_{q, \,\ell}$ vanish to infinite order at $\infty$, since the latter are related to the former by the outer automorphism (\ref{outer1}).

In the present setting the parity conditions asserted by \lemref{lem2} hold by construction, so
that (\ref{star2a}) for $q$ and $-q$ are equivalent. We therefore may and shall suppose that $q>0$. Lemma
\ref{melinvconv} already establishes the assertions (\ref{star2a}) whenever $a_k=\chi(k)$ for some
primitive Dirichlet character $\chi$ of conductor $q$. In particular, (\ref{star2a}) is correct when $q=1$.
Arguing inductively, we may suppose that (\ref{star2a}) is satisfied for every smaller choice of $q$.

For the inductive step, we identify the vector space of sequences $(a_k)_{k\in \Z/q\Z}$ indexed by $\Z/q\Z$ with $\operatorname{Fun}(\Z/q\Z,\C)$, the space of $\C$-valued functions on $\Z/q\Z$. If $r>0$ divides $q$, the projection $p(r,q):\Z/q\Z \to \Z/r\Z$ induces the pull-back of functions
\begin{equation}
\label{pullback}
p(r,q)^* :\, \operatorname{Fun}(\Z/r\Z,\C)\ \longrightarrow\ \operatorname{Fun}(\Z/q\Z,\C)\,.
\end{equation}
We shall also use the operation
\begin{equation}
\label{division}
\begin{aligned}
m(r,q) :\, \operatorname{Fun}(\Z/r\Z,\C)& \longrightarrow\ \operatorname{Fun}(\Z/q\Z,\C)\,,
\\
  (m(r,q)\,b)_k& = \ \begin{cases} \, b_{kr\!/q} &\text{if}\ \ \ kr/q \in \Z \\ \,0 &\text{otherwise}\,,\end{cases}
\end{aligned}
\end{equation}
the adjoint of the linear map $\operatorname{Fun}(\Z/q\Z,\C) \ni (a_\ell)\mapsto (a_{\ell q/r})\in\operatorname{Fun}(\Z/r\Z,\C)$ induced by $\Z/r\Z\ni \ell \mapsto \ell q/r \in \Z/q\Z$. The finite Fourier transform (\ref{fourierint3}) relates these two operations: a simple computation shows
\begin{equation}
\label{fourierint4}
a\ = \ p(r,q)^*\,b\ \ \ \Longrightarrow\ \ \ \widehat a \ = \ \f qr \ m(r,q)\,\widehat b
\end{equation}
if $r$ divides $q$, $\,a\in \operatorname{Fun}(\Z/q\Z,\C)\,$, and $\,b\in \operatorname{Fun}(\Z/r\Z,\C)$.

The crux of the inductive argument is the compatibility of (\ref{star2a}) with the two operations 
(\ref{pullback}), (\ref{division}). Specifically:
\begin{equation}
\label{compatibility}
\begin{aligned}
{}&\text{If $\,1\leq r<q\,$ is a divisor of $\,q\,$, the assertions (\ref{star2a}) hold for}
\\
{}&\text{all $\,a\in p(r,q)^*\operatorname{Fun}(\Z/r\Z,\C)\,$ and $\,a\in
m(r,q)\operatorname{Fun}(\Z/r\Z,\C)$}.
\end{aligned}
\end{equation}
Recall that the outer automorphism (\ref{outer1}) switches the roles of the $\sigma_{n,k}$ and the $\rho_{n,k}$. In particular, $a_k$ and $\widehat a_k$ play symmetric roles in (\ref{star2a}). 
Thus,   in view of (\ref{fourierint4}), we only need to verify (\ref{compatibility}) for $a\in
m(r,q)\operatorname{Fun}(\Z/r\Z,\C)$. Also note that $m(r,q)=m(s,q)\circ m(r,s)$ if $r|s$ and $s|q$, so we
only need to treat the case when $q/r$ is prime. Finally, it suffices to consider $a=m(r,q)b$ when $b\in
\operatorname{Fun}(\Z/r\Z,\C)$ is a ``delta function"~-- i.e., $b_k=1$ for exactly one $k\in\Z/r\Z$ and
$b_j=0$ for $j\neq k$ -- since delta functions span $\operatorname{Fun}(\Z/r\Z,\C)$.

To establish this remaining case, we suppose that $q=rp$, for some prime $p$, and $a=m(r,q)b$ with
 $b=$ the delta function at $k\in \Z/r\Z$. Then $\,a$ is the delta function at $pk \in \Z/q\Z$, and $\,\widehat
a_\ell = e(\ell k/r)$. In short, we must equate $\rho_{q,pk}(y)$ with ${\sum}_{\ell\in\Z/q\Z}\, e(\ell
k/r)\,\widehat\sigma_{q,\ell}(qy)$ and show that $\rho_{q,pk}$ vanishes to infinite order at $\infty$. Since
$e(\ell k/r) = e((\ell + rj)k/r)$ for $j\in \Z$, we can re-phrase the assertion to be proved as follows:
\begin{equation}
\label{star3}
\begin{aligned}
\rho_{q,pk}(y)\, = \, {\sum}_{\ell\in\Z/r\Z} \, e(\ell k/r)\, {\sum}_{j\in\Z/p\Z}\ \widehat\sigma_{q,\,\ell+rj}(qy)\,,\ \text{and}
\\
\rho_{q,pk}(y)\ \ \text{vanishes to infinite order at $y=\infty$}\,.
\end{aligned}
\end{equation}
This, finally, is what remains to be shown, for every $k\in\Z/r\Z$.

We shall do so using the inductive hypothesis and two identities which follow from \lemref{lemhecke}:
\begin{equation}
\label{heckerhosigma}
\begin{aligned}
a_{p,1}\, \sigma_{n,k}(x)\   =& \ \ p^{\l_3}\,\sigma_{n/p,\,k\!/p}(x) 
\\
& \ +\ p^{\l_2}\,\sigma_{n,\,pk}(px)\ + \ p^{\l_1-1}\,{\sum}_{j\in\Z/p\Z}\ \sigma_{np,\,k+nj}(x/p)\,,
\end{aligned}
\end{equation}
and
\begin{equation}
\label{heckerhosigma1}
\begin{aligned}
a_{p,1}\,\rho_{n,k}(y) \ =&\ \  p^{\l_3}\,\rho_{n/p,\,k}(p y) 
\\
&\  +\ p^{\l_2-1}\, {\sum}_{j\in\Z/p\Z}\ \rho_{n,\f{k-nj}{p}}(y/p) \ + \ p^{\l_1}\rho_{pn,\,pk}(y)\,.
\end{aligned}
\end{equation}
To deduce these identities from \lemref{lemhecke}, we appeal to the formula (\ref{t1pp1}) for the action of $\,T_{p,1}$ and the fact, which we have used before, that $T_{p,1}$ acts separately on the abelian and non-abelian parts of the Fourier decomposition of $\tau$. Substituting $n=r$ in (\ref{heckerhosigma1}), we find
\begin{equation}
\label{heckerhosigma2}
\begin{aligned}
a_{p,1}\,\rho_{r,k}(y)\   =& \ \ p^{\l_3}\,\rho_{r\!/p,\,k}(p y) 
\\
&\ +\ p^{\l_2-1}\, {\sum}_{j\in\Z/p\Z}\ \rho_{r,\f{k-rj}{p}}(y/p) \ + \ p^{\l_1}\rho_{q,\,pk}(y)\,.
\end{aligned}
\end{equation}
On the other hand,
\begin{equation}
\label{heckerhosigma3}
\begin{aligned}
\sum_{\ell\in\Z/rZ}\ e(k\ell/r)\,a_{p,1}\, \widehat\sigma_{r,\ell}(ry)\ =\ \sum_{\ell\in\Z/rZ}\ e(k\ell/r)
 \biggl(p^{\l_3}\,\widehat\sigma_{r/p,\,\ell/p}(ry) \biggr.
\\
\qquad\qquad\qquad +\ \biggl. p^{\l_2-1}\,\widehat\sigma_{r,\,p\ell}(ry/p)\ + \ p^{\l_1}\,\sum_{j\in\Z/p\Z}\ \widehat\sigma_{q,\,\ell+rj}(qy)\biggr)\,,
\end{aligned}
\end{equation}
which follows from (\ref{heckerhosigma}) when we take Fourier transforms, replace $k$ by $\ell$, $n$ by $r$, $x$ by $ry$, multiply the resulting equation by $e(k\ell/r)$, and take the sum over $\ell\in\Z/r\Z$. As usual, terms with nonintegral indices are to be dropped. By induction,
\begin{equation}
\label{star4}
\begin{aligned}
\rho_{r,k}(y)& \  =  \,\sum_{\ell\in\Z/r\Z} e(\ell k/r)\,\widehat\sigma_{r,\ell}(ry)\,,
\\
\rho_{r\!/p,\,k}(p y)& \ = \, \sum_{\ell\in\Z/\f rp\Z}\! e\!\(\f{\ell k}{r/p}\)\widehat\sigma_{r/p,\,\ell}(ry)\ =\,
\sum_{\ell\in\Z/r\Z}\! e(\ell k/r)\,\widehat\sigma_{r/p,\,\ell/p}(ry)\,;
\end{aligned}
\end{equation}
the second equation is to be disregarded unless $p|r$, of course. Since terms with
 nonintegral subscripts are to be dropped, again by induction,
\allowdisplaybreaks{
\begin{align*}
&\sum_{j\in\Z/p\Z}\ \rho_{r,\f{k-rj}p}(y/p) \ = \ p^{-1} \sum_{i,j\in\Z/p\Z}\ e(i(k-rj)/p)\,\rho_{r,\f{k-rj}p}(y/p)
\\
&\qquad\qquad =\ p^{-1}  \sum_{i,j\in\Z/p\Z}\ \sum_{\ell\in\Z/r\Z}\ e(i(k-rj)/p)\,e(\ell (k-rj)/rp)\,\widehat\sigma_{r,\ell}(ry/p)
\\
&\qquad\qquad = \ p^{-1}  \sum_{i,j\in\Z/p\Z}\ \sum_{0\leq \ell<r}\ e(-j(\ell+ir)/p)\,e(k(\ell+ir)/rp)\,\widehat\sigma_{r,\ell+ir}(ry/p)
\\
&\qquad\qquad = \, \sum_{i\in\Z/p\Z}\ \sum_{\stackrel{\scriptstyle{0\,\leq\, \ell\, < \,r}}{p|(\ell+ir)}}\ e(k(\ell+ir)/rp)\,\widehat\sigma_{r,\ell+ir}(ry/p)\,.
\end{align*}
}
The projection $\Z\to\Z/r\Z$ maps $\{\ell+ir  \mid 0 \leq \ell < r,\,\, 0\leq i <p,\,\,p|(\ell+ir)\}$ onto the set $\{\ell p \in \Z/r\Z \mid \ell\in \Z/r\Z\}$ and preserves multiplicities, hence
\begin{equation}
\label{star5}
{\sum}_{j\in\Z/p\Z}\ \rho_{r,\f{k-rj}p}(y/p) \ =  \  {\sum}_{\ell\in\Z/r\Z}\ e(k\ell/r)\,\widehat\sigma_{r,p\ell}(ry/p)\,.
\end{equation}
The identities (\ref{star4}), (\ref{star5}) match three of the four terms in (\ref{heckerhosigma2}) to the
corres\-ponding terms in (\ref{heckerhosigma3}). Equating what remains, we obtain the expression
(\ref{star3}) for $\rho_{q,pk}$ in terms of the $\widehat\sigma_{q,\ell+rj}$. By induction, we know that
three of the four terms in (\ref{heckerhosigma2}) vanish to infinite order at $y=\infty$, and so
$\rho_{q,pk}$ must vanish to infinite order, too. That establishes (\ref{star3}) and completes the verification
of (\ref{compatibility}).

Lemma \ref{melinvconv} asserts (\ref{star2a}) when $\,a = \chi\,$ for some primitive Dirichlet cha\-racter $\chi$ of conductor $q$. Next we consider the case of an imprimitive Dirichlet character $\chi$ modulo $q$, which is then induced from a primitive character $\chi_r$ whose conductor $r$ properly divides $q$. By definition,
\begin{equation}
\label{impriv}
\chi(\ell)\ \ = \ \
\begin{cases}\ \chi_r(\ell)\ \ \ &\text{if}\ \ \ (\ell,q) = 1
\\
\ 0\ \ \ &\text{if}\ \ \ (\ell,q)\neq 1\ ,
\end{cases}
\end{equation}
so that $\,\chi - p(r,q)^*\chi_r\,$ lies in the subspace
\begin{equation}
\label{units}
\{\,a\in\operatorname{Fun}(\Z/q\Z,\C)\, \mid \,a_\ell = 0\ \ \text{if $\,\ell\,$ is a unit modulo $q$}\,\}\,.
\end{equation}
But (\ref{star2a}) is satisfied for $\,a=p(r,q)^*\chi_r$ because of (\ref{compatibility}). The Dirichlet characters modulo $q$, both primitive and imprimitive, constitute a basis of the space of functions on $(\Z/q\Z)^*$. It therefore suffices to prove (\ref{star2a}) for functions in the space (\ref{units}). Every nonunit $\,\ell\in\Z/q\Z\,$ can be expressed uniquely as the product $\,\ell=qu/r\,$ in terms of a divisor $\, r\,$ of $\,q\,$, $\,1\leq r< q\,$, and a unit $\,u\in (\Z/r\Z)^*$. Hence
\begin{equation}
\label{nonunits}
\begin{aligned}
{}&\{\,a\in\operatorname{Fun}(\Z/q\Z,\C)\, \mid \,a_\ell = 0 \ \text{if $\,\ell\,$ is a unit modulo $q$}\,\} 
\\
{}&\ \ = \bigoplus_{r|q,\,\,1\leq r < q} m(r,q)\{\,a\in\operatorname{Fun}(\Z/r\Z,\C)\, \mid \,a_\ell = 0\ \
\text{if
$\,\ell\notin (\Z/r\Z)^*$}\,\}\,.
\end{aligned}
\end{equation}
In view of (\ref{compatibility}), the assertion (\ref{star2a}) is correct for every $\,a\,$ in 
this space, so (\ref{star2a}) holds in complete  generality.\hfill\qed

\demo{Proof of  Lemma {\rm \ref{leminforder}}}  The outer automorphism (\ref{outer1}) relates the two
statements, and  so it suffices to prove the assertion about $\tau^{(1,2)}$. According to (\ref{convergence1}),
the first term on the right in (\ref{tau1}) vanishes to infinite order along $\{x=\infty\}$. We shall use 
\cite[Lemma 7.2]{MilSch:2003b} to show that also the second term vanishes there to infinite order. It
remains to verify the hypothesis of the lemma in the present setting: we must show that
\begin{equation}
\label{leminforder1}
\begin{aligned}
{}&\psi_{\lambda_1 - \lambda_2,\d_3}\ttwo{1}{-\ell/n}{}{1}\sigma_{n,\ell}\ ,\ \ \ \text{with}\ \
n\in\Z_{\neq 0}\,,\ \ell\in\Z\,,
\\
{}&\text{vanishes to infinite order at $x=\infty$, uniformly in $(n,\ell)$},
\end{aligned}
\end{equation}
in the sense of \cite[Def.~7.1]{MilSch:2003b}.

The distributions $\tau_{y,q}$ are periodic, without constant term.
 As such they can be expressed as $k$-th derivatives of continuous periodic, hence bounded, functions, for
every sufficiently large $k$. Because the Fourier coefficients $c_{n,q}$ have polynomial bounds, it is
possible to choose $k$ independently of $q$, in which case the sup norms of the bounded continuous
functions grow at most polynomially with $q$. The family of distributions $\tau_{y,q}(c^{-2}y-a/c)$,
indexed by $c,\,q\in\Z_{\neq 0}$ and $a\in\Z$, is therefore bounded in the sense of
\cite[(7.14)]{MilSch:2003b}. The members of the family have canonical extensions across $\infty$, which
vanish there to infinite order, uniformly in $(q,c,a)$ in the sense of 
\cite[Def.~7.1]{MilSch:2003b}. According to 
\cite[Lemma
7.12]{MilSch:2003b} and \cite[Lemma 7.16]{MilSch:2003b},
\begin{equation}
\label{leminforder2}
\begin{aligned}
{}&\rho_{n,k}(y)\ = \ (\sg y)^{\d_1}\,|y|^{\l_2-\l_3-1}\,\tau_{y,q}(c^{-2}y^{-1}-a/c)\,,\ \ \text{with}
\\
{}&n=cq,\ k=\bar aq,\ q>0,\ c\neq 0,\ (a,c)=1,\ a\bar a \equiv 1\ (\operatorname{mod}c)\,,
\end{aligned}
\end{equation}
is therefore also bounded as a  family indexed by $(n,k)$ and vanishes to infinite order at $y=0$, uniformly
in
$(n,k)$. Note that (\ref{leminforder2}) is consistent with our earlier definition of the $\rho_{n,k}$.

In the proof of \lemref{lemcoinc}, specifically (\ref{star2a}) and the comment following (\ref{star2a}), we argued that the $\sigma_{n,\ell}$ vanish to infinite order at $\infty$. In this statement we may regard the $\sigma_{n,\ell}$ either as distribution sections of the appropriate line bundle or as scalar distributions on $\,\R\cup\{\infty\}$, since the two interpretations are related by the factor $(\sg y)^{\d_3}|x|^{\l_1-\l_2-1}$. Initially we shall regard them as scalar distributions. Inverting the Fourier transforms in (\ref{star2a}), we find
\begin{equation}
\label{leminforder3}
\sigma_{n,\ell}(x+\ell/n)\ = \ \mathcal F\biggl({\sum}_{k=0}^{n-1}\ e\bigl(\ell(y-k/n)\bigr)\,\rho_{n,k}(y)\biggr)(-nx)
\end{equation}
as an identity of scalar distributions on $\,\R\cup\{\infty\}$, with both sides vanishing to infinte order at $x=\infty$ -- the $\sigma_{n,\ell}$ as was just remarked, and the right-hand side as the Fourier transform of a tempered distribution which vanishes to infinite order at the origin. The boundedness of the family $\rho_{n,k}$ and the fact that they vanish to infinite order at the origin in the uniform sense implies these same two properties for the expression between the large parentheses in (\ref{leminforder3}), viewed as
a family indexed by $(n,\ell)$. The Fourier transforms then vanish to infinite order at $x=\infty$ in the
uniform sense 
\cite[Lemma 7.15]{MilSch:2003b}, and that remains true when the variable $x$ is scaled by the factor $-n$
\cite[Lemma 7.16]{MilSch:2003b}. We have shown that the $\sigma_{n,\ell}(x+\ell/n)$, regarded as a
family of scalar distributions on $\,\R\cup\{\infty\}$, vanish to infinite order at $x=\infty$, uniformly in
$(n,\ell)$. Alternatively we may regard these as a family of vectors in $W_{\lambda_1 - \lambda_2,\d_3}$.
The passage between the two viewpoints involves multiplication by the function $(\sg
x)^{\d_3}|x|^{\l_1-\l_2-1}$, which does not affect the vanishing to infinite order in the uniform sense  
\cite[Lemma 7.12]{MilSch:2003b}. That, in effect, establishes (\ref{leminforder1}) and completes the proof.
\hfill\qed

\demo{Proof  of Lemma {\rm \ref{lemexist}}} We shall construct $\,\tau$ in three stages. To simplify the
discussion, we introduce the notation
\begin{equation}
\label{seto}
\begin{aligned}
{}&X\ = \ G/MAN_-\,,\ \ \ S_x\ = \ \{x=\infty\}\,,\ \ \ S_y\ = \ \{y=\infty\}\,,
\\
{}&\qquad\  \bar S_x\ = \ \text{closure of}\ S_x\,,\ \ \ \ \ \ \ \bar S_y\ = \ \text{closure of}\ S_y\,.
\end{aligned}
\end{equation}
The intersection $\bar S_x\cap\bar S_y$ consists of three Schubert cells, $S_{x,y}$ and $S_{y,x}$, both of dimension one, and $S_0$, which is a single point. The indexing of $S_{x,y}$, $S_{y,x}$ is pinned down by the following statement:
\begin{equation}
\label{Porbits}
\begin{aligned}
{}&P^{(1,2)}\,\ \text{acts on $X$ with orbits}\ \ N\cup S_x\,,\ \ S_y\cup S_{y,x}\,,\ \ S_{x,y}\cup S_0\,,
\\
{}&\text{and}\ \ P^{(2,3)}\,\ \text{acts with orbits}\ \ N\cup S_y\,,\ \ S_x\cup S_{x,y}\,,\ \ S_{y,x}\cup S_0\,.
\end{aligned}
\end{equation}
We shall also use the closure relations
\begin{equation}
\label{orbitclosures1}
\begin{aligned}
{}&\bar S_{x,y}\ = \ S_{x,y} \cup S_0\,,\ \ \ \ \bar S_x \ = \ S_x \cup S_{x,y} \cup S_{y,x} \cup S_0\,,
\\
{}&\bar S_{y,x}\ = \ S_{y,x} \cup S_0\,,\ \ \ \ \bar S_y \ = \ S_y \cup S_{x,y} \cup S_{y,x} \cup S_0\,.
\end{aligned}
\end{equation}
One way to see this is to identify $X$ with the manifold of flags in $\R^3$, i.e.,
 the manifold of nested subspaces $F_1 \subset F_2 \subset \R^3$, with $\operatorname{dim}F_j = j$.
Evidently $G$ acts transitively on this manifold. Let $\{e_1,\,e_2,\,e_3\}$ denote the standard basis of
$\R^3$. The isotropy subgroup at the ``reference flag" $F_1=\R e_3$, $F_2 = \R e_2 \oplus \R e_3$, is the
minimal parabolic subgroup (\ref{minpar}), so that  $X\cong G/MAN_-$, as claimed. The group $N$ acts
freely on its orbit through the reference flag, and can therefore be identified with this orbit:
\begin{equation}
\label{Norbit}
N \ \cong \ \{\,F_1 \not\subset \R e_1 \oplus \R e_2\,,\ F_2\not\supset \R e_1\,\}\,;
\end{equation}
this is the open Schubert cell. Working out the $P^{(1,2)}$- and $P^{(2,3)}$-orbits and their closures, one finds
\begin{equation}
\label{orbitclosures2}
\begin{array}{rlrl}
\bar S_x\hskip-4pt  &=\ \{\,F_2 \supset \R e_1\,\}\,, &\qquad \bar S_y \hskip-4pt &=\ \{\,F_1 \subset \R e_1
\oplus \R e_2
\,\}\,,
\\
\bar S_{y,x}\hskip-3pt &=\ \{\,F_1 = \R e_1 \,\}\,, &\qquad\bar S_{x,y}\hskip-4pt &=\ \{\,F_2 = \R e_1
\oplus \R e_2\,\}\,,
\end{array}
\end{equation}
and
\begin{equation}
S_0\ =\ \bar S_{x,y} \cap \bar S_{y,x}\ = \  \{\,F_1 = \R e_1\,,\ F_2 = \R e_1 \oplus \R e_2 \,\}\,,
\end{equation}
making (\ref{Porbits}), (\ref{orbitclosures1}) completely explicit. Note also that all the orbit closures are
smooth.

Because of \lemref{lemcoinc}, there exists a distribution section $\,\tau_1$ of $\mathcal L_{\lambda,\delta}$ over the open set
\begin{equation}
\label{seto1}
\mathcal O_1 \ = \ N\cup S_x\cup S_y \ = \ X - \bar S_{x,y} - \bar S_{y,x}
\end{equation}
which agrees with $\,\tau^{(1,2)}$ over $N\cup S_x$ and with $\,\tau^{(2,3)}$ over $N\cup S_y$. 
In particular, $\,\tau_1$ is $P^{(1,2)}_\Z$-invariant over $N\cup S_x$ and $P^{(2,3)}_\Z$-invariant over
$N\cup S_y$. According to \lemref{leminforder}, $\,\tau_1$ vanishes to infinite order along $S_y$.
Differently put, the restriction of $\,\tau_1$ to $\mathcal O_1 - S_y$ has a canonical extension across $S_y$,
and $\,\tau_1$ coincides with that canonical extension. Since $\,\tau_1$ is $P^{(1,2)}_\Z$-invariant over
$\mathcal O_1 - S_y=N\cup S_x$, the restriction of $\,\tau_1 = \pi_{\lambda,\delta}(\gamma)\tau_1$ to
$\mathcal O_1 - S_y$ extends canonically also across $\,\gamma S_y$, for every $\gamma\in
P^{(1,2)}_\Z$. The various translates $\gamma S_y$, $\gamma\in P^{(1,2)}_\Z$, cover the
$P^{(1,2)}$-orbit $S_y \cup S_{y,x}$, which is a closed codimension one submanifold of $\mathcal O_1
\cup S_{y,x} = X - \bar S_{x,y}$. By construction, the canonical extension of $\,\tau_1$ from $\mathcal
O_1 - S_y$ to $\mathcal O_1 \cup S_{y,x}$ agrees with $\,\tau_1$ on $\mathcal O_1$ and is
$P^{(1,2)}_\Z$-invariant. We can argue similarly to extend $\,\tau_1$ across $S_{x,y}$. The result is a
distribution section $\,\tau_2$ of $\mathcal L_{\lambda,\delta}$ over
\begin{equation}
\label{seto2}
\mathcal O_2 \ = \ \mathcal O_1 \cup S_{x,y}\cup S_{y,x} \ = \ X - S_0
\end{equation}
which vanishes to infinite order along both $S_x\cup S_{x,y}$ and $S_y\cup S_{y,x}$, and which is $P^{(1,2)}_\Z$-invariant over $\mathcal O_1 \cup S_{y,x} = X-\bar S_{x,y}$, as well as $P^{(2,3)}_\Z$-invariant over $\mathcal O_1 \cup S_{x,y} = X-\bar S_{y,x}$.

The final extension, from $\mathcal O_2 = X - S_0$ to $X$, uses the action of the Weyl group. Together with $MA$ -- recall
(\ref{minpar})!
-- the matrices
\begin{equation}
\label{s12s23}
s_{1,2} \, =  \, \ttt{}{-1}{}{1}{}{}{}{}{\ 1}\in P^{(1,2)}_\Z ,\ \ \  s_{2,3} \, = \, \ttt{1}{}{}{}{}{1}{}{-1}{} \in P^{(2,3)}_\Z
\end{equation}
generate the normalizer $N_G(A)$ of $A$. Conjugation by $s_{1,2}$ and $s_{2,3}$ generates the Weyl group $W(A)=N_G(A)/MA$, and
\begin{equation}
\label{wlong}
w_l \ =  \ s_{1,2}\,s_{2,3}\,s_{1,2}\ = \ s_{2,3}\,s_{1,2}\,s_{2,3}
\end{equation}
represents the longest element of $W(A)$. Because of the known invariance properties of $\,\tau_2$, $\,\pi_{\lambda,\delta}(s_{1,2})\tau_2 = \tau_2 = \pi_{\lambda,\delta}({s_{2,3}}^{-1})\tau_2\,$ on $X - \bar S_{x,y} - \bar S_{y,x}$. We apply $s_{2,3}$ to both sides and note that $s_{2,3}\bar S_{y,x}=\bar S_{y,x}$, to find
\begin{equation}
\label{weyl1}
\pi_{\lambda,\delta}(s_{2,3}\,s_{1,2})\tau_2 \ = \ \tau_2 \ \ \ \text{on}\ \ \ X - s_{2,3}\bar S_{x,y} - \bar S_{y,x}\,.
\end{equation}
Since $\,\pi_{\lambda,\delta}({s_{1,2}}^{-1})\tau_2$ agrees with $\,\tau_2\,$ on $\,X-\bar S_{x,y}\,$, this implies the identity $\pi_{\lambda,\delta}({s_{1,2}}^{-1})\tau_2=\pi_{\lambda,\delta}(s_{2,3}\,s_{1,2})\tau_2$ on the complement of $\,s_{2,3}\bar S_{x,y} \cup \bar S_{y,x}\cup\bar S_{x,y}\,$, hence $\,\pi_{\lambda,\delta}(w_l)\tau_2=\tau_2\,$ on the complement of $\,s_{1,2}s_{2,3}\bar S_{x,y} \cup \,s_{1,2}\bar S_{y,x}\cup\bar S_{x,y}$. We get the same identity on the complement of $\,s_{2,3}s_{1,2}\bar S_{y,x} \cup \,s_{2,3}\bar S_{x,y}\cup\bar S_{y,x}$ after interchanging the roles of $s_{1,2}$ and $s_{2,3}$, and of $S_{x,y}$ and $S_{x,y}$, respectively. We have shown:
\begin{equation}
\label{weyl2}
\begin{aligned}
{}&\pi_{\lambda,\delta}(w_l)\tau_2 \ \ \ \text{and}\ \ \ \tau_2 \ \ \ \text{coincide on the complement of the
set}
\\
{}&(s_{1,2}s_{2,3}\bar S_{x,y} \cup \,s_{1,2}\bar S_{y,x}\cup\bar S_{x,y}) \cap (s_{2,3}s_{1,2}\bar
S_{y,x} \cup \,s_{2,3}\bar S_{x,y}\cup\bar S_{y,x})\,.
\end{aligned}
\end{equation}
The exceptional set can be described more simply as the union $S_0 \cup w_l S_0$, as we shall argue next.

The action of $s_{1,2}$ interchanges the spaces $\R e_1$, $\R e_2$ and fixes $\R e_3$; similarly $s_{2,3}$ fixes $\R e_1$ but interchanges $\R e_2$ and $\R e_3$. With this description of the action one readily identifies the exceptional set in (\ref{weyl2}) with $S_0 \cup w_lS_0$, which is also the complement of $\mathcal O_2 \cap w_l\mathcal O_2$. Thus $\pi_{\lambda,\delta}(w_l)\tau_2 = \tau_2$ on $\mathcal O_2 \cap w_l\mathcal O_2$. We can therefore ``glue together" $\pi_{\lambda,\delta}(w_l)\tau_2$ and $\tau_2$ to a well-defined distribution section $\tau$ of $\mathcal L_{\lambda,\delta}$ over $\mathcal O_2 \cup w_l\mathcal O_2 = X$. In other words, there exists $\tau\in V_{\lambda,\delta}^{-\infty}$ which is $w_l$-invariant and agrees with $\tau_2$ over $X-S_0$.

On $X-S_0$, $\tau$ agrees  with $\tau_2$, which is $P^{(1,2)}_\Z$-invariant on $X-\bar S_{x,y}$ and $P^{(2,3)}_\Z$-invariant on $X-\bar S_{y,x}$. Since $\bar S_{x,y}\cap\bar S_{y,x}=S_0$ and $P^{(1,2)}_\Z\cap P^{(2,3)}_\Z = M N_\Z$,
\begin{equation}
\label{weyl4}
\tau\ \ \text{is $MN_\Z$-invariant on}\,\ X - S_0\,.
\end{equation}
But $w_l$ normalizes $M$ and $w_l S_0 \cap S_0 = \emptyset$, so $\tau$ is globally $M$-invariant. Since $\tau=\pi_{\lambda,\delta}(w_l)\tau$ is $P^{(1,2)}_\Z$-invariant on the complement of $\bar S_{x,y}$, it is invariant under $w_lP^{(1,2)}_\Z w_l^{-1}$ on the complement of $w_l\bar S_{x,y}$. Since $w_l S_0\subset N$ and $\bar S_{x,y}$ are disjoint, so are $w_l\bar S_{x,y}$ and $S_0$. At this point (\ref{weyl4}) allows us to conclude the global invariance of $\tau$ under $N_\Z \cap w_lP^{(1,2)}_\Z w_l^{-1}$, and also under $N_\Z \cap w_lP^{(2,3)}_\Z w_l^{-1}$, as can be seen by interchanging the roles of the subscripts $x$ and $y$. Since $w_l$, $M$, $N_\Z \cap w_lP^{(1,2)}_\Z w_l^{-1}$ and $N_\Z \cap w_lP^{(2,3)}_\Z w_l^{-1}$ collectively generate $\G$, we have shown that $\tau$ is $\G$-invariant.

To establish the cuspidality of $\tau$, we use the notation of (\ref{cuspidality1}) to define two projection operators
\begin{equation}
\label{cuspidality3}
\tau\ \mapsto \ P_{x,z}\tau\ =_{\rm def}\ \tau_{x,z}\ \ \ \ \text{and}\ \ \ \ \tau\ \mapsto \ P_{y,z}\tau\ =_{\rm def}\ \tau_{y,z}\,,
\end{equation}
by taking the average of, respectively, the $N_{x,z}$- and $N_{y,z}$-translates of $\tau$ over a fundamental domain for the action of the subgroups of integral points. Since the abelian Fourier coefficients $c_{m,0}$ and $c_{0,n}$ vanish by construction, the proof of \lemref{lem1} shows that $P_{x,z}\tau=0=P_{y,z}\tau$ at least on the open Schubert cell $N$. In the course of the construction of $\tau$ we saw that $\tau$ vanishes to infinite order along $S_x\cup S_{x,y}$ and $S_y\cup S_{y,x}$. These are $N$-invariant, closed submanifolds of $X-\bar S_{y,x}$, respectively $X-\bar S_{x,y}$, which are $N$-invariant and open in $X$. We can therefore appeal to
\cite[Prop.~7.20]{MilSch:2003b} and conclude that both $P_{x,z}\tau$ and $P_{y,z}\tau$ vanish to infinite
order along $S_x\cup S_{x,y}$ and along $S_y\cup S_{y,x}$. We already know that these projections
vanish on $N$, which is the complement of $S_x\cup S_{x,y}\cup S_y\cup S_{y,x}$ in $(X-\bar
S_{x,y})\cup (X-\bar S_{x,y})= X - S_0$. Thus $P_{x,z}\tau$ and $P_{y,z}\tau$ must vanish on $X-S_0$ 
\cite[Lemma 2.8]{MilSch:2003b}.

Note that $S_y\cup S_{y,x}= \bar S_y \cap(X-\bar S_{x,y})$. As was just remarked, $\tau$ vanishes to infinite order along this closed submanifold of $X-\bar S_{x,y}$. Thus $\tau=\pi_{\lambda,\delta}(s_{2,3})\tau$ vanishes to infinite order along the $s_{2,3}N s_{2,3}^{-1}$-invariant closed submanifold $s_{2,3}\bar S_y \cap(X-s_{2,3}\bar S_{x,y})$ of the open, $s_{2,3}N s_{2,3}^{-1}$-invariant subset $X-s_{2,3}\bar S_{x,y}\subset X$. Since $N_{x,z} \subset s_{2,3}N s_{2,3}^{-1}$, we can appeal once more to 
\cite[Prop.~7.20]{MilSch:2003b} and conclude that $P_{x,z}\tau$ vanishes to infinite order along
$s_{2,3}\bar S_y \cap(X-s_{2,3}\bar S_{x,y})$. But $P_{x,z}\tau=0$ on $X-S_0$ and
\begin{equation}
\label{cuspidality4}
s_{2,3}\bar S_y \cap(X-s_{2,3}\bar S_{x,y})\ \ = \ \ \{\,F_1 \subset \R e_1 \oplus \R e_3\,,\ F_2\neq \R e_1 \oplus \R e_3\,\}
\end{equation}
contains $S_0=\{F_1= \R e_1,\,F_2= \R e_1 \oplus \R e_2\}$; hence $P_{x,z}\tau=0$ on all of $X$. In this
argument we can switch the roles of the subscripts $x$, $y$, and interchange $s_{1,2}$ and $s_{2,3}$,
allowing us to conclude 
 also that $P_{y,z}\tau$ vanishes globally. This makes $\tau$ cuspidal.
\hfill\qed

\end{document}